\documentclass[11pt,reqno]{amsart}
\usepackage[margin=25mm]{geometry}
\usepackage[T1]{fontenc}
\usepackage{lmodern}\usepackage[protrusion=true,expansion=false]{microtype}
\usepackage{amsmath,amssymb,amsthm,mathtools,mathrsfs}
\usepackage{booktabs,tabularx,longtable,array,enumitem}
\usepackage[unicode,hidelinks]{hyperref}
\hypersetup{
  pdfborder={0 0 0},
  colorlinks=false,
  pdftitle={Nonsymmetric Newton Interpolation and Relative Habiro Lifts of Type-A Knot Invariants},
  pdfauthor={Honghuai Fang and Tian Zhou},
  pdfsubject={Integral Newton interpolation and relative Habiro lifts for type-A knot invariants},
  pdfkeywords={Newton interpolation, Habiro ring, knot invariants, Alexander polynomial}
}
\setlist{itemsep=0.25em,topsep=0.4em}
\allowdisplaybreaks[2]
\numberwithin{equation}{section}
\theoremstyle{plain}
\newtheorem{theorem}{Theorem}[section]
\newtheorem{proposition}[theorem]{Proposition}
\newtheorem{lemma}[theorem]{Lemma}
\newtheorem{corollary}[theorem]{Corollary}
\theoremstyle{definition}
\newtheorem{definition}[theorem]{Definition}
\newtheorem{example}[theorem]{Example}
\theoremstyle{remark}
\newtheorem{remark}[theorem]{Remark}
\theoremstyle{plain}
\newtheorem{maintheorem}{Theorem}

\newcommand{\A}{\mathbb Z[Q^{\pm1}]}
\newcommand{\Z}{\mathbb Z}
\newcommand{\Q}{\mathbb Q}
\newcommand{\N}{\mathbb N}
\newcommand{\F}{\mathbb F}
\newcommand{\Habi}{\mathcal H}
\newcommand{\Frac}{\operatorname{Frac}}
\newcommand{\Alt}{\operatorname{Alt}}
\newcommand{\sgn}{\operatorname{sgn}}
\newcommand{\hgt}{\operatorname{ht}}
\newcommand{\Tr}{\operatorname{Tr}}
\newcommand{\ord}{\operatorname{ord}}

\title[Newton Interpolation and Relative Habiro Lifts]{Nonsymmetric Newton Interpolation and Relative Habiro Lifts of Type-$A$ Knot Invariants}
\author{Honghuai Fang}
\address{Institute for Theoretical Sciences, Westlake University, No. 600 Dunyu Road, Xihu District, Hangzhou, Zhejiang 310030, China}
\email{fanghonghuai@westlake.edu.cn}
\author{Tian Zhou}
\address{Peking University, No. 5 Yiheyuan Road, Haidian District, Beijing 100871, China}
\email{2201110034@pku.edu.cn}
\begin{document}
\begin{abstract}
The reduced $\mathfrak{sl}_n$ invariants of a knot form a family of
Laurent polynomials $J_{K,\lambda}(Q)$ indexed by the dominant weights
$\lambda=\sum_ir_i\omega_i$. We interpolate this family in independent
color variables $t_i=Q^{r_i}$ by Newton series with coefficients in
$\mathbb Z[Q^{\pm1}]$. An integral interpolant exists exactly when all
mixed divided differences of the colored values are integral, and it is
determined by its joint Taylor expansion at $Q=1$, $t_i=1$. Inside the
resulting Newton ring we describe the image of a relative Habiro ring
with denominator $D$: a Newton series lies in this image if and only
if, at every root of unity $\zeta$ of order $m$, the coefficients of
its expansion in $Q-\zeta$ are rational functions of the colors whose
denominators are powers of $D(t_1^m,\ldots,t_d^m)$. For knots, the
integral interpolant comes from the cyclotomic expansion of Beliakova
and Gorsky, and rationality is proved by computing traces of braid
operators on Verma modules at roots of unity and, at $Q=1$, by the
rational Kontsevich invariant. As a result, for every $n\ge2$ the
reduced $\mathfrak{sl}_n$ invariants in all dominant weights have a
unique relative Habiro lift. Its denominator is the product of
$\Delta_K(t_a t_{a+1}\cdots t_{b-1})$ over the positive roots
$\alpha_a+\cdots+\alpha_{b-1}$, and at $Q=1$ its leading term is the
reciprocal of this product.
\end{abstract}
\maketitle

\section{Introduction}
\label{sec:introduction}

Let $K$ be a knot in $S^3$ and write a dominant weight of
$\mathfrak{sl}_n$ as
\[
 \lambda=r_1\omega_1+\cdots+r_{n-1}\omega_{n-1},
 \qquad \mathbf r\in\mathbb N^{n-1},
\]
where the $\omega_i$ are the fundamental weights and
$\mathbb N=\mathbb Z_{\ge0}$. We write $J_{K,\lambda}^{(n)}(Q)$ for the
reduced, zero-framed quantum invariant of $K$ colored by the
irreducible representation of highest weight $\lambda$, normalized to
be one on the unknot. With $Q=q^2$ it lies in $A=\mathbb Z[Q^{\pm1}]$
\cite{RT,Lusztig,LeIntegrality}.

This paper treats the whole family $\{J^{(n)}_{K,\lambda}\}_\lambda$ as
a single function of the color. We introduce color variables
\[
 t_i=Q^{r_i},\qquad 1\le i\le n-1,
\]
and look for a function $F(\mathbf t;Q)$, with $Q$ and the $t_i$
independent, whose value at $t_i=Q^{r_i}$ is the invariant colored by
$\sum_ir_i\omega_i$. For $n\ge3$ the colors include mixed weights such
as $r\omega_1+s\omega_2$, and not only the symmetric powers $r\omega_1$,
whose Newton interpolation was studied in \cite{FangZhou,FangNewton}.

Two questions organize the paper. The first is integrality: which
families of Laurent polynomials admit an interpolating series with
coefficients in $A$, and does the family of knot invariants admit one?
The second concerns roots of unity. As $Q$ approaches a root of unity,
the interpolation nodes $Q^r$ collide, and an integral interpolating
series acquires a Taylor expansion in the color variables that is not
visible in any single colored invariant. The figure-eight knot in rank
one shows this already. Each of its colored Jones polynomials equals
one at $Q=1$, whereas the interpolating function, with $t$ kept free,
specializes at $Q=1$ to
\[
 \frac{1}{3-t-t^{-1}}=\frac{1}{\Delta_{4_1}(t)}.
\]
The Taylor coefficients of this function at $t=1$ are lost if one
first substitutes $t=Q^r$ and then sets $Q=1$. At a root of unity of
order $m$ the Alexander polynomial reappears, evaluated at $t^m$.

Our main result, Theorem~\ref{main:type-a}, shows that this is the
general pattern in type $A$. For every knot and every $n\ge2$ the
reduced $\mathfrak{sl}_n$ invariants in all dominant weights are the
color values of a unique element of a relative Habiro ring. At each
root of unity its local expansion is a power series in $Q-\zeta$ whose
coefficients are rational functions of the colors, with poles only
along the Alexander divisors attached to the positive roots. The proof
has two independent parts. The first is a recognition criterion,
Theorem~\ref{main:recognition}, valid for arbitrary integral Newton
series: such a series comes from a relative Habiro ring as soon as its
local expansions are rational over the cyclotomic field. The second
part applies this criterion to knots. The integral interpolant is
constructed by factorial Schur interpolation
(Theorem~\ref{main:symmetric}), and the rationality of its expansions
is proved by a trace computation for quantum braid operators and, at
$Q=1$, by the rational Kontsevich invariant.

\subsection{Integral interpolation of colored values}
\label{subsec:intro-integral-interpolation}

In one variable we use the Newton polynomials for the nodes
$1,Q,Q^2,\ldots$,
\[
 P_k(t;Q)=\prod_{j=0}^{k-1}(t-Q^j),\qquad P_0(t;Q)=1.
\]
Since $P_k(Q^r;Q)=0$ for $k>r$, a series
\begin{equation}\label{eq:intro-one-variable-series}
 F(t;Q)=\sum_{k\ge0}c_k(Q)P_k(t;Q),\qquad c_k(Q)\in A,
\end{equation}
has a Laurent polynomial value at every color:
\[
 F(Q^r;Q)=\sum_{k=0}^r c_k(Q)P_k(Q^r;Q).
\]
The $k$th term changes the value at the node $Q^k$ and leaves the
values at $1,Q,\ldots,Q^{k-1}$ unchanged. Integrality of the
coefficients is therefore a divisibility condition on the colored
values. The first instance is
\begin{equation}\label{eq:intro-first-difference}
 F(Q;Q)-F(1;Q)=(Q-1)c_1(Q),
\end{equation}
and in general $c_k$ is the $k$th Newton divided difference of the
values \cite{deBoor,HarmanHopkins}. Divisibility of all pairwise
differences is not enough. The values $f(0)=f(1)=0$ and $f(2)=Q^2-1$
at the nodes $1,Q,Q^2$ satisfy $Q^a-Q^b\mid f(a)-f(b)$ for all $a,b$,
but their second Newton coefficient is $1/(Q(Q-1))$.

For $d$ independent colors we use the products
\[
 P_{\boldsymbol\alpha}(\mathbf t;Q)
   =\prod_{i=1}^dP_{\alpha_i}(t_i;Q),
 \qquad \boldsymbol\alpha\in\mathbb N^d,
\]
and define the dominant Newton ring
\begin{equation}\label{eq:intro-dominant-ring}
 \mathcal N_d^+(A)
   =\varprojlim_{\mathbf m}
    A[t_1^{\pm1},\ldots,t_d^{\pm1}]
      /(P_{m_1}(t_1),\ldots,P_{m_d}(t_d)),
\end{equation}
where the positive integers $m_i$ increase independently. At each
finite level, division by the monic polynomials $P_{m_i}$ shows that
the $P_{\boldsymbol\alpha}$ with $\alpha_i<m_i$ form a basis over $A$,
and the transition maps forget the coefficients outside the smaller
box. Every element therefore has a unique expansion
\begin{equation}\label{eq:intro-tensor-newton}
 F(\mathbf t;Q)=\sum_{\boldsymbol\alpha\in\mathbb N^d}
       c_{\boldsymbol\alpha}(Q)P_{\boldsymbol\alpha}(\mathbf t;Q),
 \qquad c_{\boldsymbol\alpha}(Q)\in A,
\end{equation}
and products are computed in the finite quotients. At the color
$\mathbf r$ only the terms with $\alpha_i\le r_i$ for all $i$ are
nonzero.

Coefficients with several nonzero indices couple different color
directions. For $f(r,s)=F(Q^r,Q^s;Q)$ the first of them is
\begin{equation}\label{eq:intro-mixed-divided-difference}
 c_{1,1}(Q)=
 \frac{f(1,1)-f(1,0)-f(0,1)+f(0,0)}{(Q-1)^2}.
\end{equation}
The numerator is what remains of the value at $(1,1)$ after the
changes along the two axes have been subtracted, and an integral
interpolant requires it to be divisible by $(Q-1)^2$. Iterated divided
differences give the corresponding condition for every mixed index.

\begin{maintheorem}[Integral interpolation and Taylor rigidity]
\label{main:newton}
Evaluation identifies $\mathcal N_d^+(A)$ with the functions
$\mathbb N^d\to A$ that are represented by a Laurent polynomial on
every finite set of colors.  Its joint Taylor map at
$(Q,\mathbf t)=(1,\mathbf1)$ is injective.  In particular,
$\mathcal N_d^+(A)$ is a domain.
\end{maintheorem}

In the first statement the Laurent polynomial may depend on the finite
set; what is required is that one polynomial represents all values on
that set. The proof is short. On a box of colors the evaluation matrix
in the Newton basis is triangular with nonzero diagonal, so the
values determine the Newton coefficients, and the coefficients
computed on nested boxes agree. The second statement concerns the
joint expansion in $Q-1$ and $t_i-1$, which we describe next.

\subsection{Taylor expansions at colliding nodes}
\label{subsec:intro-taylor-reconstruction}

At $Q=1$ all interpolation nodes equal $1$. Put $h=Q-1$ and
$u_i=t_i-1$. Every factor $t_i-Q^j=u_i-((1+h)^j-1)$ of
$P_{\boldsymbol\alpha}$ lies in the ideal $(h,u_1,\ldots,u_d)$, so the
term of index $\boldsymbol\alpha$ has order at least
$|\boldsymbol\alpha|=\sum_i\alpha_i$. The Newton sum therefore
converges to a formal power series
\begin{equation}\label{eq:intro-joint-taylor}
 F(\mathbf1+\mathbf u;1+h)
       \in\mathbb Z[[h,u_1,\ldots,u_d]],
\end{equation}
the joint Taylor germ of $F$ at the identity.

The substitution $u_i=(1+h)^{r_i}-1$ turns this germ into the Taylor
series at $Q=1$ of the colored value $F(Q^{\mathbf r};Q)$. If the germ
is zero, all colored values have zero Taylor series, so they vanish,
and $F=0$ by the first part of Theorem~\ref{main:newton}. This proves
the rigidity statement.

Keeping the colors free while setting $h=0$ has a different effect.
The Newton basis becomes the monomial basis, and
\[
 F(\mathbf1+\mathbf u;1)
   =\sum_{\boldsymbol\alpha}c_{\boldsymbol\alpha}(1)
          \mathbf u^{\boldsymbol\alpha}.
\]
Every formal power series in $\mathbf u$ with integer coefficients
arises in this way. For example, the element $t-1$ specializes to $u$,
although each of its colored values $Q^r-1$ vanishes at $Q=1$. More
precisely, for every $M\ge1$,
\begin{equation}\label{eq:intro-quantum-neighborhood}
 \mathcal N_d^+(A)/(Q-1)^M\mathcal N_d^+(A)
   \simeq (\mathbb Z[h]/(h^M))[[u_1,\ldots,u_d]]
\end{equation}
(Proposition~\ref{prop:newton-quantum-jets}). Hence the first $M$ quantum
coefficients of a Newton series may be arbitrary formal series in the
colors. Whether they are rational functions of the colors is an
additional property, and it is this property that singles out the
relative Habiro ring inside the Newton ring.

At a primitive $m$th root of unity $\zeta$ the nodes repeat with
period $m$, and the Newton basis factors as
\begin{equation}\label{eq:intro-period-factorization}
 P_{mk+b}(t;\zeta)=(t^m-1)^kP_b(t;\zeta),
 \qquad 0\le b<m.
\end{equation}
Each complete period contributes a factor $t^m-1$, which vanishes at
every $m$th root of unity. For this reason the Newton sum also has a
joint Taylor germ at each center $(Q,\mathbf t)=(\zeta,\zeta^{\mathbf b})$.

These germs are determined, to any finite order, by finitely many
colors. The total-order $L$ jet at such a center consists of the terms
of total degree at most $L$ in $Q-\zeta$ and the $t_i-\zeta^{b_i}$.
Theorem~\ref{thm:newton-finite-jets} shows that this jet is determined
by the Taylor expansions through order $L$ at $Q=\zeta$ of the colored
values in the box
\[
 0\le r_i\le m(L+1)-1,\qquad 1\le i\le d.
\]
The vanishing order of $P_{\boldsymbol\alpha}$ at the center grows
linearly with $\boldsymbol\alpha$, so only finitely many Newton
coefficients enter the jet. Each of them is a divided difference of
colored values, and the division costs at most as many orders of
$Q-\zeta$ as the vanishing order of the corresponding Newton
polynomial; this is why colored values through order $L$ suffice.

\subsection{Singular colors and Habiro boundaries}
\label{subsec:intro-singular-boundary}

The first node outside the dominant range is $t=Q^{-1}$, that is,
$r=-1$. There the Newton basis takes the values
\begin{equation}\label{eq:intro-wall-factorial}
 P_k(Q^{-1};Q)=Q^{-k}(Q;Q)_k,
 \qquad (Q;Q)_k=\prod_{j=1}^k(1-Q^j).
\end{equation}
These are the products that define Habiro's cyclotomic completion
\[
 \mathcal H=\varprojlim_k A/((Q;Q)_k)
\]
\cite{Habiro,Habiro2002}. Consequently every integral Newton series
has a value at the missing node,
\[
 F(Q^{-1};Q)=\sum_{k\ge0}c_k(Q)Q^{-k}(Q;Q)_k\in\mathcal H,
\]
although this value need not be a Laurent polynomial. An element of
$\mathcal H$ has a Taylor expansion at every root of unity, because
$(Q;Q)_k$ vanishes to order $\lfloor k/m\rfloor$ at a root of order
$m$. In particular, the value of the displayed sum at a root of order
$m$ only involves the terms with $k<m$.

In several variables, restriction to $t_1=Q^{-1}$ gives a Newton
series in the remaining colors with coefficients in $\mathcal H$, and
the same holds for all normal derivatives along this wall. Write
$u=t_1-Q^{-1}$, and let $\mathcal N_j(B)$ denote the Newton completion
in $j$ variables with coefficients in an $A$-algebra $B$.
Theorem~\ref{thm:newton-wall-jets} states that for every $M\ge1$
\begin{equation}\label{eq:intro-wall-jets}
 \mathcal N_d^+(A)/u^M\mathcal N_d^+(A)
    \simeq\mathcal N_{d-1}(\mathcal H)[u]/(u^M).
\end{equation}
Every finite list of Habiro-valued normal coefficients is therefore
realized by an integral Newton series, and the complete list of normal
coefficients determines the series.

The boundary value also decides whether $F$ extends one step beyond
the dominant chamber. By Theorem~\ref{thm:newton-wall-fiber-product}, a
Newton series extends to an integral interpolant on
$\{-1,0,1,\ldots\}\times\mathbb N^{d-1}$ exactly when its boundary
value lies in $\mathcal N_{d-1}(A)$. In one variable this asks whether
the Habiro element $F(Q^{-1};Q)$ is a Laurent polynomial. The knot
elements constructed below have such boundary values along each simple
wall.

\subsection{Factorial Schur interpolation and the knot element}
\label{subsec:intro-factorial-schur}

An integral knot element in independent colors requires integrality of
all mixed divided differences of the colored invariants. We obtain it
from factorial Schur interpolation and from the integral elements of
the completed quantum group center constructed by Beliakova and Gorsky
\cite{BG}.

Put $R_n=A[x_1,\ldots,x_n]$ and $\Lambda_n=R_n^{S_n}$. A partition
$\mu=(\mu_1\ge\cdots\ge\mu_n\ge0)$ determines the node
\[
 p_\mu=(Q^{\mu_1+n-1},Q^{\mu_2+n-2},\ldots,Q^{\mu_n}).
\]
The factorial Schur polynomials
\[
 S_\nu(\mathbf x;Q)=
 \frac{\det[P_{\nu_i+n-i}(x_j;Q)]_{i,j=1}^n}
      {\prod_{i<j}(x_i-x_j)}
\]
form an $A$-basis of $\Lambda_n$ and have the interpolation property
\[
 S_\nu(p_\mu;Q)=0\quad\text{unless}\quad
 \nu_i\le\mu_i\ \text{for every }i
\]
\cite{MolevSagan}. Hence an infinite sum $\sum_\nu b_\nu S_\nu$ has a
finite value at each partition node. We write $\mathcal B_n$ for the
completion of $\Lambda_n$ formed by such sums. The integral knot
element of \cite{BG} lies in $\mathcal B_n$, and its values at the
partition nodes are the reduced invariants.

The exponents of a partition node are strictly decreasing, and
permuting them gives exactly the ordered regular grid
\[
 (Q^{a_1},\ldots,Q^{a_n}),\qquad
 a_i\in\mathbb N,\quad a_i\ne a_j\quad(i\ne j).
\]
Let $\mathcal N_n^{\mathrm{reg},+}$ be the interpolation completion of
$R_n$ on this grid. Unlike $\mathcal B_n$, it keeps track of the order
of the coordinates, and it is a finite free extension of the symmetric
completion.

\begin{maintheorem}[Symmetric interpolation and dominant colors]
\label{main:symmetric}
There are natural identifications
\[
 (\mathcal N_n^{\mathrm{reg},+})^{S_n}=\mathcal B_n,
 \qquad
 \mathcal N_n^{\mathrm{reg},+}
       \simeq R_n\otimes_{\Lambda_n}\mathcal B_n.
\]
This extension is finite free of rank $n!$ over $\mathcal B_n$.
Under the Dynkin substitution
\[
 x_i=Q^{n-i}\prod_{j=i}^{n-1}t_j\quad(i<n),\qquad x_n=1,
\]
the integral Beliakova--Gorsky knot element gives an element
$\mathcal J_{K,n}\in\mathcal N_{n-1}^+(A)$ whose values are the
reduced invariants of all dominant $\mathfrak{sl}_n$ colors.
\end{maintheorem}

A free basis of the extension is given by the monomials
\[
 x_1^{a_1}\cdots x_n^{a_n},\qquad 0\le a_i\le n-i.
\]
The point of the theorem is that all of this holds over $A$, not only
over the field $\mathbb Q(Q)$. Integral divided-difference operators
lift invariant classes to symmetric polynomials, and the successive
monic equations satisfied by ordered roots give the free basis
\cite{Demazure,EkedahlLaksov}.

The Dynkin substitution sends $\mathbf t=Q^{\mathbf r}$ to $p_\mu$
with $\mu_i=r_i+\cdots+r_{n-1}$ and $\mu_n=0$. In $\mathfrak{sl}_3$,
for example, it reads
\[
 (x_1,x_2,x_3)=(Q^2t_1t_2,Qt_2,1)
              =(Q^{r_1+r_2+2},Q^{r_2+1},1),
\]
which corresponds to the partition $(r_1+r_2,r_2,0)$ and the weight
$r_1\omega_1+r_2\omega_2$. On a finite box of Dynkin labels the
partition sizes are bounded, so by the vanishing property only
finitely many $S_\nu$ are nonzero there. This makes the substitution
continuous and produces the integral Newton element $\mathcal J_{K,n}$.

Its first mixed Newton coefficient gives the congruence
\begin{equation}\label{eq:intro-knot-mixed-congruence}
 J_{K,\omega_i+\omega_j}^{(n)}(Q)
 -J_{K,\omega_i}^{(n)}(Q)-J_{K,\omega_j}^{(n)}(Q)+1
       \in(Q-1)^2A\qquad(i\ne j),
\end{equation}
and all higher mixed divided differences of the colored invariants are
integral as well. The knot element also extends to integer exponent
tuples with distinct entries and is invariant under permutations of
the shifted exponents; in Dynkin coordinates this is the shifted Weyl
symmetry (Corollary~\ref{cor:sym-shifted-weyl-knot}). On the
symmetric-power axis it agrees with the construction of
\cite{FangNewton}. Section~\ref{subsec:sym-comparison}
identifies the two functions and relates their singular boundaries.

\subsection{Rational color expansions and the relative Habiro ring}
\label{subsec:intro-relative-image}

By Theorem~\ref{main:newton}, an element of the Newton ring is the
same thing as a family of colored Laurent polynomials with integral
mixed divided differences. The next question is which of these
families have local expansions at roots of unity that are rational in
the colors. Fix
$D\in\mathbb Z[T_1^{\pm1},\ldots,T_d^{\pm1}]$ with $D(\mathbf1)=1$. At
a primitive $m$th root $\zeta$ put $x=Q-\zeta$ and $T_i=t_i^m$, and
consider the local series
\begin{equation}\label{eq:intro-relative-local-ring}
 g_\zeta(\mathbf t,x)=\sum_{h\ge0}g_{\zeta,h}(\mathbf t)x^h,
 \qquad
 g_{\zeta,h}\in
 \mathbb Z[\zeta][\mathbf t^{\pm1},D(\mathbf t^m)^{-1}].
\end{equation}
Each coefficient is a Laurent polynomial divided by a power of
$D(\mathbf t^m)$, and the power may grow with $h$. The normalization
$D(\mathbf1)=1$ guarantees a Taylor expansion at $\mathbf t=\mathbf1$.
When $D=1$ all quantum coefficients are Laurent polynomials in the
colors.

The components at different roots must be compatible, in the sense of
the relative gluing condition of Garoufalidis and Wheeler \cite{GW}.
First, Galois conjugation acts on the root and on the coefficients
together. Second, let $p$ be a prime, and let $\zeta$ and $\xi$ have
orders $m$ and $pm$ and be congruent modulo a prime above $p$. The
series at $\zeta$ is re-expanded at $\xi$ by
\[
 Q-\zeta=(Q-\xi)+(\xi-\zeta),
\]
and with the $t_i$ held fixed the result must equal the component at
$\xi$. In terms of $T_i=t_i^m$, passing from $\zeta$ to $\xi$ replaces
$T_i$ by $T_i^p$; this is the Frobenius twist in the gluing condition
of \cite{GW}. The comparison takes place in a $p$-adically completed
coefficient ring, where the congruence
\[
 D(\mathbf T^p)\equiv D(\mathbf T)^p\pmod p
\]
makes the second denominator invertible. For $d=0$ these conditions
reduce to the compatibility of the Taylor expansions of an element of
$\mathcal H$ \cite{Habiro,GSWZ}.

Write $\mathcal H_{\mathrm{rel}}(D)$ for the ring of such families.
Substituting $t_i=Q^{r_i}$ at a dominant color gives an element of
$\mathcal H$, and the families all of whose color values lie in $A$
form a subring $\mathcal H_{\mathrm{rel}}(D)_{\mathrm{pol}}$.

\begin{maintheorem}[Recognition of relative Habiro lifts]
\label{main:recognition}
Color evaluation embeds $\mathcal H_{\mathrm{rel}}(D)$ in
$\mathcal N_d(\mathcal H)$ and
$\mathcal H_{\mathrm{rel}}(D)_{\mathrm{pol}}$ in
$\mathcal N_d^+(A)$. A Newton series $F\in\mathcal N_d^+(A)$ lies in
the image of $\mathcal H_{\mathrm{rel}}(D)_{\mathrm{pol}}$ if and only
if, at every primitive $m$th root of unity $\zeta$, its joint Taylor
germ at $\mathbf t=\mathbf1$ is the germ of a series in
\[
 \mathbb Q(\zeta)[\mathbf t^{\pm1},D(\mathbf t^m)^{-1}]
       [[Q-\zeta]].
\]
In that case the relative preimage is unique.
\end{maintheorem}

Thus, for an integral Newton series, rationality over the cyclotomic
field is all one has to check: integrality of the numerators and the
gluing between different roots then follow. This is what makes the
criterion usable for trace computations, which produce rational
functions over a field.

The proof rests on two properties of Newton germs. First, their Taylor
coefficients are integral. If a rational representative is multiplied
by a power of $D(\mathbf t^m)$, the result is a Laurent polynomial
with integral Taylor coefficients at $\mathbf t=\mathbf1$, hence with
integral coefficients. Second, a Newton series can be re-expanded at
$p$-adically close roots, and its germs agree after re-expansion.
Since color Taylor expansion is injective on the completed
localization, this agreement passes to the rational representatives.
Galois equivariance follows from the same uniqueness.

The criterion can also be read off from the Newton coefficients.
Group the indices as $\boldsymbol\alpha=m\mathbf k+\mathbf b$ with
$0\le b_i<m$, and form the generating series
\[
 C_{\zeta,\mathbf b,h}(\mathbf z)
   =\sum_{\mathbf k\in\mathbb N^d}
      [x^h]c_{m\mathbf k+\mathbf b}(\zeta+x)
      \mathbf z^{\mathbf k}.
\]
In quantum order zero, the period factorization
\eqref{eq:intro-period-factorization}, with $z_i=t_i^m-1$, expresses
the local function through these residue blocks. In higher orders the
motion of the nodes contributes polynomial differential operators
applied to blocks of lower order, and the resulting change of
coordinates is triangular with identity diagonal
(Theorem~\ref{thm:rel-block-transport}). The image criterion is
therefore equivalent to
\[
 C_{\zeta,\mathbf b,h}(\mathbf z)
  \in\mathbb Q(\zeta)
  [(1+z_1)^{\pm1},\ldots,(1+z_d)^{\pm1},
                         D(1+\mathbf z)^{-1}]
\]
for every $\zeta,\mathbf b,h$. Clearing denominators turns this into
linear recurrences for the Newton coefficients, with recurrence
polynomials that depend on $D$ but not on the order of the root
(Corollary~\ref{cor:rel-block-recurrence}).

\subsection{The Alexander product and the type-\texorpdfstring{$A$}{A} lift}
\label{subsec:intro-alexander-divisors}

For a knot $K$ we normalize the Alexander polynomial by
\[
 \Delta_K(T^{-1})=\Delta_K(T),\qquad \Delta_K(1)=1.
\]
In rank one, the Melvin--Morton--Rozansky theorem identifies the
leading term of the large-color expansion with $\Delta_K^{-1}$
\cite{MelvinMorton,BarNatanGaroufalidis}. In free color coordinates
this is the constant quantum coefficient at $Q=1$. Rational forms of
the Kontsevich invariant, combined with Lie algebra weight systems,
give higher-rank expansions with Alexander denominators
\cite{RozanskyRational,Kricker,GK,Beads}.

The positive roots of $\mathfrak{sl}_n$ are the intervals
\[
 \alpha_{a,b}=\alpha_a+\cdots+\alpha_{b-1},
 \qquad 1\le a<b\le n.
\]
We attach to $\alpha_{a,b}$ the color monomial
$T_{a,b}=T_a\cdots T_{b-1}$ and set
\begin{equation}\label{eq:intro-alexander-divisor}
 D_{K,n}(\mathbf T)
    =\prod_{1\le a<b\le n}\Delta_K(T_{a,b}).
\end{equation}
For $\mathfrak{sl}_3$, for example,
\begin{equation}\label{eq:intro-alexander-sl3}
 D_{K,3}(T_1,T_2)
   =\Delta_K(T_1)\Delta_K(T_2)\Delta_K(T_1T_2).
\end{equation}
The third factor belongs to the root $\alpha_1+\alpha_2$ and couples
the two color directions; in higher rank the longer intervals play the
same role. Since $D_{K,n}(\mathbf1)=1$, the relative Habiro ring
$\mathcal H_{\mathrm{rel}}(D_{K,n})$ is defined.

\begin{maintheorem}[Relative lift for all dominant type-$A$ colors]
\label{main:type-a}
For every $n\ge2$ and every zero-framed knot $K\subset S^3$, there
is a unique element
\[
 \mathscr J_{K,n}
       \in\mathcal H_{\mathrm{rel}}(D_{K,n})_{\mathrm{pol}}
\]
whose value at $\mathbf r\in\mathbb N^{n-1}$ is
$J_{K,\sum_i r_i\omega_i}^{(n)}(Q)$.  Its Newton image is
$\mathcal J_{K,n}$.  Its local components satisfy
\[
 g_{K,n,1}(\mathbf t,0)=D_{K,n}(\mathbf t)^{-1},\qquad
 g_{K,n,\zeta}(\mathbf t,0)
   =\frac{A_{K,n,\zeta}(\mathbf t)}{D_{K,n}(\mathbf t^m)},
\]
where $\zeta$ has order $m$ and
$A_{K,n,\zeta}\in\mathbb Z[\zeta][\mathbf t^{\pm1}]$.
\end{maintheorem}

The colored polynomials therefore determine rational color expansions
at all roots of unity, linked by the relative gluing relations. At
$Q=1$ the leading coefficient is exactly the reciprocal Alexander
product. At a root of order $m$ the same product appears in the $m$th
powers of the colors, with an integral Laurent numerator, and higher
quantum orders allow higher powers of the denominator.

The component at $Q=1$ also recovers the Alexander polynomial. Exactly
$n-1$ positive roots contain $\alpha_1$, so on the first Dynkin axis
\begin{equation}\label{eq:intro-symmetric-alexander-recovery}
 g_{K,n,1}(t,1,\ldots,1;0)=\Delta_K(t)^{-(n-1)}.
\end{equation}
On the $i$th axis the exponent is $i(n-i)$. Because $\Delta_K(1)=1$, a
power series with constant term one has a unique root of each order
with constant term one, so $\Delta_K$ is determined by any of these
axis germs.

Mixed derivatives count overlaps of root intervals. Alexander symmetry
gives $\Delta_K'(1)=0$, and for type $A_2$
\[
 g_{K,3,1}((1+u,1+v),0)
   =1-\Delta_K''(1)(u^2+uv+v^2)+O((u,v)^3),
\]
where the mixed term $uv$ comes from $\alpha_1+\alpha_2$. In general
(Corollary~\ref{cor:cons-alexander-hessian})
\begin{equation}\label{eq:intro-cartan-hessian}
 \left.
 \left(\frac{\partial^2g_{K,n,1}(\mathbf t,0)}
 {\partial t_i\partial t_j}\right)_{i,j=1}^{n-1}
 \right|_{\mathbf t=\mathbf1}
   =-n\Delta_K''(1)C_{n-1}^{-1},
\end{equation}
where $C_{n-1}$ is the Cartan matrix of type $A_{n-1}$, with diagonal
entries $2$, entries $-1$ next to the diagonal, and zeros elsewhere.
The $(i,j)$ entry of $nC_{n-1}^{-1}$ is the number of root intervals
that contain both $i$ and $j$. Since the Newton basis becomes the
monomial basis at $Q=1$, this formula also gives the degree-two Newton
coefficients at $Q=1$.

\subsection{Rationality from braid operators and traces}
\label{subsec:intro-trace-method}

To prove Theorem~\ref{main:type-a} we start from the integral knot
element of Theorem~\ref{main:symmetric} and apply
Theorem~\ref{main:recognition}, so it remains to prove that the local
expansions are rational. At a nontrivial root of unity we use a braid
presentation: we cut one strand, fix its highest weight vector, and
take the quantum trace over the other strands. In rank one, the
relation between such traces, Burau matrices and the Alexander
polynomial is developed in \cite{RozanskyBurau,HuynhLe,MW}. In higher
rank each positive root carries its own variables, and the main task
is to control how different roots interact.

Fix a primitive $m$th root $\zeta$ with $m>1$. In an ordered
Poincar\'e--Birkhoff--Witt basis we write each root exponent as
\[
 N_\alpha=mA_\alpha+r_\alpha,\qquad 0\le r_\alpha<m.
\]
The residues $r_\alpha$ take finitely many values and span a
finite-dimensional factor, while the quotients $A_\alpha$ become
exponents of polynomial variables. After a fixed diagonal change of
basis, each quantum coefficient of a braid operator is a finite matrix
of polynomial differential operators followed by a triangular
polynomial substitution (Theorem~\ref{thm:braid-all-rank-insertion}).

The trace of a linear substitution is computed by
\[
 \sum_{k\ge0}s^k\operatorname{Tr}(\operatorname{Sym}^k B)
      =\det(I-sB)^{-1}.
\]
In one dimension this is the geometric series
$\sum_{k\ge0}(sb)^k=(1-sb)^{-1}$. Polynomial differential operators
inserted into the trace differentiate this expression and raise the
power of the denominator.

The nonlinear part of the substitution is controlled by root height,
$\operatorname{ht}(\alpha_{a,b})=b-a$. A nonlinear term replaces the
variable of an interval by a product of variables of shorter
intervals, and so it raises the ordinary polynomial degree, whereas a
differential operator of fixed order can lower the degree only by a
bounded amount. Only finitely many nonlinear terms therefore
contribute to a diagonal matrix coefficient, and what remains is a
finite combination of traces of linear substitutions. After the cut,
the linear
block of the root $\alpha_{a,b}$ has determinant
$\Delta_K(t_{a,b}^m)$ up to a Laurent unit, and the product over all
positive roots is $D_{K,n}(\mathbf t^m)$.

The parameter $s$ regularizes these infinite traces, and it has to be
removed before the comparison with colored invariants. Each
coefficient of $s$ is a trace on a finite-dimensional space. At an
integral color the cut linear block is nilpotent, so the regularized
determinant equals one there. Consequently every fixed joint jet in
the color and quantum variables is a polynomial in $s$, and it can be
evaluated at $s=1$.

For a dominant highest weight, the finite-dimensional irreducible
module is the quotient of the Verma module by its maximal proper
submodule. The permutation of the braid is a single cycle, because
its closure is a knot. In a basis adapted to the kernel of the
quotient map, a diagonal matrix coefficient can meet that kernel only
at a set of tensor positions that is a union of cycles of the
permutation, and the cut at the highest weight vector excludes the
full cycle. Only the irreducible quotient therefore contributes to the
partial trace, which is the reduced colored invariant
(Theorem~\ref{thm:trace-quotient-selection}). Agreement on the dominant
color grid then identifies the rational trace expression with the
Newton germ.

At $Q=1$ we use instead the rational Kontsevich invariant, its Lie
group weight system, and the comparison with quantum group invariants
\cite{GK,Beads,LM}. The leading matrix contribution is the inverse
square root of the product of the Alexander factors over all root
spaces. Opposite roots give equal factors, so the square root
normalized to one at the identity leaves one reciprocal factor for
each positive root. Diagrams of positive degree have poles along the
same divisors. This gives the component at $Q=1$ and completes the
rationality input for Theorem~\ref{main:recognition}.

The estimates on differential orders are uniform in the order of the
root. For fixed $n$ and a braid presentation $\beta$ of $K$ there is a
constant $C_{n,\beta}$, independent of $m>1$, such that
\[
 D_{K,n}(\mathbf t^m)^{1+C_{n,\beta}h}
       [x^h]g_{K,n,\zeta}(\mathbf t,x)
       \in\mathbb Z[\zeta][\mathbf t^{\pm1}]
\]
(Theorem~\ref{thm:cons-uniform-pole-order}). The exponent grows linearly
in the quantum order $h$. Its independence of $m$ comes from the fact
that the number of residue-lowering insertions that can survive in a
product is bounded in terms of $n$ alone.

\subsection{Examples and consequences}
\label{subsec:intro-examples}

For the figure-eight knot in rank one, Habiro's cyclotomic expansion
gives the Newton element
\[
 \mathcal J_{4_1,2}(t;Q)
  =\sum_{k\ge0}\prod_{j=1}^k
       \bigl(Qt+Q^{-1}t^{-1}-Q^j-Q^{-j}\bigr).
\]
The $k$th summand is divisible by $P_k(t;Q)$, so the series
terminates at every dominant color. At a primitive $m$th root $\zeta$,
a complete period has product
\[
 \prod_{j=1}^{m}
  \bigl(\zeta t+\zeta^{-1}t^{-1}-\zeta^j-\zeta^{-j}\bigr)
     =t^m+t^{-m}-2.
\]
Summing the geometric series over complete periods gives the
denominator $3-t^m-t^{-m}=\Delta_{4_1}(t^m)$, and the incomplete period
contributes a finite numerator. For instance,
\[
 g_{4_1,2,1}(t,0)=\frac1{\Delta_{4_1}(t)},\qquad
 g_{4_1,2,-1}(t,0)
      =\frac{\Delta_{4_1}(t)}{\Delta_{4_1}(t^2)}.
\]
In quantum order $h$ the exponent of the denominator is at most
$h+\lfloor h/2\rfloor+1$ at every root of unity, and at $Q=1$ this
exponent is the exact pole order at both zeros of the Alexander
polynomial. At the singular color the same element takes the value
\[
 \mathcal J_{4_1,2}(Q^{-1};Q)
   =\sum_{k\ge0}(-1)^kQ^{-k(k+1)/2}(Q;Q)_k^2\in\mathcal H.
\]

Torus knots show that the diagonal and mixed quadratic coefficients,
which agree at $Q=1$, separate at first order in $Q-1$. For
$K=T(2,\ell)$ with $\ell\ge1$ odd, the Rosso--Jones formula \cite{RJ}
computes the $\mathfrak{sl}_3$ invariants at $\omega_1$, $2\omega_1$
and $\omega_1+\omega_2$. The resulting Newton coefficients satisfy
\[
 c_{2,0}^{(\ell)}(1)=c_{1,1}^{(\ell)}(1)
       =-\frac{\ell^2-1}{4},
\]
as predicted by the Alexander Hessian, but
\[
 (c_{2,0}^{(\ell)})'(1)-(c_{1,1}^{(\ell)})'(1)
       =\frac{7(\ell^2-1)}8,
\]
so their first quantum corrections differ for every nontrivial knot in
the family. Section~\ref{sec:examples} also computes the corresponding
coefficients for $T(3,4)$ and the degree-two coefficients in type
$A_3$.

Connected sums give distinct elements with the same classical color
germ in all degrees. Let $T$ be the trefoil, $G=T\#T$, and
$S=T\#\overline T$. Their Alexander polynomials agree, so for every
$n\ge2$
\[
 g_{G,n,1}(\mathbf t,0)=g_{S,n,1}(\mathbf t,0)
       =D_{T,n}(\mathbf t)^{-2}.
\]
Their $\mathfrak{sl}_3$ invariants at the fundamental color
nevertheless give
\[
 c^G_{1,0}(Q)-c^S_{1,0}(Q)
       =48(Q-1)^2+O((Q-1)^3)
\]
in the conventions of Section~\ref{subsec:ex-connected-sums}, so
the two elements are separated by the quantum dependence of a single
Newton coefficient.

The Kinoshita--Terasaka knot, in the representative fixed in
Section~\ref{subsec:ex-alexander-trivial}, has Alexander
polynomial one. Its classical color germ is therefore one in every
rank, and all its local quantum coefficients are Laurent polynomials
in the colors. The lift is nevertheless not constant: its first
$\mathfrak{sl}_2$ Newton coefficient is
\[
 c_1(Q)=-Q^{-6}(Q-1)(Q^2+1)(Q^6-1)
       =-12(Q-1)^2+O((Q-1)^3).
\]

Finally, the Alexander divisors constrain difference equations in the
colors. Expanding an annihilating operator whose coefficients are
regular at $Q=1$ and Laurent polynomial through first order produces a
vector field in the color variables. Because the leading term is
$D_{K,n}^{-1}$, this vector field must preserve the ideal $(D_{K,n})$,
that is, it must be tangent to every irreducible component of the
divisor $D_{K,n}=0$ (Theorem~\ref{thm:cons-logarithmic-annihilator}).
In type $A_2$ this includes the mixed divisor $\Delta_K(t_1t_2)=0$.
For the trefoil and the figure-eight knot we give a symbol that passes
the two simple-root tests but fails the mixed one, and we construct an
operator that annihilates the local component modulo $(Q-1)^2$.

\subsection{Organization of the paper}
\label{subsec:intro-organization}

Section~\ref{sec:newton} develops Newton interpolation in
independent colors: Taylor expansions, singular boundaries, and
reconstruction from finitely many colors.
Section~\ref{sec:symmetric} identifies the symmetric invariants of the
regular-grid completion, proves that the extension is finite free, and
constructs the knot element in all dominant colors; it ends with the
comparison with symmetric-color Newton rings.

Section~\ref{sec:relative} proves the relative Habiro embedding
and the recognition theorem, together with the criterion in terms of
coefficient blocks. Sections~\ref{sec:braid}
and~\ref{sec:trace} prove local rationality for the knot
element. The first expresses braid operators at a root of unity
through polynomial differential operators and substitutions; the
second computes their traces and compares them with the colored
invariants. The argument at $Q=1$ at the end of
Section~\ref{sec:trace} completes the proof of
Theorem~\ref{main:type-a}.

Section~\ref{sec:consequences} contains the uniform pole
bound, the recovery of the Alexander polynomial, the inverse Cartan
Hessian, and the results on knot operations and difference equations.
Section~\ref{sec:examples} contains the explicit calculations for
specific knots.

\section{Newton interpolation in independent colors}
\label{sec:newton}

This section develops the Newton ring in $d$ independent colors. Its
elements are the integral interpolants of families of Laurent
polynomials on $\mathbb N^d$ (Theorem~\ref{thm:newton-interpolation}). We
then study their Taylor expansions at $Q=1$ and at other roots of
unity, their values on the first singular wall, and their
reconstruction from finitely many colors.

Let $d\ge1$, let $Q$ be an indeterminate, and put $A=\A$.
For a multiindex $\alpha\in\mathbb N^d$, write
$|\alpha|=\sum_i\alpha_i$ and $\mathbf t^\alpha=\prod_i t_i^{\alpha_i}$.
The symbol $\mathbf e_i$ denotes the $i$th standard basis vector.
For $k\in\mathbb N$ set
\[
 P_k(t)=\prod_{j=0}^{k-1}(t-Q^j),\qquad P_0(t)=1,
 \qquad
 P_{\boldsymbol\alpha}(\mathbf t)
   =\prod_{i=1}^dP_{\alpha_i}(t_i).
\]
For $\mathbf m=(m_1,\ldots,m_d)$ with $m_i\ge1$, let
\[
 A_{\mathbf m}=\A[t_1^{\pm1},\ldots,t_d^{\pm1}]/
       (P_{m_1}(t_1),\ldots,P_{m_d}(t_d)).
\]

\begin{definition}\label{def:newton-dominant}
The Newton ring on the dominant grid is
\[
 \mathcal N_d^+=\varprojlim_{\mathbf m}A_{\mathbf m},
\]
with the natural quotient maps as transition maps and the
inverse-limit topology, $\A$ being discrete. For a commutative
$\A$-algebra $B$ we write $\mathcal N_d(B)$ for the same construction
with coefficients in $B$. Thus $\mathcal N_d^+=\mathcal N_d(\A)$, which
we also write $\mathcal N_d^+(A)$, and $\mathcal N_0(B)=B$.
\end{definition}

The constant term of $P_m(t)$ is $\pm Q^{m(m-1)/2}$, a unit of $\A$, so
$t$ is invertible modulo $P_m(t)$. Hence the Laurent quotient in the
definition coincides with the polynomial quotient
$\A[\mathbf t]/(P_{m_i}(t_i))_i$, and successive division by the monic
polynomials $P_{m_i}$ shows that $A_{\mathbf m}$ is free over $\A$ of
rank $m_1\cdots m_d$.

\begin{theorem}\label{thm:newton-interpolation}
Every element of $\mathcal N_d^+$ has a unique expansion
\[
 F=\sum_{\boldsymbol\alpha\in\mathbb N^d}
       c_{\boldsymbol\alpha}(Q)P_{\boldsymbol\alpha}(\mathbf t),
 \qquad c_{\boldsymbol\alpha}\in\A.
\]
Evaluation on the dominant grid defines an injective homomorphism
\[
 \operatorname{ev}:\mathcal N_d^+
       \longrightarrow\prod_{\mathbf r\in\mathbb N^d}\A,
 \qquad F\longmapsto \bigl(F(Q^{r_1},\ldots,Q^{r_d})\bigr)_{\mathbf r}.
\]
Its image consists of the functions whose restriction to every finite
subset of $\mathbb N^d$ is represented by a single Laurent polynomial
over $\A$.
\end{theorem}

\begin{proof}
Fix a level $\mathbf m$ and put
\[
 I_{\mathbf m}=(P_{m_1}(t_1),\ldots,P_{m_d}(t_d)),\qquad
 \mathcal B_{\mathbf m}
 =\{P_{\boldsymbol\alpha}:0\le\alpha_i<m_i\}.
\]
Since $P_{\boldsymbol\alpha}$ is monic of multidegree
$\boldsymbol\alpha$, the change of basis from the monomials
$\mathbf t^{\boldsymbol\alpha}$, $0\le\alpha_i<m_i$, to
$\mathcal B_{\mathbf m}$ is unitriangular, so $A_{\mathbf m}$ is free
on $\mathcal B_{\mathbf m}$. If $\mathbf m'\le\mathbf m$, then every
$P_{\boldsymbol\alpha}$ with $\alpha_i\ge m_i'$ for some $i$ is
divisible by $P_{m_i'}(t_i)$. The transition map
$A_{\mathbf m}\to A_{\mathbf m'}$ therefore kills these basis elements
and keeps the others. An element of the inverse limit is the same as
an arbitrary family of coefficients
$(c_{\boldsymbol\alpha})_{\boldsymbol\alpha\in\mathbb N^d}$ in $\A$,
which gives the stated expansion.

The value of $P_{\boldsymbol\alpha}$ at the node $Q^{\mathbf r}$ is
zero unless $\boldsymbol\alpha\le\mathbf r$ componentwise. In one
variable the evaluation matrix
\[
 E_m=(P_a(Q^r))_{0\le r,a<m},\qquad
 P_a(Q^a)=\prod_{j<a}(Q^a-Q^j)\ne0,
\]
is therefore lower triangular with nonzero diagonal, and the
evaluation matrix on the box $0\le r_i<m_i$ is
$\bigotimes_i E_{m_i}$. Since $A$ is a domain, evaluation on this box
is injective on $A_{\mathbf m}$.

Let $f$ be a Laurent polynomial that vanishes on the box. After
multiplication by a monomial, which is a unit modulo
$I_{\mathbf m}$, we may divide $f$ by the monic polynomials
$P_{m_i}(t_i)$. The remainder is an $\A$-linear combination of
$\mathcal B_{\mathbf m}$ with zero values on the box, so it is zero by
the injectivity just proved. Consequently
\[
 \ker\!\left(A[\mathbf t^{\pm1}]
       \longrightarrow\prod_{0\le r_i<m_i}A\right)
 =I_{\mathbf m},
\]
so $A_{\mathbf m}$ is exactly the ring of functions on the box that
are represented by Laurent polynomials. A function on $\mathbb N^d$
that is represented by Laurent polynomials on all boxes therefore
defines compatible elements of all $A_{\mathbf m}$, that is, an
element of $\mathcal N_d^+$, and every element of $\mathcal N_d^+$
arises in this way. Since every finite set of colors lies in a box,
this proves the description of the image. Injectivity also follows:
an element with all values zero has zero image in every
$A_{\mathbf m}$.
\end{proof}

The coefficients in Theorem~\ref{thm:newton-interpolation} are iterated
divided differences of the colored values; see
\cite[Section~7, equation~(41)]{deBoor} for the classical formula,
which we apply in each coordinate.

\begin{proposition}\label{prop:newton-divided-differences}
Let $f:\mathbb N^d\to\mathbb Q(Q)$.  Its unique tensor Newton
coefficients are
\[
 c_{\boldsymbol\alpha}=
 \sum_{0\le\beta_i\le\alpha_i}
 \frac{f(\boldsymbol\beta)}
 {\displaystyle\prod_{i=1}^d
  \prod_{\substack{0\le j\le\alpha_i\\j\ne\beta_i}}
     (Q^{\beta_i}-Q^j)}.
\]
In particular an $\A$-valued function belongs to the image in
Theorem~\ref{thm:newton-interpolation} if and only if all these mixed
divided differences belong to $\A$.
\end{proposition}

\begin{proof}
In one variable, the divided difference of $f$ at
$Q^0,\ldots,Q^a$ is the leading coefficient of the interpolating
polynomial of degree at most $a$, that is, its coefficient of $P_a$,
and the displayed formula is the standard expression for it. Divided
differences in different variables commute, so applying the
one-variable formula in each variable in turn gives the coefficient of
$P_{\boldsymbol\alpha}$. If all coefficients lie in $\A$, the finite
sums $\sum_{\alpha_i<m_i}c_{\boldsymbol\alpha}P_{\boldsymbol\alpha}$
are integral Laurent polynomials representing $f$ on the boxes.
Conversely, if $f$ is represented on a box by a Laurent polynomial
over $\A$, division by the monic polynomials $P_{m_i}$ shows that its
Newton coordinates lie in $\A$.
\end{proof}

\begin{remark}\label{rem:newton-pairwise}
Pairwise congruences do not suffice, even on this quantum grid. The
values
\[
 f(0)=f(1)=0,\qquad f(2)=Q^2-1
\]
at the nodes $1,Q,Q^2$ satisfy $Q^a-Q^b\mid f(a)-f(b)$ in $\A$ for all
$a,b$, but their second Newton coefficient is
\[
 \frac{Q^2-1}{(Q^2-1)(Q^2-Q)}
       =\frac1{Q(Q-1)}\notin\A.
\]
\end{remark}

\subsection{Taylor expansions and ring structure}

Put $h=Q-1$ and $u_i=t_i-1$.

\begin{theorem}\label{thm:newton-taylor-domain}
Joint Taylor expansion defines an injective homomorphism
\[
 \mathcal N_d^+\hookrightarrow
       \mathbb Z[[h,u_1,\ldots,u_d]].
\]
In particular $\mathcal N_d^+$ is a domain.  For every rational prime
$p$, coefficient reduction induces an isomorphism
\[
 \mathcal N_d^+/p\mathcal N_d^+
 \simeq\varprojlim_{\mathbf m}
 \mathbb F_p[Q^{\pm1},\mathbf t^{\pm1}]/(P_{m_i}(t_i))_{i=1}^d,
\]
and this quotient is a domain as well.
\end{theorem}

\begin{proof}
Put $\mathfrak m=(h,u_1,\ldots,u_d)$. The substitution $Q=1+h$,
$t_i=1+u_i$ defines a ring homomorphism
$\tau:A[\mathbf t^{\pm1}]\to\mathbb Z[[h,u_1,\ldots,u_d]]$, because
$Q$ and the $t_i$ become units. Each factor
$t_i-Q^j=u_i-((1+h)^j-1)$ of a Newton polynomial lies in
$\mathfrak m$, so
\[
 c_{\boldsymbol\alpha}(1+h)
 P_{\boldsymbol\alpha}(\mathbf1+\mathbf u;1+h)
 \in\mathfrak m^{|\boldsymbol\alpha|}
\]
and $\tau(I_{\mathbf m})\subseteq\mathfrak m^{\min_im_i}$. Hence $\tau$
is continuous for the Newton topology and the $\mathfrak m$-adic
topology, and it extends to a ring homomorphism on $\mathcal N_d^+$,
given on Newton expansions by termwise substitution.

Suppose $\tau(F)=0$ and fix $\mathbf r\in\N^d$. The substitution
\[
 u_i=(1+h)^{r_i}-1
\]
is a continuous homomorphism $\mathbb Z[[h,\mathbf u]]\to\mathbb Z[[h]]$.
It sends $\tau(P_{\boldsymbol\alpha})$ to the Taylor series of
$P_{\boldsymbol\alpha}(Q^{\mathbf r};Q)$, which is zero unless
$\boldsymbol\alpha\le\mathbf r$, and therefore sends $\tau(F)$ to the
Taylor series at $Q=1$ of the Laurent polynomial
$F(Q^{r_1},\ldots,Q^{r_d})$. Since $A\hookrightarrow\Z[[h]]$, this
grid value is zero. As $\mathbf r$ was arbitrary,
Theorem~\ref{thm:newton-interpolation} gives $F=0$.

For the statement modulo $p$ we use that $p\mathcal N_d^+$ consists of
the elements whose Newton coefficients all lie in $pA$. Reduction of
the Newton coefficients modulo $p$ is therefore surjective with kernel
\[
 \{F:c_{\boldsymbol\alpha}\in pA\text{ for every }
       \boldsymbol\alpha\}=p\mathcal N_d^+,
\]
and it is compatible with multiplication at every finite level, which
gives the asserted isomorphism. Over $\F_p[Q^{\pm1}]$ the nodes $Q^r$
are still distinct and $\F_p[Q^{\pm1}]\hookrightarrow\F_p[[h]]$, so
Theorem~\ref{thm:newton-interpolation} and the argument above hold
verbatim and embed the quotient in $\F_p[[h,u_1,\ldots,u_d]]$. Both
rings in the statement are thus subrings of power series rings over
domains, hence domains.
\end{proof}

\begin{proposition}\label{prop:newton-quantum-jets}
For every $M\ge1$, joint Taylor expansion induces an isomorphism
\[
 \mathcal N_d^+/(Q-1)^M\mathcal N_d^+
 \simeq (\mathbb Z[h]/(h^M))[[u_1,\ldots,u_d]].
\]
The $(Q-1)$-adic completion of $\mathcal N_d^+$ is therefore
$\mathbb Z[[h,u_1,\ldots,u_d]]$.
\end{proposition}

\begin{proof}
As in the previous proof, $(Q-1)^M\mathcal N_d^+$ consists of the
elements whose Newton coefficients lie in $(Q-1)^MA$, so the quotient
is the Newton completion over $B_M=A/(Q-1)^M\simeq\mathbb Z[h]/(h^M)$.
It remains to identify this completion. In one variable,
\[
 P_k(1+u)=\prod_{j<k}\bigl(u-(Q^j-1)\bigr),\qquad Q^j-1\in hB_M.
\]
When the product is expanded, every term with at least $M$ factors of
the form $Q^j-1$ vanishes in $B_M$. Conversely, $P_k(1+u)=u^k+hg$ for
some $g\in B_M[u]$, so $u^k\in hB_M[u]$ modulo $P_k(1+u)$, and hence
$u^{kM}\in h^MB_M[u]=0$ modulo $P_k(1+u)$. Therefore
\[
 P_k(1+u)\in(u^{k-M+1})\quad(k\ge M),\qquad
 u^{kM}\in(P_k(1+u)),
\]
so the Newton ideals and the powers of $(u)$ define the same topology
on $B_M[u]$. Applying this in each variable gives
\[
 \mathcal N_d(B_M)=B_M[[u_1,\ldots,u_d]],
\]
and the identification is induced by the Taylor map of
Theorem~\ref{thm:newton-taylor-domain}. Taking the inverse limit over $M$
proves the final assertion.
\end{proof}

\begin{proposition}\label{prop:newton-color-degree}
For $F\in\mathcal N_d^+$ write
\[
 F(Q^{r_1},\ldots,Q^{r_d};Q)
       =\sum_{j\ge0}v_j(\mathbf r)(Q-1)^j.
\]
The function $v_j$ on $\mathbb N^d$ is the restriction of an
integer-valued polynomial on $\mathbb Z^d$, with rational
coefficients and total degree at most $j$.  In particular
\[
 v_0=c_{\mathbf0}(1),\qquad
 v_1=c'_{\mathbf0}(1)+\sum_{i=1}^d c_{\mathbf e_i}(1)r_i.
\]
\end{proposition}

\begin{proof}
Put $h=Q-1$. On the dominant grid every factor of a Newton polynomial
is
\[
 (1+h)^{r_i}-(1+h)^s
 =\sum_{a\ge1}\left(\binom{r_i}{a}-\binom{s}{a}\right)h^a.
\]
Its coefficient of $h^a$ is an integer-valued polynomial in $r_i$ of
degree at most $a$, and the factor is divisible by $h$. A term of
Newton index $\boldsymbol\alpha$ is therefore divisible by
$h^{|\boldsymbol\alpha|}$, and
\[
 v_j(\mathbf r)
 =\sum_{|\boldsymbol\alpha|\le j}
 [h^j]\,c_{\boldsymbol\alpha}(1+h)
 \prod_i\prod_{s<\alpha_i}
       \bigl((1+h)^{r_i}-(1+h)^s\bigr).
\]
This formula holds for every $\mathbf r\in\mathbb N^d$: a term with
$\alpha_i>r_i$ contains the factor with $s=r_i$, which is zero. In the
coefficient of $h^j$, each factor contributes a polynomial in $r_i$ of
degree at most its own order in $h$, so the total degree in
$\mathbf r$ is at most $j$, and sums of products of integer-valued
polynomials are integer-valued. The terms with
$|\boldsymbol\alpha|\le1$ give the formulas for $v_0$ and $v_1$.
\end{proof}

\begin{theorem}\label{thm:newton-units}
The unit group is
\[
 (\mathcal N_d^+)^\times
  =\{\pm Q^a t_1^{b_1}\cdots t_d^{b_d}:
          a,b_i\in\mathbb Z\}.
\]
\end{theorem}

\begin{proof}
Each $t_i$ is a unit in every $A_{\mathbf m}$, so these monomials are
units. Conversely, let $F$ be a unit. Its value at each grid point is
a unit of $\A$, say $f(\mathbf r)=\sigma(\mathbf r)Q^{a(\mathbf r)}$
with $\sigma(\mathbf r)=\pm1$. Expanding at $Q=1$ gives
$v_0(\mathbf r)=\sigma(\mathbf r)$ and
$v_1(\mathbf r)=\sigma(\mathbf r)a(\mathbf r)$. By
Proposition~\ref{prop:newton-color-degree} the sign is a constant
$\sigma$, and
\[
 \sigma a(\mathbf r)
      =c'_{\mathbf0}(1)+\sum_i c_{\mathbf e_i}(1)r_i.
\]
Hence $a(\mathbf r)=a_0+\sum_i b_i r_i$ with $a_0=a(\mathbf0)$ and
$b_i=a(\mathbf e_i)-a(\mathbf0)$, which are integers. The element
$\sigma Q^{a_0}\mathbf t^{\mathbf b}$ has the same values as $F$, so
the two agree by Theorem~\ref{thm:newton-interpolation}.
\end{proof}

\begin{proposition}\label{prop:newton-nonadic}
The Newton topology on $\mathcal N_d^+$ is not the $J$-adic topology
for any ideal $J$.
\end{proposition}

\begin{proof}
Suppose the topology were $J$-adic. Evaluation at a node $\mathbf r$
factors through a finite quotient and takes values in the discrete
ring $\A$, so it is continuous, and its kernel contains
$J^{a(\mathbf r)}$ for some $a(\mathbf r)\ge1$. For $x\in J$ this
gives $\operatorname{ev}_{\mathbf r}(x)^{a(\mathbf r)}=0$, hence
$\operatorname{ev}_{\mathbf r}(x)=0$ because $\A$ is a domain. Since
$\mathbf r$ is arbitrary, Theorem~\ref{thm:newton-interpolation} gives
$J=0$, and the $J$-adic topology would be discrete. The Newton
topology is not discrete: the kernel of the projection to any
$A_{\mathbf m}$ contains the nonzero element $P_{m_1}(t_1)$.
\end{proof}

\subsection{Changes of color}

\begin{proposition}\label{prop:newton-operations}
The following operations are defined on $\mathcal N_d^+$.
\begin{enumerate}
\item The forward shifts
\[
 E_iF(\mathbf t)=F(t_1,\ldots,Qt_i,\ldots,t_d)
\]
are continuous injective endomorphisms.
\item For every integer $a\ge1$, the substitutions
\[
 \Psi_a(Q)=Q^a,\qquad\Psi_a(t_i)=t_i^a
\]
are continuous endomorphisms, satisfy
$\Psi_a\Psi_b=\Psi_{ab}$, and obey
$\Psi_p(F)\equiv F^p\pmod{p\mathcal N_d^+}$ for every prime $p$.
\item Simultaneous inversion $Q\mapsto Q^{-1}$ and
$t_i\mapsto t_i^{-1}$ is a continuous involution.
\end{enumerate}
\end{proposition}

\begin{proof}
(1) We use the description of $\mathcal N_d^+$ in
Theorem~\ref{thm:newton-interpolation}. If a Laurent polynomial $f$
represents $F$ on $S+\mathbf e_i$, then
$f(t_1,\ldots,Qt_i,\ldots,t_d)$ represents
$\mathbf r\mapsto F(\mathbf r+\mathbf e_i)$ on $S$. So $E_iF$ lies in
$\mathcal N_d^+$, $E_i$ is a ring homomorphism, and it is continuous
because the values of $E_iF$ on a box only involve the values of $F$
on the shifted box. For injectivity, suppose $E_iF=0$, so that $F$
vanishes at every node with $r_i\ge1$, and fix the other coordinates
of $\mathbf r$. The resulting one-variable function $g$ is represented
on $\{0,k\}$ by a Laurent polynomial, so $g(k)-g(0)$ is divisible by
$Q^k-1$. Since $g(k)=0$ for $k\ge1$, the Laurent polynomial $g(0)$ is
divisible by $Q^k-1$ for every $k\ge1$. A nonzero Laurent polynomial
has no divisors of arbitrarily large degree, so $g(0)=0$, and all
values of $F$ vanish.

(2) Since $t-Q^j$ divides $t^a-Q^{aj}$, the polynomial
$\Psi_a(P_k(t;Q))=\prod_{j<k}(t^a-Q^{aj})$ is divisible by
$P_k(t;Q)$. Hence $\Psi_a$ maps each ideal $I_{\mathbf m}$ into itself
and extends continuously to the inverse limit. The rule
$\Psi_a\Psi_b=\Psi_{ab}$ holds on Laurent polynomials and therefore on
$\mathcal N_d^+$. For a Laurent polynomial $f$ with integer
coefficients, $f(Q^p,\mathbf t^p)\equiv f(Q,\mathbf t)^p\pmod p$.
Writing $F$ as a limit of its Newton truncations $f_k$, we get
$\Psi_p(F)-F^p=\lim_k(\Psi_p(f_k)-f_k^p)$, a limit of elements of
$p\mathcal N_d^+$. This subset is closed, since it consists of the
elements whose Newton coefficients all lie in $pA$.

(3) The identity
\[
 P_k(t^{-1};Q^{-1})
    =(-1)^k t^{-k}Q^{-\binom{k}{2}}P_k(t;Q)
\]
shows that simultaneous inversion maps each $I_{\mathbf m}$ to
itself. It is an involution of the Laurent polynomial ring, and hence
of $\mathcal N_d^+$.
\end{proof}

\begin{remark}\label{rem:newton-weyl}
The shifted Weyl group does not preserve the dominant chamber. In
type $A$ the simple reflection at $i$ acts on Dynkin labels by
$r_i\mapsto-r_i-2$ and by $r_j\mapsto r_j+r_i+1$ at the adjacent
vertices $j$, which corresponds to the substitutions
$t_i\mapsto Q^{-2}t_i^{-1}$ and $t_j\mapsto Qt_it_j$.
These substitutions move dominant nodes out of the chamber, so they do
not act on the dominant completion. One can instead complete on a
Weyl-stable regular grid, as in Section~\ref{sec:symmetric}, or show
that a particular element extends beyond the chamber. Along the
singular wall $r_i=-1$ every element has the Habiro-valued boundary
value constructed in Section~\ref{subsec:newton-walls}, and requiring this
value to lie in $\A$ is a genuine condition
(Theorem~\ref{thm:newton-wall-fiber-product}).
\end{remark}

\begin{proposition}\label{prop:newton-padic-colors}
Let $p$ be a prime and fix $q_0\ne1$, with $q_0\in1+p\Z_p$
if $p$ is odd and $q_0\in1+4\Z_2$ if $p=2$.
Put $a=v_p(q_0-1)$. Every
$F\in\mathcal N_d^+$ defines a restricted analytic function on $\Z_p^d$
by $Q=q_0$ and $t_i=q_0^{z_i}$. The part of its Newton expansion with
$|\alpha|\ge N$ has Tate-algebra norm at most $p^{-aN}$.
\end{proposition}

\begin{proof}
Under these hypotheses $\log q_0\in p^a\Z_p$, and
$q_0^z=\exp(z\log q_0)$ converges in the Tate algebra
$\Z_p\langle z\rangle$. Indeed, since $v_p(k!)\le(k-1)/(p-1)$ and
$a>1/(p-1)$,
\[
 v_p\!\left(\frac{(\log q_0)^k}{k!}\right)
 =ka-v_p(k!)\ge a\quad(k\ge1),
\]
and these valuations tend to infinity. Thus $q_0^z-1$ has norm at most
$p^{-a}$, and so does $q_0^z-q_0^j=(q_0^z-1)-(q_0^j-1)$. It follows
that
\[
 \|q_0^z-q_0^j\|\le p^{-a},\qquad
 \|c_\alpha(q_0)P_\alpha(q_0^{\mathbf z};q_0)\|
       \le p^{-a|\alpha|},
\]
where $c_\alpha(q_0)\in\Z_p$ because $q_0$ is a unit. There are only
finitely many indices of each total degree, so the Newton series
converges in $\Z_p\langle z_1,\ldots,z_d\rangle$, and the
non-Archimedean triangle inequality gives the tail bound.
\end{proof}

The specialization need not be injective. If $q_0$ is algebraic over
$\mathbb Q$, any nonzero polynomial in $\Z[Q]$ that vanishes at $q_0$
is a nonzero element of its kernel.

\subsection{Singular colors}
\label{subsec:newton-walls}

Write $\Phi_m(Q)$ for the $m$th cyclotomic polynomial. The Habiro
ring \cite{Habiro} is
\[
 c_N=(Q;Q)_N=\prod_{j=1}^N(1-Q^j),\qquad
 \Habi=\varprojlim_N\A/(c_N).
\]
The natural map $\A\to\Habi$ is injective, because a nonzero Laurent
polynomial is not divisible by every $c_N$.
For $d\ge1$ we consider the first simple singular wall $t_1=Q^{-1}$,
with normal coordinate
\[
 u=t_1-Q^{-1}.
\]

\begin{theorem}\label{thm:newton-wall-jets}
For every $M\ge1$, normal expansion induces a natural isomorphism
\[
 \mathcal N_d^+/u^M\mathcal N_d^+
   \simeq\mathcal N_{d-1}(\Habi)[u]/(u^M).
\]
In particular there is a surjective boundary homomorphism
\[
 \beta_1:\mathcal N_d^+\longrightarrow\mathcal N_{d-1}(\Habi),
 \qquad \ker\beta_1=u\mathcal N_d^+.
\]
The full normal Taylor map
\[
 \mathcal N_d^+\hookrightarrow\mathcal N_{d-1}(\Habi)[[u]]
\]
is injective, and its target is the $u$-adic completion of the source.
\end{theorem}

\begin{proof}
We first express the quotient as a limit of finite quotients. At a
box level $\mathbf m$, multiplication by $u=t_1-Q^{-1}$ becomes
diagonal after extending scalars to $\Q(Q)$, with eigenvalues
$Q^{r_1}-Q^{-1}\ne0$. Since $A_{\mathbf m}$ is free over $A$,
multiplication by $u^M$ is injective, and
\[
 0\longrightarrow A_{\mathbf m}
 \xrightarrow{\ u^M\ }A_{\mathbf m}
 \longrightarrow A_{\mathbf m}/u^MA_{\mathbf m}
 \longrightarrow0
\]
is exact. The levels $\mathbf m=(k,\ldots,k)$ are cofinal, and the
left-hand system has surjective transition maps, so the inverse limit
of these sequences is exact and identifies
$\mathcal N_d^+/u^M\mathcal N_d^+$ with the limit of the finite
cokernels.

Next we treat the first variable modulo $u^M$. Since
$t_1=Q^{-1}+u$, we have $A[t_1]=A[u]$, and we work in $A[u]/(u^M)$.
The constant term of $P_k(Q^{-1}+u)$ is $P_k(Q^{-1})=Q^{-k}c_k$. The
coefficient of $u^j$, $j<M$, is a sum of products obtained from
$\prod_{j'<k}(Q^{-1}-Q^{j'})$ by deleting $j$ factors, and up to units
these factors are the $1-Q^i$, $1\le i\le k$. The product $c_k$
contains $\Phi_\ell(Q)$ with multiplicity $\lfloor k/\ell\rfloor$, and
each factor $1-Q^i$ contains it at most once. Every coefficient of
$P_k(Q^{-1}+u)$ in $A[u]/(u^M)$ therefore contains $\Phi_\ell$ with
multiplicity at least
\[
 \max\!\left(\left\lfloor\frac{k}{\ell}\right\rfloor-M+1,0\right).
\]
For $k=MN$ and $\ell\le N$ this is at least $\lfloor N/\ell\rfloor$,
because $\lfloor MN/\ell\rfloor-M+1\ge M(\lfloor N/\ell\rfloor-1)+1
\ge\lfloor N/\ell\rfloor$. As $c_N=\pm\prod_{\ell\le N}\Phi_\ell^{\lfloor
N/\ell\rfloor}$, unique factorization in $A$ gives
\[
 (P_{MN}(Q^{-1}+u))\subseteq(c_N)
       \quad\text{in }A[u]/(u^M).
\]
In the other direction, $P_k(Q^{-1}+u)=Q^{-k}c_k+ug$ with
$g\in A[u]$, so $c_k$ is a multiple of $u$ modulo
$P_k(Q^{-1}+u)$, and $c_k^M$ lies in the ideal generated by
$P_k(Q^{-1}+u)$ in $A[u]/(u^M)$. Moreover $c_k^M$ divides $c_{Mk}$,
since $\lfloor Mk/\ell\rfloor\ge M\lfloor k/\ell\rfloor$. Hence
\[
 (c_{Mk})\subseteq(c_k^M)\subseteq(P_k(Q^{-1}+u))
 \quad\text{in }A[u]/(u^M).
\]
The two families of ideals are cofinal, so the completion of
$A[u]/(u^M)$ along the Newton ideals in $t_1$ coincides with its
completion along the ideals $(c_N)$. Since $A[u]/(u^M)$ is free of
rank $M$ over $A$, this completion is $\mathcal H[u]/(u^M)$.

The other variables keep their tensor Newton coordinates: the finite
quotient $A_{\mathbf m}/u^MA_{\mathbf m}$ is the tensor product of
$A[u]/(u^M,P_{m_1}(Q^{-1}+u))$ with the free $A$-module on the
$P_{\boldsymbol\alpha'}(t_2,\ldots,t_d)$, $\alpha'_i<m_i$. Passing to
the limit gives
\[
 \mathcal N_d^+/u^M\mathcal N_d^+
       \simeq\mathcal N_{d-1}(\mathcal H)[u]/(u^M).
\]
The case $M=1$ is the assertion about the boundary map.

Finally, let $F\in\bigcap_{M\ge1}u^M\mathcal N_d^+$. At every dominant
node
\[
 F(Q^{r_1},\ldots,Q^{r_d})
 \in\bigcap_{M\ge1}(Q^{r_1}-Q^{-1})^MA=0,
\]
because $Q^{r_1}-Q^{-1}$ is a nonzero nonunit of $A$ and its powers
have unbounded degree. Injectivity of grid evaluation gives $F=0$.
Thus the map from $\mathcal N_d^+$ to the inverse limit of the jet
rings, which is $\mathcal N_{d-1}(\Habi)[[u]]$, is injective, and by
the jet isomorphisms this inverse limit is the $u$-adic completion of
$\mathcal N_d^+$.
\end{proof}

Explicitly, the boundary is
\[
 \beta_1(F)=
 \sum_{\boldsymbol\alpha'\in\mathbb N^{d-1}}
 \left(\sum_{k\ge0}c_{(k,\boldsymbol\alpha')}(Q)Q^{-k}c_k(Q)\right)
 P_{\boldsymbol\alpha'}(t_2,\ldots,t_d),
\]
where each inner sum converges in $\Habi$.

\begin{theorem}\label{thm:newton-wall-fiber-product}
Let $\mathcal N_{\rm ext}$ be the tensor Newton completion on the
grid $\{-1,0,1,\ldots\}\times\mathbb N^{d-1}$, with the first
node sequence $Q^{-1},1,Q,Q^2,\ldots$.  There is a natural isomorphism
\[
 \mathcal N_{\rm ext}\simeq
 \mathcal N_{d-1}(\A)
       \times_{\mathcal N_{d-1}(\Habi)}\mathcal N_d^+.
\]
An element of $\mathcal N_d^+$ therefore extends across this wall
with Laurent-valued interpolation if and only if its boundary lies
in the embedded subring $\mathcal N_{d-1}(\A)$.
\end{theorem}

\begin{proof}
For the node sequence $Q^{-1},1,Q,\ldots$, the Newton basis in the
first variable is
\[
 1,\ u,\ uP_1(t_1),\ uP_2(t_1),\ \ldots,\qquad u=t_1-Q^{-1}.
\]
Separating the terms that contain the basis element $1$ from the
others, we can write every element of $\mathcal N_{\rm ext}$ uniquely
as
$F_0+uG$ with $F_0\in\mathcal N_{d-1}(\A)$ and $G\in\mathcal N_d^+$.
Its restriction to the wall $t_1=Q^{-1}$ is $F_0$, and its restriction
to the dominant grid is $F_0+uG\in\mathcal N_d^+$, whose boundary is
the image of $F_0$ in $\mathcal N_{d-1}(\Habi)$. This defines a map to
the fiber product. It is injective, because $F_0=0$ and $uG=0$ imply
$G=0$; multiplication by $u$ is injective on $\mathcal N_d^+$ by
Theorem~\ref{thm:newton-wall-jets}.

Conversely, suppose that $F_0\in\mathcal N_{d-1}(\A)$ and
$F\in\mathcal N_d^+$ have the same image in
$\mathcal N_{d-1}(\Habi)$. Regard $F_0$ as an element of
$\mathcal N_d^+$ independent of $t_1$. Then $\beta_1(F-F_0)=0$, so
Theorem~\ref{thm:newton-wall-jets} gives $F-F_0=uG$ with
$G\in\mathcal N_d^+$, and $F_0+uG\in\mathcal N_{\rm ext}$ has the
prescribed restrictions. This proves surjectivity and the extension
criterion.
\end{proof}

\begin{remark}\label{rem:newton-wall-range}
The same proof works for the wall $t_1=Q^{-s}$ for every $s\ge1$.
Up to a Laurent unit, $P_k(Q^{-s})$ equals
$\prod_{j=s}^{s+k-1}(1-Q^j)$, which contains $\Phi_\ell$ with
multiplicity at least $\lfloor k/\ell\rfloor$ and divides $c_{s+k-1}$.
These products therefore define the Habiro topology, and the
multiplicity estimates in the proof of
Theorem~\ref{thm:newton-wall-jets} hold as stated.
\end{remark}

\subsection{Reconstruction from finitely many colors}

Let $d\ge1$ and let a Newton series be written as
\[
 F(\mathbf t,Q)=\sum_{\alpha\in\N^d}c_\alpha(Q)
       \prod_{i=1}^dP_{\alpha_i}(t_i;Q),\qquad c_\alpha\in\A.
\]
For $a\ge0$ and $0\le r\le a$, set
\[
 \begin{bmatrix}a\\r\end{bmatrix}_Q=
 \frac{\prod_{j=1}^a(1-Q^j)}
 {\prod_{j=1}^r(1-Q^j)\prod_{j=1}^{a-r}(1-Q^j)}\in A,
\]
and put $\begin{bmatrix}a\\r\end{bmatrix}_Q=0$ otherwise. These are the
Gaussian binomial coefficients in the variable $Q$. They are
polynomials, by the Pascal identity
\[
 \begin{bmatrix}a\\r\end{bmatrix}_Q
 =\begin{bmatrix}a-1\\r-1\end{bmatrix}_Q
       +Q^r\begin{bmatrix}a-1\\r\end{bmatrix}_Q
\]
and the value one at $r=0,a$. The normalized Newton factors satisfy
\[
 \frac{P_k(Q^r;Q)}{P_k(Q^k;Q)}
       =\begin{bmatrix}r\\k\end{bmatrix}_Q
       \qquad(r\ge0).
\]
They are the quantum binomial polynomials of \cite{HarmanHopkins}
under the change of variable $t=1+(Q-1)y$:
if $B_k(y)$ denotes the basis polynomial of that paper, then
$B_k((t-1)/(Q-1))=P_k(t;Q)/P_k(Q^k;Q)$.

Fix a primitive $m$th root $\zeta$ and residues $0\le b_i<m$, and use
the local coordinates $x=Q-\zeta$ and $w_i=t_i-\zeta^{b_i}$.
For $L\ge0$, the total $L$ jet at $(\zeta,\zeta^{\mathbf b})$ is the
reduction modulo $(x,w_1,\ldots,w_d)^{L+1}$. For $a\ge0$ put
\[
 \nu_{m,b}(a)=\#\{0\le j<a:j\equiv b\pmod m\},\qquad
 \nu_{\mathbf b}(\alpha)=\sum_i\nu_{m,b_i}(\alpha_i).
\]

\begin{theorem}\label{thm:newton-finite-jets}
For every integer $L\ge0$, the joint total-order $L$ jet of $F$ at
$(\zeta,\zeta^{\mathbf b})$ depends only on the coefficients
$c_\alpha$ with $\nu_{\mathbf b}(\alpha)\le L$, and every such index
satisfies
\begin{equation}\label{eq:newton-jet-cutoff}
 \alpha_i\le mL+b_i,\qquad
 |\alpha|\le mL+\sum_i b_i\le mL+d(m-1).
\end{equation}
This joint jet is determined by the ordinary quantum $L$ jets at
$Q=\zeta$ of the finitely many colored values
\[
 J_{\mathbf r}(Q)=F(Q^{r_1},\ldots,Q^{r_d},Q),\qquad
 0\le r_i\le m(L+1)-1.
\]
\end{theorem}

\begin{proof}
At the chosen center, $t_i-Q^j$ vanishes exactly when
$j\equiv b_i\pmod m$. In that case
$t_i-Q^j=w_i+\zeta^{b_i}-(\zeta+x)^j$ has the nonzero linear part
$w_i-j\zeta^{j-1}x$, and otherwise it is a unit. Hence
\[
 \ord_{(x,\mathbf w)}P_\alpha=\nu_{\mathbf b}(\alpha).
\]
The $L$ jet of $c_\alpha P_\alpha$ is zero if
$\nu_{\mathbf b}(\alpha)>L$; otherwise it depends only on the Taylor
coefficients of $c_\alpha$ at $\zeta$ of order at most
$L-\nu_{\mathbf b}(\alpha)$. If $a=mu+s$ with $0\le s<m$, then
\[
 \nu_{m,b}(a)=u+\mathbf1_{s>b},\qquad
 a\le m\nu_{m,b}(a)+b.
\]
So $\nu_{\mathbf b}(\alpha)\le L$ implies
$\alpha_i\le m\nu_{m,b_i}(\alpha_i)+b_i\le mL+b_i$, and summing the
inequalities $\alpha_i\le m\nu_{m,b_i}(\alpha_i)+b_i$ over $i$ gives
\eqref{eq:newton-jet-cutoff}.

It remains to recover the Taylor coefficients of $c_\alpha$. Set
\[
 d_a=P_a(Q^a;Q)
    =Q^{\binom a2}\prod_{s=1}^a(Q^s-1),\qquad
 d_\alpha=\prod_i d_{\alpha_i}.
\]
Tensor Lagrange interpolation gives
\begin{equation}\label{eq:newton-gaussian-difference}
 d_\alpha c_\alpha=
 \sum_{0\le\mathbf r\le\alpha}
 \prod_i\left(
 (-1)^{\alpha_i-r_i}Q^{\binom{\alpha_i-r_i}{2}}
 \begin{bmatrix}\alpha_i\\r_i\end{bmatrix}_Q\right)
 J_{\mathbf r}(Q).
\end{equation}
In one variable the Lagrange denominator at the node $Q^r$ is
\[
 (-1)^{a-r}Q^{\binom r2+r(a-r)}
 \prod_{s=1}^r(Q^s-1)\prod_{s=1}^{a-r}(Q^s-1),
\]
and dividing $d_a$ by it gives the weight in
\eqref{eq:newton-gaussian-difference}.

Exactly the factors $Q^s-1$ with $m\mid s$ vanish at $\zeta$, each to
first order, so $d_\alpha$ vanishes at $\zeta$ to order
\[
 v_m(d_\alpha)=\sum_i\left\lfloor\frac{\alpha_i}{m}\right\rfloor
       \le\nu_{\mathbf b}(\alpha).
\]
If the right side of \eqref{eq:newton-gaussian-difference} is known
modulo $(Q-\zeta)^{L+1}$, division by $d_\alpha$ in
$\mathbb Q(\zeta)[[Q-\zeta]]$ determines $c_\alpha$ modulo
$(Q-\zeta)^{L+1-v_m(d_\alpha)}$. Since
$L+1-v_m(d_\alpha)\ge L+1-\nu_{\mathbf b}(\alpha)$, this is all the jet
requires. The right side modulo $(Q-\zeta)^{L+1}$ is determined by the
$L$ jets of the $J_{\mathbf r}$, and the colors involved satisfy
$r_i\le\alpha_i\le mL+b_i\le m(L+1)-1$.
\end{proof}

\begin{remark}\label{rem:newton-jet-conventions}
For $m=1$ the conditions become $|\alpha|\le L$ and
$0\le r_i\le L$. For $L=0$ only the indices with $\alpha_i\le b_i$
contribute, and for them $d_\alpha$ does not vanish at $\zeta$. The
rescaling $Q-\zeta=\zeta\epsilon$, $w_i=\zeta^{b_i}u_i$ changes
coordinates by cyclotomic units and does not affect the bounds. If the
colored values are only known modulo a power of a prime $p$, the
division by $d_\alpha$ may require more $p$-adic precision, because the
unit part of $d_\alpha$ at $\zeta$, namely
$d_\alpha/(Q-\zeta)^{v_m(d_\alpha)}$ evaluated at $\zeta$, can be
divisible by primes dividing $m$.
\end{remark}

\section{Symmetric interpolation and regular colors}
\label{sec:symmetric}

In this section we pass from symmetric interpolation at partition
nodes to interpolation on the ordered regular grid, prove that the
second is a finite free extension of the first, and use it to
construct the integral knot element in all dominant colors.

Throughout the section, $n\ge2$, $A=\A$,
$R_n=\A[x_1,\ldots,x_n]$, and
$\Lambda_n=R_n^{\mathfrak S_n}$, where $\mathfrak S_n$ acts by
permuting the variables. Put
\[
 \boldsymbol\delta=(n-1,n-2,\ldots,0),\qquad
 \Delta(\mathbf x)=\prod_{i<j}(x_i-x_j).
\]
Here $\Delta(\mathbf x)$ is the Vandermonde polynomial; the Alexander
polynomial of a knot always carries a subscript $K$.
A partition of length at most $n$ is written
$\lambda=(\lambda_1\ge\cdots\ge\lambda_n\ge0)$, with zero parts
included; we set $|\lambda|=\sum_i\lambda_i$ and write
$\lambda\subseteq\mu$ when $\lambda_i\le\mu_i$ for every $i$.
Define
\begin{equation}\label{eq:sym-monic-schur}
 S_\lambda(\mathbf x;Q)=
 \frac{\det[P_{\lambda_i+n-i}(x_j;Q)]_{i,j=1}^n}
      {\Delta(\mathbf x)}.
\end{equation}
The numerator is alternating, so the quotient is a symmetric
polynomial over $\A$, and its top homogeneous component is the Schur
polynomial $s_\lambda$; see \cite[Chapter~I]{Macdonald}.
Since the $s_\lambda$ form a basis of the symmetric polynomials, the
$S_\lambda$ form an $\A$-basis of $\Lambda_n$. They are the factorial
Schur polynomials for the parameter sequence $a_j=Q^{j-1}$ in
\cite[Section~2]{MolevSagan}.

The polynomials in \eqref{eq:sym-monic-schur} are the monic
normalization of the interpolation polynomials of Beliakova and
Gorsky \cite[Section~10.2, equation~(39)]{BG}, whose one-variable
building blocks are
$f_a(X;q)=(X;q)_a=\prod_{j=0}^{a-1}(1-Xq^j)$. Since
$f_a(X;Q^{-1})=(-1)^aQ^{-\binom a2}P_a(X;Q)$, these factors can be
pulled out of the determinant row by row, which gives
\begin{equation}\label{eq:sym-unit-conversion}
 F_\lambda^{\mathrm{BG}}(\mathbf x;Q^{-1})
 =(-1)^{|\lambda|+\binom n2}
  Q^{-\sum_i\binom{\lambda_i+n-i}{2}}S_\lambda(\mathbf x;Q).
\end{equation}
The two normalizations therefore differ only by Laurent units.

\begin{lemma}\label{lem:sym-partition-vanishing}
For $p_\mu=(Q^{\mu_1+n-1},\ldots,Q^{\mu_n})$ one has
$S_\lambda(p_\mu;Q)=0$ unless $\lambda\subseteq\mu$, and
$S_\mu(p_\mu;Q)\ne0$. If $I_d$ is the ideal of symmetric polynomials
vanishing at every $p_\mu$ with $|\mu|<d$, then
\[
 I_d=\bigoplus_{|\lambda|\ge d}\A S_\lambda.
\]
Hence the partition interpolation completion is the topological ring
\begin{equation}\label{eq:sym-pure-completion}
 \mathcal B_n=\varprojlim_d\Lambda_n/I_d
             =\prod_{\ell(\lambda)\le n}\A S_\lambda,
\end{equation}
and its evaluations at all partition nodes are jointly injective.
\end{lemma}

\begin{proof}
The vanishing is the interpolation property of factorial Schur
polynomials \cite[Theorem~2.1]{MolevSagan}; in the present
normalization, see \cite[Corollary~10.9]{BG}. It can also be read off
from the determinant. The entry in row $i$ and column $k$ is
$P_{\lambda_i+n-i}(Q^{\mu_k+n-k})$, which vanishes when
$\lambda_i+n-i>\mu_k+n-k$. If $\lambda_j>\mu_j$ for some $j$, this
inequality holds whenever $i\le j\le k$, so the first $j$ rows vanish
in the columns $j,\ldots,n$. These $j$ rows are then supported in
$j-1$ columns, and every term of the determinant is zero. For
$\lambda=\mu$ the same inequality shows that the matrix is lower
triangular. Its diagonal entries $P_{\mu_i+n-i}(Q^{\mu_i+n-i})$ and
the Vandermonde factor at $p_\mu$ are nonzero in
$\operatorname{Frac}(\A)$.

The description of $I_d$ is a special case of the following
statement. Let $\mathscr P$ be a set of partitions that is closed
under passing to smaller partitions, and let
$f=\sum_\lambda b_\lambda S_\lambda$ vanish at every $p_\mu$ with
$\mu\in\mathscr P$. Then $b_\mu=0$ for all $\mu\in\mathscr P$: arguing
by induction on $|\mu|$, we may assume $b_\lambda=0$ for
$\lambda\subsetneq\mu$, and then
$0=f(p_\mu)=\sum_{\lambda\subseteq\mu}b_\lambda S_\lambda(p_\mu)
=b_\mu S_\mu(p_\mu)$, so $b_\mu=0$ because $\A$ is a domain.
Conversely, $S_\lambda$ with $\lambda\notin\mathscr P$ vanishes at
every $p_\mu$ with $\mu\in\mathscr P$, since $\lambda\not\subseteq\mu$.
So the symmetric polynomials vanishing on these nodes are exactly
$\bigoplus_{\lambda\notin\mathscr P}\A S_\lambda$. Taking for
$\mathscr P$ the partitions of size less than $d$ gives the formula
for $I_d$. Since the $I_d$ are vanishing ideals, multiplication passes
to the inverse limit, and the same induction, applied to an infinite
sum $\sum_\lambda b_\lambda S_\lambda$ whose value at $p_\mu$ is the
finite sum over $\lambda\subseteq\mu$, gives joint injectivity.
\end{proof}

\subsection{The ordered regular grid}

The Weyl translates of the shifted highest-weight points $p_\mu$
form an ordered grid of tuples with distinct coordinates.
For $L\ge n$ put
\[
 \Omega_L=\{(Q^{a_1},\ldots,Q^{a_n}):0\le a_i<L,\ a_i\ne a_j\},
\]
and let
\[
 T_L=(P_L(x_1),\ldots,P_L(x_n)),\qquad
 K_L=\ker(R_n\longrightarrow\A^{\Omega_L}).
\]
Regularity is a condition on the exponent tuples; $Q$ remains an
indeterminate.

\begin{theorem}\label{thm:sym-regular-finite-presentation}
Let $P_L[x_1,\ldots,x_i]$ denote an ordinary polynomial divided
difference. Then
\begin{equation}\label{eq:sym-regular-ideal}
 K_L=T_L:\Delta
 =\bigl(P_L[x_1],P_L[x_1,x_2],\ldots,P_L[x_1,\ldots,x_n]\bigr).
\end{equation}
Moreover, $R_n/K_L$ is free over $\A$, of rank
$(L)_n=L(L-1)\cdots(L-n+1)$, with basis
\begin{equation}\label{eq:sym-regular-finite-basis}
 x_1^{a_1}\cdots x_n^{a_n},\qquad 0\le a_i<L-i+1.
\end{equation}
\end{theorem}

\begin{proof}
Let $K=\Frac(A)$. Division by the monic polynomials $P_L(x_i)$ makes
$R_n/T_L$ a free $A$-module, and after extending scalars to $K$,
evaluation identifies it with the ring of functions on the full
Cartesian grid $\{(Q^{a_1},\ldots,Q^{a_n}):0\le a_i<L\}$. Hence
evaluation on the Cartesian grid is injective over $A$, and $T_L$ is
the ideal of polynomials vanishing on that grid. Since $\Delta(p)$ is
nonzero exactly at the regular points,
\[
 \Delta f\in T_L
 \ \Longleftrightarrow\
 \Delta(p)f(p)=0\text{ on the Cartesian grid}
 \ \Longleftrightarrow\ f\in K_L.
\]
This proves $K_L=T_L:\Delta$.

Put $g_i=P_L[x_1,\ldots,x_i]$ and $J=(g_1,\ldots,g_n)$. The divided
difference in $i$ variables of a monic polynomial of degree $L$ is a
symmetric polynomial whose top homogeneous part is the complete
symmetric polynomial of degree $L-i+1$; as a polynomial in $x_i$ it is
therefore monic of degree $L-i+1$:
\[
 g_i=x_i^{L-i+1}
       +\sum_{r=0}^{L-i}a_{ir}(x_1,\ldots,x_{i-1})x_i^r.
\]
Successive monic division makes $R_n/J$ free on
\eqref{eq:sym-regular-finite-basis}, of rank $(L)_n$.
A divided difference of $P_L$ at distinct roots of $P_L$ vanishes, so
each $g_i$ vanishes on $\Omega_L$ and $J\subseteq K_L$.

Over $K$, the equations $g_1=\cdots=g_i=0$ successively select $i$
distinct roots of $P_L$: if $x_1,\ldots,x_{i-1}$ are distinct
roots of $P_L$, then
$P_L[x_1,\ldots,x_{i-1},X]=P_L(X)/\prod_{k<i}(X-x_k)$, whose roots are
the remaining roots of $P_L$. The common zeros of $J$ are therefore
exactly the $(L)_n$ points of $\Omega_L$, all defined over $K$. The
$K$-algebra $(R_n/J)\otimes_AK$ has dimension $(L)_n$, the number of
its points, so each of its local factors is one-dimensional and
\[
 (R_n/J)\otimes_A K\ \xrightarrow{\ \mathrm{ev}\ }\ K^{\Omega_L}
\]
is an isomorphism. Since $R_n/J$ is $A$-free, evaluation is injective
on $R_n/J$ itself, which means $J=K_L$.
\end{proof}

\begin{definition}\label{def:sym-regular-completion}
The regular nonsymmetric Newton completion is
\[
 \mathcal N_n^{\mathrm{reg},+}=\varprojlim_{L\ge n}R_n/K_L.
\]
\end{definition}

This ring carries a continuous action of $\mathfrak S_n$, and
evaluation at all nonnegative regular nodes is injective. Each $x_i$
is invertible in every finite quotient, since $P_L(x_i)\in K_L$ and
$P_L$ has a unit constant coefficient, so defining the completion with
Laurent polynomials gives the same ring.

\subsection{Symmetric invariants and freeness}

\begin{theorem}\label{thm:sym-integral-invariants}
There is a natural isomorphism of topological rings
\begin{equation}\label{eq:sym-invariant-isomorphism}
 (\mathcal N_n^{\mathrm{reg},+})^{\mathfrak S_n}
 \cong\mathcal B_n.
\end{equation}
At finite level, if $J_L=K_L\cap\Lambda_n$, then
\[
 J_L=\bigoplus_{\lambda_1>L-n}\A S_\lambda,
 \qquad (R_n/K_L)^{\mathfrak S_n}=\Lambda_n/J_L.
\]
The last quotient is free over $\A$ of rank $\binom Ln$.
\end{theorem}

\begin{proof}
Each regular orbit in $\Omega_L$ contains a unique point
$(Q^{a_1},\ldots,Q^{a_n})$ with $a_1>\cdots>a_n\ge0$. Writing
$a_i=\mu_i+n-i$ identifies these orbits with the partitions satisfying
$\mu_1\le L-n$, that is, with the nodes $p_\mu$ for $\mu$ in the
$n\times(L-n)$ rectangle. A symmetric polynomial vanishes on
$\Omega_L$ if and only if it vanishes at these nodes. The rectangle is
closed under passing to smaller partitions, so the argument in the
proof of Lemma~\ref{lem:sym-partition-vanishing} gives
\[
 J_L=\bigoplus_{\lambda_1>L-n}A S_\lambda.
\]
The remaining basis elements are indexed by the partitions in the
rectangle, so $\Lambda_n/J_L$ is free of rank $\binom Ln$.

To lift invariant classes integrally we use the type-$A$
divided-difference symmetrizer \cite{Demazure}. With
$\boldsymbol\delta=(n-1,n-2,\ldots,0)$ define
\[
 E(f)=\frac{\Alt(\mathbf x^{\boldsymbol\delta}f)}{\Delta(\mathbf x)},
 \qquad
 \Alt(f)=\sum_{w\in\mathfrak S_n}\sgn(w)wf.
\]
The alternating numerator is divisible by $\Delta$ over $A$, so
$E(f)\in\Lambda_n$, and
\[
 \Alt(\mathbf x^{\boldsymbol\delta})=\Delta,
 \qquad E(g)=g\quad(g\in\Lambda_n).
\]
Suppose that the class of $f$ in $R_n/K_L$ is invariant. Then
$(wf)(p)=f(p)$ for every $w$ and every $p\in\Omega_L$, that is, $f$ is
constant on each regular orbit, and
\begin{align*}
 \Delta(p)E(f)(p)
 &=\sum_w\sgn(w)(w\mathbf x^{\boldsymbol\delta})(p)(wf)(p)\\
 &=f(p)\sum_w\sgn(w)(w\mathbf x^{\boldsymbol\delta})(p)
 =\Delta(p)f(p).
\end{align*}
Since $\Delta(p)\ne0$ in the domain $A$, the symmetric polynomial
$E(f)$ represents the same class as $f$. For $f\in K_L$ the same
computation shows $E(f)\in J_L$. The natural map
$\Lambda_n/J_L\to(R_n/K_L)^{\mathfrak S_n}$ is therefore surjective,
and it is injective because $J_L=K_L\cap\Lambda_n$. This gives
\[
 (R_n/K_L)^{\mathfrak S_n}=\Lambda_n/J_L.
\]
Finally, taking invariants commutes with inverse limits, and the
rectangle cutoffs $\lambda_1>L-n$ and the size cutoffs $|\lambda|\ge d$
define the same topology on $\prod_\lambda\A S_\lambda$, because
$\lambda_1\le|\lambda|\le n\lambda_1$. Passing to the limit proves
\eqref{eq:sym-invariant-isomorphism}.
\end{proof}

\begin{theorem}\label{thm:sym-free-cover}
The regular nonsymmetric ring is a finite free algebra of rank $n!$
over its invariant ring:
\begin{equation}\label{eq:sym-artin-cover}
 \mathcal N_n^{\mathrm{reg},+}
 \cong R_n\otimes_{\Lambda_n}\mathcal B_n
 =\bigoplus_{0\le a_i\le n-i}
       \mathcal B_n\,x_1^{a_1}\cdots x_n^{a_n}.
\end{equation}
\end{theorem}

\begin{proof}
Let $\mathcal A$ be the set of Artin monomials in
\eqref{eq:sym-artin-cover}. Over $\Lambda_n$, the ordered roots of the
generic polynomial
\[
 X^n-e_1X^{n-1}+\cdots+(-1)^ne_n
\]
can be adjoined one at a time: $x_1$ satisfies this monic equation of
degree $n$, $x_2$ satisfies its quotient by $X-x_1$, of degree $n-1$,
and so on. The resulting algebra is $R_n$, and successive monic
division gives
\[
 R_n=\bigoplus_{a\in\mathcal A}\Lambda_n a,
 \qquad |\mathcal A|=n!.
\]
This is \cite[Proposition~1.10]{EkedahlLaksov} with the order of the
roots reversed.

For a regular orbit $p_\mu$, put
$\mathfrak j_\mu=(e_i-e_i(p_\mu):1\le i\le n)\subset\Lambda_n$, the
kernel of evaluation at $p_\mu$. The quotient
$R_n/\mathfrak j_\mu R_n=R_n\otimes_{\Lambda_n}A$ is the ordered
splitting algebra of the polynomial whose roots are the coordinates of
$p_\mu$, and it is free of rank $n!$ over $A$. Over $\Frac(A)$ these
roots are distinct, so evaluation at the $n!$ orderings of $p_\mu$ is
an isomorphism, by the same count of points as in the proof of
Theorem~\ref{thm:sym-regular-finite-presentation}. Over $A$ the kernel
of evaluation is therefore torsion, and it is zero because the source
is free. Thus $\mathfrak j_\mu R_n$ is exactly the ideal of
polynomials vanishing on the orbit of $p_\mu$.

Since $R_n$ is free over $\Lambda_n$ with basis $\mathcal A$, an ideal
$I\subseteq\Lambda_n$ satisfies $IR_n=\bigoplus_{a\in\mathcal A}Ia$,
so finite intersections of extended ideals are computed coefficient by
coefficient. The ideal $K_L$ is the intersection of the vanishing
ideals of the orbits contained in $\Omega_L$, and the intersection of
the corresponding $\mathfrak j_\mu$ is $J_L$. Consequently
\[
 \begin{split}
 K_L
 &=\bigcap_{\mu_1\le L-n}\mathfrak j_\mu R_n
 =\left(\bigcap_{\mu_1\le L-n}\mathfrak j_\mu\right)R_n
 =J_LR_n,\\
 R_n/K_L
 &=\bigoplus_{a\in\mathcal A}(\Lambda_n/J_L)a.
 \end{split}
\]
The finite direct sum commutes with the inverse limit, and the limit
of the $\Lambda_n/J_L$ is $\mathcal B_n$ by
Theorem~\ref{thm:sym-integral-invariants}. This gives
\eqref{eq:sym-artin-cover}.
\end{proof}

\subsection{Repeated coordinates}

Let $\mathcal N_n^{\mathrm{cart}}=\varprojlim_L R_n/T_L$ be the
ordinary Cartesian Newton completion, whose elements have expansions
$\sum_{\alpha\in\N^n}c_\alpha\prod_iP_{\alpha_i}(x_i;Q)$ with
$c_\alpha\in\A$. Unlike the regular completion, it also interpolates
at tuples with repeated coordinates.

\begin{proposition}\label{prop:sym-cartesian-injection}
Restriction to regular nodes defines continuous injections
\[
 \mathcal N_n^{\mathrm{cart}}\hookrightarrow
 \mathcal N_n^{\mathrm{reg},+},\qquad
 (\mathcal N_n^{\mathrm{cart}})^{\mathfrak S_n}
 \hookrightarrow\mathcal B_n.
\]
\end{proposition}

\begin{proof}
Continuity follows from $T_L\subseteq K_L$. Let $f$ be a Cartesian
Newton element vanishing at every regular node, and let
$\mathbf a\in\N^n$ be an exponent tuple, possibly with repetitions.
For all large integers $k$ the tuple $b_i=a_i+(i-1)k$ is regular. Let
$g$ be a Laurent polynomial representing $f$ on a box containing both
$\mathbf a$ and $\mathbf b$. Each coordinate difference
$Q^{b_i}-Q^{a_i}=Q^{a_i}(Q^{(i-1)k}-1)$ is divisible by $Q^k-1$, so
\[
 f(Q^{b_1},\ldots,Q^{b_n})
 \equiv f(Q^{a_1},\ldots,Q^{a_n})\pmod{Q^k-1}.
\]
The left side is zero, so the Laurent polynomial on the right is
divisible by $Q^k-1$ for arbitrarily large $k$, and it vanishes. All
Cartesian evaluations of $f$ are therefore zero, and $f=0$ by Newton
triangularity. The assertion about invariants follows from
Theorem~\ref{thm:sym-integral-invariants}.
\end{proof}

To compare the values at repeated coordinates we use the mixed
completion
\[
 \mathcal M_n=\varprojlim_k
 \frac{R_n}{((Q;Q)_k,P_k(x_1),\ldots,P_k(x_n))},\qquad
 \mathcal H=\varprojlim_k\A/((Q;Q)_k).
\]
Monic Newton division identifies $\mathcal M_n$ with the ring of
tensor Newton series with coefficients in $\mathcal H$. For each fixed
$k$, $S_\lambda$ lies in the mixed ideal above once $|\lambda|$ is
large enough \cite[Lemma~10.30]{BG}, so $\mathcal B_n$ maps
continuously to $\mathcal M_n$. This map is injective, because
evaluation at a partition node is continuous on $\mathcal M_n$, and on
the image of $F=\sum b_\lambda S_\lambda$ it gives the finite sum
$\sum_{\lambda\subseteq\mu}b_\lambda S_\lambda(p_\mu)\in\A$; since
$\A\hookrightarrow\mathcal H$, an element with zero image has zero
partition values, and Lemma~\ref{lem:sym-partition-vanishing} applies.

\begin{example}\label{ex:sym-singular-boundary}
Already for $n=2$ the inclusion of Cartesian invariants is proper.
Here
\[
 S_{(k,0)}(x,y;Q)
 =\frac{P_{k+1}(x;Q)-P_{k+1}(y;Q)}{x-y},\qquad
 S_{(k,0)}(1,1;Q)=(Q;Q)_k\ne0.
\]
Thus $S_{(k,0)}\to0$ in the partition topology but not in the
Cartesian topology, in which evaluation at $(1,1)$ is continuous with
values in the discrete ring $\A$. The element
$F=\sum_{k\ge0}(-1)^kS_{(k,0)}$ of $\mathcal B_2$ does not come from
the Cartesian completion. Its value at the repeated node is
\[
 F(1,1)=\sum_{k\ge0}\prod_{j=1}^k(Q^j-1)\in\mathcal H.
\]
At $Q=1+h$ every factor $Q^j-1$ has nonnegative coefficients and
leading term $jh$, so all coefficients of $F(1,1)$ are nonnegative,
and the coefficient of $h^r$ is at least $r!$, already from the term
$k=r$. The Taylor coefficients at $Q=1$ of a Laurent polynomial grow at
most polynomially: positive powers of $Q$ contribute finitely many
terms, and $Q^{-a}$ contributes binomial coefficients of fixed degree.
So $F(1,1)\notin\A$, whereas a Cartesian Newton element has its value
at $(1,1)$ in $\A$.
\end{example}

The symmetric completion therefore allows Habiro-valued values at
repeated coordinates, while the Cartesian completion has Laurent
polynomial values there.

\subsection{Knot elements}

The parameter denoted $q$ in \cite{BG} is $Q$ here.
For a zero-framed knot $K$,
\cite[Theorems~1.3 and~8.2, Proposition~8.7]{BG} and
\eqref{eq:sym-unit-conversion} provide an element
\begin{equation}\label{eq:sym-knot-element}
 Z_{K,n}=\sum_\lambda b_{K,n,\lambda}(Q)S_\lambda
 \in\mathcal B_n,\qquad b_{K,n,\lambda}(Q)\in\A.
\end{equation}
In the proof of Theorem~8.2 of that paper, the eigenvalue of the
central basis element $\sigma_\lambda$ on $V(\mu)$ is identified with
$F_\lambda^{\mathrm{BG}}(p_\mu;Q^{-1})$. Hence $Z_{K,n}(p_\mu)$ is the
reduced invariant of $K$ colored by $V(\mu)$, and Proposition~8.7 of
that paper shows that the coefficients in \eqref{eq:sym-knot-element}
are Laurent polynomials.

\begin{theorem}\label{thm:sym-dominant-knot-newton}
For the fundamental weights $\omega_1,\ldots,\omega_{n-1}$ of
$\mathfrak{sl}_n$, the substitution
\begin{equation}\label{eq:sym-dynkin-substitution}
 x_i=Q^{n-i}\prod_{j=i}^{n-1}t_j\quad(i<n),\qquad x_n=1
\end{equation}
extends continuously to a ring homomorphism
$\mathcal B_n\to\mathcal N^+_{n-1}$, where
$\mathcal N^+_{n-1}$ is the dominant Newton ring of
Section~\ref{sec:newton} in the variables $t_1,\ldots,t_{n-1}$,
with nodes $(Q^{r_1},\ldots,Q^{r_{n-1}})$, $r_i\ge0$.
The image of $Z_{K,n}$ is the unique integral Newton element whose
values are the reduced $\mathfrak{sl}_n$ invariants of $K$ in the
dominant highest weights $\sum_i r_i\omega_i$.
\end{theorem}

\begin{proof}
At $t_j=Q^{r_j}$ the substitution gives $p_\mu$, where
\[
 \mu_i=\sum_{j=i}^{n-1}r_j,\qquad \mu_n=0.
\]
On the box $0\le r_j<M$ one has
\[
 |\mu|=\sum_{j=1}^{n-1}j r_j
       \le(M-1)\frac{n(n-1)}2=:d_M.
\]
If $|\lambda|>d_M$, the substituted $S_\lambda$ is a polynomial in
$\mathbf t$ that vanishes at every node of this box by
Lemma~\ref{lem:sym-partition-vanishing}, so it lies in the box ideal by
Theorem~\ref{thm:newton-interpolation}. The substitution therefore maps
$I_{d_M+1}$ into the box ideal of level $(M,\ldots,M)$, is continuous,
and extends to a ring homomorphism over $A$.

At every dominant node, the image of $Z_{K,n}$ takes the value
$Z_{K,n}(p_\mu)$, the $\mathfrak{gl}_n$ invariant with color $\mu$. For
a zero-framed knot this equals the $\mathfrak{sl}_n$ invariant with
color $\sum_ir_i\omega_i$ by \cite[Corollary~6.5]{BG}. Uniqueness
follows from injectivity of dominant-grid evaluation.
\end{proof}

The substitution is not injective, since
$S_{(1^n)}=\prod_i(x_i-1)$ restricts to zero. It is injective,
however, on a subring that contains the knot elements.

For $L\ge0$, let $K_L^{\mathbb Z}$ be the ideal of Laurent
polynomials vanishing at all nodes $Q^{\mathbf a}$ with
$-L\le a_i\le L$ and pairwise distinct exponents. Write
\[
 \mathcal N_n^{\mathrm{reg},\mathbb Z}
 =\varprojlim_L
 A[x_1^{\pm1},\ldots,x_n^{\pm1}]/K_L^{\mathbb Z}
\]
for the completion on the integer regular grid.

\begin{proposition}\label{prop:sym-determinant-extension}
The operation $\tau f(\mathbf x)=f(Qx_1,\ldots,Qx_n)$ is continuous
on $\mathcal B_n$. Its fixed subring injects under
\eqref{eq:sym-dynkin-substitution}. Moreover, every $\tau$-fixed
element extends uniquely to an element of
$\mathcal N_n^{\mathrm{reg},\mathbb Z}$ invariant under Weyl
permutations and common translations. In particular, the element
$Z_{K,n}$ of a
zero-framed knot is fixed by $\tau$ and has this extension.
\end{proposition}

\begin{proof}
On partition nodes,
\[
 (\tau f)(p_\mu)=f(p_{\mu+(1^n)}).
\]
If $f$ vanishes at all $p_\mu$ with $|\mu|<d$, then $\tau f$ vanishes
at all $p_\mu$ with $|\mu|<d-n$, so $\tau$ is continuous. If
$\tau f=f$, then
\[
 f(p_\mu)=f(p_{\mu-\mu_n(1^n)}),
\]
so $f$ is determined by its values on the slice $\mu_n=0$. These are
exactly the nodes reached by the Dynkin substitution, which proves the
injectivity.

For an integer regular tuple $\mathbf a$, choose $c$ with
$a_i+c\ge0$ for all $i$ and define
\[
 \widetilde f(Q^{a_1},\ldots,Q^{a_n})
       =f(Q^{a_1+c},\ldots,Q^{a_n+c}),
\]
where $f$ is evaluated on the nonnegative regular grid through
Theorem~\ref{thm:sym-integral-invariants}. The fixed-point condition
makes this independent of $c$. To obtain a Laurent representative on a
finite integer box, choose one such $c$ for the whole box and a
partition truncation $f_N$ representing $f$ at all shifted nodes; then
$f_N(Q^cx_1,\ldots,Q^cx_n)$ represents $\widetilde f$ on the box.
These classes are compatible, so they define an element of the
integer regular completion. It is invariant under Weyl permutations
and common translations by construction, and conversely these
invariances determine every integer value from a nonnegative one,
which gives uniqueness.

Tensoring a $\mathfrak{gl}_n$ color with the determinant does not
change the underlying $\mathfrak{sl}_n$ color. For a zero-framed knot,
\cite[Corollary~6.5]{BG} therefore gives
\[
 Z_{K,n}(p_{\mu+(1^n)})=Z_{K,n}(p_\mu).
\]
Injectivity of partition evaluation yields $\tau Z_{K,n}=Z_{K,n}$,
so the construction applies to the knot element.
\end{proof}

\begin{corollary}\label{cor:sym-shifted-weyl-knot}
The knot element extends to $A$-valued functions on the regular grid,
represented by Laurent polynomials on every finite subset of
\[
 \mathcal C_{\mathrm{reg}}=
 \left\{\mathbf r\in\Z^{n-1}:
   \sum_{j=a}^{b-1}(r_j+1)\ne0
   \text{ for every }1\le a<b\le n\right\}.
\]
It is invariant under the shifted Weyl action. In Dynkin variables its
simple reflection sends $t_i$ to $Q^{-2}t_i^{-1}$, sends
$t_j$ to $Qt_it_j$ for vertices $j$ adjacent to $i$, and fixes the other
variables.
\end{corollary}

\begin{proof}
For integer Dynkin labels the exponents in
\eqref{eq:sym-dynkin-substitution} are
$a_i=n-i+\sum_{j=i}^{n-1}r_j$ and $a_n=0$,
and their differences are $a_a-a_b=\sum_{j=a}^{b-1}(r_j+1)$.
The displayed set is therefore the slice $a_n=0$ of the integer
regular grid, and we restrict the extension of
Proposition~\ref{prop:sym-determinant-extension} to it. Under the
monomial substitution \eqref{eq:sym-dynkin-substitution}, Laurent
representatives on finite sets remain Laurent polynomials. A Weyl
permutation of the exponents, followed by the common translation that
restores $a_n=0$, is the shifted Weyl action on Dynkin labels, and the
extension is invariant under both operations. For the transposition of
$a_i$ and $a_{i+1}$ the new labels are $r_i'=-r_i-2$ and
$r_j'=r_j+r_i+1$ for $j=i\pm1$, which gives the stated formulas.
\end{proof}

\begin{corollary}\label{cor:sym-knot-consequences}
Let $\mathcal J_{K,n}$ be the dominant Newton element. All mixed
divided differences of its colored values lie in $A$; in particular,
for two distinct Dynkin directions,
\[
 J_{K,\omega_i+\omega_j}^{(n)}
 -J_{K,\omega_i}^{(n)}-J_{K,\omega_j}^{(n)}+1
 \in (Q-1)^2A.
\]
Along each simple singular wall the knot element also has a canonical
Habiro-valued boundary value, with normal jets as in
Theorem~\ref{thm:newton-wall-jets}.
\end{corollary}

\begin{proof}
The integrality of divided differences is
Proposition~\ref{prop:newton-divided-differences}, and the congruence is
the case of the index $\mathbf e_i+\mathbf e_j$, with $J^{(n)}_{K,0}=1$.
The statement about walls is Theorem~\ref{thm:newton-wall-jets} applied
to this element.
\end{proof}

\subsection{Comparison with symmetric-color Newton rings}
\label{subsec:sym-comparison}

The restriction of $\mathcal J_{K,n}$ to the first Dynkin axis
interpolates the symmetric-power colors. We compare this restriction
with the one-variable Newton ring of \cite{FangNewton}, which uses a
trace coordinate in place of $t=Q^r$.

Let $R=\mathbb Z[q^{\pm1}]$, with $Q=q^2$, and set
\[
 x_r^{(n)}=q^{n+2r}+q^{-n-2r},\qquad
 B_k^{(n)}(X)=\prod_{j=0}^{k-1}(X-x_j^{(n)}).
\]
The symmetric-color Newton ring is
\begin{equation}\label{eq:sym-trace-newton-ring}
 \mathcal C_n=\varprojlim_k R[X]/(B_k^{(n)}(X)).
\end{equation}
We write $\mathcal N_1(R)$ for the dominant Newton completion with
coefficients in $R$. The substitution
\begin{equation}\label{eq:sym-trace-substitution}
 X=q^nt+q^{-n}t^{-1}
\end{equation}
takes $t=Q^r$ to $x_r^{(n)}$. At the level of Laurent polynomials,
the involution $t\mapsto q^{-2n}t^{-1}$ has invariant ring $R[X]$:
with $z=q^nt$ it becomes $z\mapsto z^{-1}$, every invariant Laurent
polynomial is an $R$-linear combination of $1$ and the $z^k+z^{-k}$,
and these are polynomials in $X=z+z^{-1}$.

\begin{proposition}\label{prop:sym-trace-embedding}
For $n\ge2$, substitution \eqref{eq:sym-trace-substitution} induces
a continuous injective homomorphism
\begin{equation}\label{eq:sym-trace-embedding}
 \iota_n:\mathcal C_n\hookrightarrow\mathcal N_1(R).
\end{equation}
Its image is a proper subring of $\mathcal N_1(R)$.
\end{proposition}

\begin{proof}
Each trace-coordinate Newton factor satisfies
\begin{equation}\label{eq:sym-trace-factorization}
 X-x_j^{(n)}
   =q^n(t-Q^j)
       \bigl(1-q^{-2n}Q^{-j}t^{-1}\bigr).
\end{equation}
Consequently $B_k^{(n)}(X(t))$ belongs to $P_k(t;Q)R[t^{\pm1}]$, so
the substitution is continuous for the Newton topologies and extends
to the completions. If $\iota_n(F)=0$, evaluation at every $t=Q^r$
gives $F(x_r^{(n)};q)=0$. The trace nodes $x_r^{(n)}$ are distinct in
the domain $R$, so the triangular Newton evaluation matrices give
$F=0$.

For any $F\in\mathcal C_n$, the difference of its first two values is
divisible in $R$ by
\[
 x_1^{(n)}-x_0^{(n)}
   =q^n(Q-1)(1-q^{-2n}Q^{-1}).
\]
For the element $t\in\mathcal N_1(R)$ this difference is $Q-1$. The
factor $1-q^{-2n}Q^{-1}$ is not a unit of $R$, so $Q-1$ is not
divisible by $x_1^{(n)}-x_0^{(n)}$, and $t$ is not in the image.
\end{proof}

The two rings thus impose different integral congruences on a
one-parameter color family. For knot invariants, however, the two
constructions give the same function on their common colors.

\begin{corollary}\label{cor:sym-symmetric-color-restriction}
Let $F_{K,n}\in\mathcal C_n$ be the symmetric-power knot element
of \cite{FangNewton}.  With the same reduced, zero-framed and
quantum-parameter conventions, the knot element of
Theorem~\ref{thm:sym-dominant-knot-newton} satisfies, over $R$,
\begin{equation}\label{eq:sym-knot-axis-comparison}
 \mathcal J_{K,n}(t,1,\ldots,1;q^2)
   =\iota_n(F_{K,n})
   =F_{K,n}(q^nt+q^{-n}t^{-1};q).
\end{equation}
\end{corollary}

\begin{proof}
Setting $t_2=\cdots=t_{n-1}=1$ in the tensor Newton expansion kills
every term with a positive index in one of these coordinates, because
$P_k(1)=0$ for $k\ge1$. It therefore defines an element of
$\mathcal N_1(R)$. At $t=Q^r$ both sides of
\eqref{eq:sym-knot-axis-comparison} equal $J_{K,r\omega_1}^{(n)}(Q)$,
and injectivity of Newton evaluation gives the identity.
\end{proof}

At $q=1$, and hence $Q=1$, the trace coordinate becomes
$X=t+t^{-1}$. Write
$\Delta_K(t)=\delta_K(t+t^{-1}-2)$, where
$\delta_K\in\mathbb Z[v]$ and $\delta_K(0)=1$. The formula at $Q=1$ in
Theorem~\ref{main:type-a}, restricted to the first Dynkin axis, reads
\[
 g_{K,n,1}(t,1,\ldots,1;0)=\Delta_K(t)^{-(n-1)},
\]
and in the trace coordinate it becomes
$F_{K,n}(2+v;1)=\delta_K(v)^{-(n-1)}$, the classical specialization
of \cite{FangNewton}.

For the figure-eight knot and $n=2$, all trace-coordinate Newton
coefficients are equal to one, and substitution in
\eqref{eq:sym-knot-axis-comparison} gives the series
\eqref{eq:ex-figure-eight-series}. At $Q=1$ its two coordinate
descriptions are
\[
 F_{4_1,2}(2+v;1)=\frac1{1-v},\qquad
 g_{4_1,2,1}(1+u,0)=\frac{1+u}{1+u-u^2},
 \qquad v=\frac{u^2}{1+u}.
\]
So the trace Newton coefficients can be constant while the
coefficients in the basis $P_k(t;Q)$ depend on $k$ and $Q$.

The same substitution relates the singular boundaries.  At
$t=Q^{-1}$ the trace coordinate takes the value
\[
 X_*^{(n)}=q^{n-2}+q^{2-n}=x_{-1}^{(n)},
\]
which is the first node removed from the trace grid. By
\eqref{eq:sym-trace-factorization} and
$P_k(Q^{-1};Q)=Q^{-k}(Q;Q)_k$, the value $B_k^{(n)}(X_*^{(n)})$ is
divisible by $(Q;Q)_k$ over $R$. Both boundary evaluations therefore
converge in
\[
 \mathcal H_R=\varprojlim_k R/((Q;Q)_k),
\]
and they agree under $\iota_n$ by continuity. This identifies the
missing-node boundary of \cite{FangNewton} with the boundary of
Section~\ref{subsec:newton-walls}, restricted to the trace subring.

The normal parameters are related by the same substitution. For
$\mathcal C_2$, put $u=t-Q^{-1}$; then
\[
 X-2=Qt+Q^{-1}t^{-1}-2
     =\frac{Q^2u^2}{1+Qu}.
\]
The trace parameter thus has order two in the normal parameter $u$:
an expansion modulo $(X-2)^h$ pulls back to an expansion modulo
$u^{2h}$, and the higher boundary coefficients are related by this
ramified change of coordinate.

\section{Relative Habiro rings and the Newton image}
\label{sec:relative}

This section proves Theorem~\ref{main:recognition}. We first
characterize Habiro's ring by its Taylor expansions at roots of unity,
then define relative Habiro collections and embed them into the Newton
ring by interpolating their color values, and finally describe the
image, both through rational germs and through the Newton coefficients.

\subsection{Scalar Taylor series}

Habiro's cyclotomic completion has injective Taylor maps at
roots of unity \cite[Theorem~5.2]{Habiro}.  The following theorem
characterizes their simultaneous image by integral gluing conditions.

\begin{theorem}\label{thm:rel-scalar-gluing}
Let $\mathcal H=\varprojlim_N A/((Q;Q)_N)$.  Suppose that for every
root of unity $\zeta$ a series
\[
 f_\zeta(x)\in\mathbb Z[\zeta][[x]]
\]
is given. Assume that the family is Galois-equivariant and that,
whenever $\zeta$ and $\xi$ have orders $m$ and $pm$ for a prime $p$
and are congruent in the residue field of a common complete
cyclotomic discrete valuation ring $\mathcal O$ above $p$,
\begin{equation}\label{eq:rel-scalar-glue}
 f_\zeta(x+\xi-\zeta)=f_\xi(x)
 \quad\hbox{in }\mathcal O[[x]].
\end{equation}
Here the left side is expanded coefficientwise in the valuation
topology. Then there is a unique $f\in\mathcal H$ with the prescribed
Taylor series at every root of unity. Conversely, every element of
$\mathcal H$ defines such a family, and the Taylor homomorphism at
$Q=1$ is injective.
\end{theorem}

\begin{proof}
We first interpolate the prescribed jets modulo a finite product of
cyclotomic polynomials. Let
\[
 \delta(Q)=\prod_{m\in S}\Phi_m(Q)^{e_m},\qquad e_m\ge1,
\]
with $S$ finite. Hermite interpolation over a cyclotomic field gives a
unique $b_\delta\in\overline{\Q}[Q]$ of degree less than $\deg\delta$
whose Taylor expansion at each primitive $m$th root, $m\in S$, agrees
with $f_\zeta$ through order $e_m-1$. By uniqueness and Galois
equivariance, $b_\delta\in\Q[Q]$.

To show that $b_\delta$ has integral coefficients, fix a prime $p$ and
a complete cyclotomic discrete valuation ring $\mathcal O$ above $p$
that contains all roots of unity of orders in $S$; let $\mathfrak m$
be its maximal ideal. Write each root as
\[
 \zeta=\zeta_0\eta,\qquad
 (\ord\zeta_0,p)=1,\qquad \ord\eta=p^a.
\]
Since $\eta\equiv1\pmod{\mathfrak m}$ and distinct roots of order
prime to $p$ stay distinct modulo $\mathfrak m$, two roots have the
same reduction exactly when their prime-to-$p$ parts $\zeta_0$ agree.
We call the set of roots with a given $\zeta_0$ a cluster. Applying
\eqref{eq:rel-scalar-glue} successively along the orders
$\ord(\zeta_0)p^j$, $j=0,1,\ldots,a$, shows that every series in the
cluster is a translate of the one at $\zeta_0$:
$f_\zeta(x)=f_{\zeta_0}(x+\zeta-\zeta_0)$.

For one cluster put
\[
 W(y)=\prod_{\zeta\text{ in the cluster}}
       \bigl(y-(\zeta-\zeta_0)\bigr)^{e_{\ord\zeta}},
 \qquad C=\mathcal O[y]/(W),\qquad s=\deg W.
\]
The roots $\zeta-\zeta_0$ of $W$ lie in $\mathfrak m$, so every
nonleading coefficient of $W$ lies in $\mathfrak m$ and
$y^s\in\mathfrak m C$. Thus $y$ is topologically nilpotent in the
finite free complete $\mathcal O$-algebra $C$, and $f_{\zeta_0}(y)$
has a well-defined image in $C$. Its unique representative of degree
less than $s$ has coefficients in $\mathcal O$. Under the map
\[
 C\longrightarrow\mathcal O[x]/(x^{e_{\ord\zeta}}),
 \qquad y\longmapsto\zeta-\zeta_0+x,
\]
this representative goes to the jet of
$f_{\zeta_0}(x+\zeta-\zeta_0)=f_\zeta(x)$, as required.

Distinct clusters have distinct reductions, so the polynomials
$W(Q-\zeta_0)$ attached to different clusters have unit resultants,
and the Chinese remainder theorem combines the representatives over
$\mathcal O$ into a polynomial of degree less than $\deg\delta$ with
all the prescribed jets. By uniqueness over the fraction field this
polynomial is $b_\delta$. This holds at every prime, whence
\[
 v_p([Q^j]b_\delta)\ge0\quad\text{for every prime }p\text{ and every }j,\qquad b_\delta\in\Z[Q].
\]

If $\delta\mid\delta'$, then $b_{\delta'}$ reduced modulo the monic
polynomial $\delta$ has the jets that characterize $b_\delta$, so
\[
 b_{\delta'}\equiv b_\delta\pmod\delta
       \quad\text{in }\Z[Q].
\]
We apply this to the monic associates of the $(Q;Q)_N$. Their constant
terms are units, so adjoining $Q^{-1}$ does not change the quotients,
and the compatible remainders define an element of $\mathcal H$ with
the prescribed Taylor series. Since each remainder is determined by
the jets, this element is unique, which also proves injectivity of the
full Taylor map.

Conversely, compatible representatives of $f\in\mathcal H$ give,
by telescoping,
\[
 f=\sum_{N\ge0}a_N(Q)(Q;Q)_N,\qquad a_N\in A.
\]
At an $m$th root the $N$th term vanishes to order at least
$\lfloor N/m\rfloor$, so the sum has a Taylor expansion there. For two
roots $\zeta,\xi$ as in \eqref{eq:rel-scalar-glue}, put $Q=\xi+x$. The
factors $1-Q^j$ with $m\mid j$ lie in $(\mathfrak m,x)$, because
$\xi^j$ is then a $p$th root of unity and hence congruent to $1$. The
$N$th term therefore lies in $(\mathfrak m,x)^{\lfloor N/m\rfloor}$,
and the sum converges in $\mathcal O[[x]]$. This single convergent sum
is both the Taylor expansion at $\xi$ and the re-expansion at $\xi$ of
the Taylor expansion at $\zeta$, which proves
\eqref{eq:rel-scalar-glue}. Galois equivariance holds termwise.

Finally, suppose the Taylor series of $f$ at $1$ is zero. We show by
induction on the order that all Taylor series of $f$ vanish. Let
$\zeta$ have order $m>1$, choose a prime $p\mid m$, and write
$\zeta=\zeta_0\eta$ as above with $a\ge1$. The root $\zeta'=\zeta_0\eta^p$
has order $m/p$ and the same reduction, so $f_{\zeta'}=0$ by induction,
and \eqref{eq:rel-scalar-glue} gives $f_\zeta=0$. Injectivity of the
full Taylor map then gives $f=0$.
\end{proof}

\subsection{Relative collections}

Put
\[
 B=\Z[T_1^{\pm1},\ldots,T_d^{\pm1}],\qquad
 D\in B,\qquad D(\mathbf1)=1,\qquad R=B[D^{-1}].
\]
For a primitive root $\zeta$ of order $m$, write
\[
 R_\zeta=\Z[\zeta,T_1^{\pm1/m},\ldots,T_d^{\pm1/m},D(\mathbf T)^{-1}].
\]
We write $D(\mathbf t^m)$ for $D(t_1^m,\ldots,t_d^m)$, and we use the
multivariable version of the relative Habiro ring of
\cite[Definition~1.1 and equation~(6)]{GW}. In the component at a root
of order $m$, the variable $t_i=T_i^{1/m}$ is the color variable: at
a dominant color it is specialized to $Q^{r_i}$.

\begin{definition}\label{def:rel-collection}
A relative Habiro collection for $R/B$ is a family
\[
 H_\zeta(\mathbf T,x)\in R_\zeta[[x]],\qquad x=Q-\zeta,
\]
indexed by the roots of unity, with the following properties.

The family is Galois-equivariant: conjugating $\zeta$ applies the
same conjugation to its coefficients and fixes the formal variables
$x,T_i^{1/m}$.

For the compatibility condition, let $\zeta$ have order $m$, let
$\xi$ have order $pm$ for a prime $p$, and suppose the two roots
are $p$-adically close. Choose a common complete cyclotomic discrete
valuation ring $\mathcal O$ above $p$ and put $t_i=T_i^{1/m}$.
The coefficient map $\phi_p$ fixes $\mathcal O$ and $x$ and sends
\[
 T_i^{1/(pm)}\longmapsto T_i^{1/m},
 \qquad\text{hence}\qquad T_i\longmapsto T_i^p.
\]
We require
\begin{equation}\label{eq:rel-glue}
 H_\zeta(\mathbf T,x+\xi-\zeta)
   =\phi_p\bigl(H_\xi(\mathbf T,x)\bigr)
\end{equation}
in the ring
\[
 \widehat{\mathcal O[\mathbf t^{\pm1},D(\mathbf t^m)^{-1}]}_p[[x]].
\]
The denominator $D(\mathbf t^{pm})$ of the right side is invertible
in this ring, as shown below. These collections form a ring under
componentwise operations, denoted $\mathcal H_{\mathrm{rel}}(D)$.
\end{definition}

In terms of the color variables, $\phi_p$ identifies the color
variable $T_i^{1/(pm)}$ of the component at $\xi$ with the color
variable $t_i=T_i^{1/m}$ of the component at $\zeta$. The condition
\eqref{eq:rel-glue} therefore says that the two
components agree, as functions of the same color $\mathbf t$, after
re-expansion of the first one at $\xi$.

The inverse of $D(\mathbf T^p)$ exists in the completion: writing
\[
 E_p(\mathbf T)=\frac{D(\mathbf T^p)-D(\mathbf T)^p}{p}\in B,
\]
one has, for $a\ge1$,
\[
 D(\mathbf T^p)^{-a}
 =\sum_{j\ge0}(-p)^j\binom{a+j-1}{j}
       \frac{E_p(\mathbf T)^j}{D(\mathbf T)^{p(a+j)}},
\]
and the series converges $p$-adically.

Three lemmas on color expansions will be used repeatedly. The first
says that a formal germ is determined by its restrictions to the
colored curves of one residue class.

\begin{lemma}\label{lem:rel-color-rigidity}
Let $\Bbbk$ be a characteristic-zero field containing a primitive $m$th
root $\zeta$. Fix $b_1,\ldots,b_d\in\Z$. If
\[
 F\in\Bbbk[[t_1-\zeta^{b_1},\ldots,t_d-\zeta^{b_d},\epsilon]]
\]
vanishes after each substitution
\[
 t_i=\bigl(\zeta(1+\epsilon)\bigr)^{ma_i+b_i},
 \qquad (a_1,\ldots,a_d)\in\N^d,
\]
then $F=0$.
\end{lemma}

\begin{proof}
Set $w_i=t_i/\zeta^{b_i}-1$. If $F\ne0$, let
$P_\ell(\mathbf w,\epsilon)$ be its lowest nonzero homogeneous
component for the total degree. Under the substitution,
$w_i=(1+\epsilon)^{ma_i+b_i}-1=(ma_i+b_i)\epsilon+O(\epsilon^2)$, so the
coefficient of $\epsilon^\ell$ is
\[
 P_\ell(ma_1+b_1,\ldots,ma_d+b_d,1).
\]
A polynomial vanishing on a product of infinite subsets of a field is
zero, so $P_\ell(\mathbf w,1)=0$. Since the homogeneous components of
$P_\ell(\mathbf w,1)$ in $\mathbf w$ are the coefficients of the
powers of $\epsilon$ in $P_\ell$, this forces $P_\ell=0$, a
contradiction.
\end{proof}

The second lemma says that a rational function with integral Taylor
coefficients has an integral numerator.

\begin{lemma}\label{lem:rel-integral-saturation}
Let $\mathcal O\subset\Bbbk$ be a subring of a field, and let
$D_0\in\mathcal O[\mathbf t^{\pm1}]$ satisfy $D_0(\mathbf1)=1$.
Inside $\Bbbk[[\mathbf t-\mathbf1]]$ one has
\[
 \Bbbk[\mathbf t^{\pm1},D_0^{-1}]
 \cap\mathcal O[[\mathbf t-\mathbf1]]
 =\mathcal O[\mathbf t^{\pm1},D_0^{-1}].
\]
\end{lemma}

\begin{proof}
Write an element of the intersection as $P/D_0^e$ with
$P\in\Bbbk[\mathbf t^{\pm1}]$. Multiplication by $D_0^e$, whose Taylor
coefficients lie in $\mathcal O$, shows that $P$ has Taylor
coefficients in $\mathcal O$, and so does $\mathbf t^{\mathbf c}P$
for $\mathbf c\in\N^d$. Choosing $\mathbf c$ so that
$\mathbf t^{\mathbf c}P$ is a polynomial, its Taylor expansion
$\mathbf t^{\mathbf c}P(\mathbf1+\mathbf w)$ is a polynomial with
coefficients in $\mathcal O$, and the inverse translation
$\mathbf w=\mathbf t-\mathbf1$ shows that $\mathbf t^{\mathbf c}P$
has coefficients in $\mathcal O$. Hence
$P/D_0^e\in\mathcal O[\mathbf t^{\pm1},D_0^{-1}]$. The reverse
inclusion holds because $D_0$ is a unit in
$\mathcal O[[\mathbf t-\mathbf1]]$.
\end{proof}

The third lemma is the $p$-adically completed version of the
injectivity of color Taylor expansion.

\begin{lemma}\label{lem:rel-completed-injection}
Let $\mathcal O$ be a complete discrete valuation ring with uniformizer
$\varpi$. If $D_1,\ldots,D_s\in\mathcal O[\mathbf t^{\pm1}]$ have unit
value at $\mathbf t=\mathbf1$, then color Taylor expansion gives an
injection
\[
 \widehat{\mathcal O[\mathbf t^{\pm1},D_1^{-1},\ldots,D_s^{-1}]}_{\varpi}
 \lhook\joinrel\longrightarrow\mathcal O[[\mathbf t-\mathbf1]].
\]
\end{lemma}

\begin{proof}
Put $R=\mathcal O[\mathbf t^{\pm1},D_1^{-1},\ldots,D_s^{-1}]$
and $\mathbf w=\mathbf t-\mathbf1$. For every $v\ge1$, each
$D_i(\mathbf1+\mathbf w)$ is a unit of
$(\mathcal O/\varpi^v)[[\mathbf w]]$, so Taylor expansion is defined
on $R/\varpi^vR$. Suppose an element of $R/\varpi^vR$ has zero Taylor
expansion. Multiplying by its localization denominator and by a
monomial reduces this to a polynomial $P$ over $\mathcal O/\varpi^v$
with $P(\mathbf1+\mathbf w)=0$. Translation is an automorphism of the
polynomial ring over $\mathcal O/\varpi^v$, so $P=0$. The Taylor map
is therefore injective at every finite level, and taking inverse limits,
with
\[
 \varprojlim_v(\mathcal O/\varpi^v)[[\mathbf w]]
       =\mathcal O[[\mathbf w]],
\]
gives the assertion.
\end{proof}

At a primitive $m$th root $\zeta$ we write $t_i=T_i^{1/m}$ and
$x=Q-\zeta$, so that a relative component becomes
\begin{equation}\label{eq:rel-local}
 g_\zeta(\mathbf t,x)
   \in\mathbb Z[\zeta,\mathbf t^{\pm1},D(\mathbf t^m)^{-1}]
                 [[x]].
\end{equation}
By Theorem~\ref{thm:rel-scalar-gluing}, integral Galois-equivariant
families satisfying the scalar gluing identities are the same as
elements of $\Habi$. To apply this to a relative collection, we
specialize at a fixed $\mathbf r\in\mathbb N^d$ by
\begin{equation}\label{eq:rel-color-evaluation}
 t_i=(\zeta+x)^{r_i},\qquad
 T_i=(\zeta+x)^{m r_i}.
\end{equation}
At $x=0$ this gives $T_i=1$, and $D(\mathbf1)=1$, so every
denominator becomes a unit and the specialized series has coefficients
in $\Z[\zeta]$. If $\xi$ has order $pm$ and is $p$-adically close to
$\zeta$, the relative relation reads
\[
 H_\zeta(\mathbf T,x+\xi-\zeta)
        =\phi_pH_\xi(\mathbf T,x),
 \qquad
 \phi_p(T_i^{1/(pm)})=T_i^{1/m}.
\]
Substituting $T_i^{1/m}=(\xi+x)^{r_i}$ turns the left side into
the re-expansion at $\xi$ of the specialization at $\zeta$, and the
right side into the specialization at $\xi$. Relative gluing therefore
implies scalar gluing, and the specializations define an element of
$\Habi$, which we denote by $f_{\mathbf r}$.

\subsection{Newton interpolation of relative collections}

The values $f_{\mathbf r}\in\mathcal H$ of a relative collection are
interpolated by a Newton series with coefficients in $\mathcal H$. Its
coefficients are divided differences, and the main point is that the
divisions can be carried out in $\mathcal H$. Divisibility in
$\mathcal H$ by products of cyclotomic polynomials is controlled by the
following lemma.

\begin{lemma}\label{lem:rel-cyclotomic-quotient}
Suppose $\delta\in\A$ is, up to a Laurent unit, a finite product
of cyclotomic polynomials with nonnegative integer multiplicities.
Then there is a natural exact sequence
\[
 0\longrightarrow\Habi\xrightarrow{\ \delta\ }\Habi
       \longrightarrow\A/(\delta)\longrightarrow0.
\]
Consequently
\[
 \Habi/\delta\Habi\simeq\A/(\delta),
 \qquad \A\cap\delta\Habi=\delta\A,
\]
and $\delta$ is a non-zero-divisor in $\Habi$. Moreover, the abelian
group $\A/(\delta)$ is torsion-free.
\end{lemma}

\begin{proof}
For sufficiently large $N$, $\delta\mid c_N$, and the sequence
\[
 0\longrightarrow A/(c_N/\delta)
 \xrightarrow{\ \delta\ }A/(c_N)
 \longrightarrow A/(\delta)\longrightarrow0
\]
is exact. The left-hand system has surjective transition maps, and the
ideals $(c_N/\delta)$ define the Habiro topology: on the one hand
$(c_N)\subseteq(c_N/\delta)$, and on the other hand, given $M$, we may
choose $N$ such that
\[
 \left\lfloor\frac N\ell\right\rfloor
 \ge v_{\Phi_\ell}(\delta)+\left\lfloor\frac M\ell\right\rfloor
\]
for every $\ell$ occurring in $\delta c_M$; then
$(c_N/\delta)\subseteq(c_M)$. Exactness of the inverse limit gives
\[
 0\longrightarrow\mathcal H\xrightarrow{\ \delta\ }\mathcal H
 \longrightarrow A/(\delta)\longrightarrow0.
\]

Up to a Laurent unit, $\delta$ is monic with constant coefficient
$\pm1$. Hence $Q$ is invertible in $\Z[Q]/(\delta)$, so
$A/(\delta)=\Z[Q]/(\delta)$, which is free over $\Z$ on its usual
monomial basis (and zero when $\delta$ is a unit). The restriction of
$\mathcal H\to A/(\delta)$ to $A$ is the natural map, whose kernel is
$\delta A$; together with exactness and $A\hookrightarrow\mathcal H$
this gives $A\cap\delta\mathcal H=\delta A$.
\end{proof}

The Taylor maps of Theorem~\ref{thm:rel-scalar-gluing} detect
divisibility by a fixed finite product of cyclotomic polynomials.

\begin{corollary}\label{cor:rel-jet-divisibility}
Let $b\in\Habi$, and let $\delta$ be as in
Lemma~\ref{lem:rel-cyclotomic-quotient}.  Suppose that at every
primitive $m$th root $\zeta$ for which $\Phi_m$ divides $\delta$,
the Taylor series of $b$ has order at least
$v_{\Phi_m}(\delta)$.  Then $b\in\delta\Habi$.
\end{corollary}

\begin{proof}
Let $b_0\in\A$ represent the class of $b$ in
$\Habi/\delta\Habi\simeq\A/(\delta)$. Since $b-b_0$ is a multiple of
$\delta$, the Taylor series of $b_0$ vanish to the same orders as
those of $b$. Over $\mathbb Q$ the Laurent polynomial $b_0$ is
therefore divisible by each $\Phi_m^{v_{\Phi_m}(\delta)}$, hence by
$\delta$. Its class in $\A/(\delta)$ thus becomes zero after tensoring
with $\mathbb Q$, and since this group is torsion-free, the class is
zero. Hence $b\in\delta\Habi$.
\end{proof}

The next lemma is the local input: divided differences of a rational
local component have no poles in $x$, even though the
interpolation nodes collide at $x=0$.

\begin{lemma}\label{lem:rel-local-divided-difference}
Fix a primitive $m$th root $\zeta$ and a multiindex
$\boldsymbol\alpha\in\mathbb N^d$.  Suppose
\[
 g(\mathbf t,x)
  \in\mathbb Q(\zeta)[\mathbf t^{\pm1},D(\mathbf t^m)^{-1}]
          [[x]].
\]
Then the mixed divided difference of its values on the rectangle
$t_i=(\zeta+x)^{j_i}$, $0\le j_i\le\alpha_i$, lies in
$\mathbb Q(\zeta)[[x]]$; that is, it has no negative powers of
$x$.
\end{lemma}

\begin{proof}
Set $B=\Q(\zeta)[[x]]$ and consider the finite free algebra
\[
 C_\alpha=B[t_1^{\pm1},\ldots,t_d^{\pm1}]
 \bigg/\left(\prod_{j=0}^{\alpha_i}(t_i-(\zeta+x)^j)
                                  :1\le i\le d\right).
\]
The nodes are units, so this is also the polynomial quotient, with
basis the monomials $\mathbf t^{\mathbf b}$, $0\le b_i\le\alpha_i$.

Over $\Frac(B)$ the nodes are distinct, so multiplication by
$D(\mathbf t^m)$ is diagonalized by evaluation, and
\[
 \det M_{D(\mathbf t^m)}
 =\prod_{0\le j_i\le\alpha_i}
       D((\zeta+x)^{mj_1},\ldots,(\zeta+x)^{mj_d})
 \in1+x B,
\]
because $(\zeta+x)^{mj}\equiv\zeta^{mj}=1$ modulo $x$. The
determinant is computed from a matrix with entries in $B$, so the
identity holds over $B$, and the adjugate formula shows that
$D(\mathbf t^m)$ is invertible in $C_\alpha$.

The localization therefore maps to $C_\alpha$. Since $C_\alpha$ is
finite free over $B$ and $x$-adically complete, the whole series
$g$ has an image there, with a unique representative
\[
 \overline g=\sum_{0\le b_i\le\alpha_i}a_{\mathbf b}\mathbf t^{\mathbf b},
 \qquad a_{\mathbf b}\in B.
\]
Over $\Frac(B)$ this is the interpolating polynomial of the values of
$g$ at the nodes, so its top coefficient $a_\alpha$ is the mixed
divided difference, which therefore lies in $B$.
\end{proof}

\begin{remark}\label{rem:rel-local-scope}
The exponent of the denominator in the lemma may grow with the order
in $x$, since completeness of $C_\alpha$ allows the series to be
evaluated coefficientwise. For integral $g$ the same proof works over
$\mathbb Z[\zeta][[x]]$.
\end{remark}

For $\boldsymbol\alpha\in\mathbb N^d$ define
\begin{equation}\label{eq:rel-mixed-denominator}
 \delta_{\boldsymbol\alpha}
       =\prod_{i=1}^d (Q;Q)_{\alpha_i}
\end{equation}
and
\begin{equation}\label{eq:rel-mixed-numerator}
 b_{\boldsymbol\alpha}=
 \sum_{0\le j_i\le\alpha_i}
 (-1)^{\sum_i j_i}
 Q^{\sum_i(\binom{j_i+1}{2}-\alpha_i j_i)}
 \prod_{i=1}^d
   \begin{bmatrix}\alpha_i\\j_i\end{bmatrix}_Q
 f_{\mathbf j}.
\end{equation}
The Gaussian coefficients are integral, so
$b_{\boldsymbol\alpha}\in\Habi$. Tensor Lagrange interpolation over a
field shows that the Newton coefficient we are looking for is
\begin{equation}\label{eq:rel-formal-coefficient}
 c_{\boldsymbol\alpha}
       =b_{\boldsymbol\alpha}/\delta_{\boldsymbol\alpha}.
\end{equation}
This rests on the one-variable identity
\[
 \prod_{\substack{0\le k\le a\\k\ne j}}(Q^j-Q^k)
   =(-1)^jQ^{aj-\binom{j+1}{2}}(Q;Q)_j(Q;Q)_{a-j},
\]
whose reciprocal, multiplied by $(Q;Q)_a$, is the weight of
$f_{\mathbf j}$ in~\eqref{eq:rel-mixed-numerator}.

\begin{proposition}\label{prop:rel-integral-coefficients}
For every $\boldsymbol\alpha$ there exists a unique
$c_{\boldsymbol\alpha}\in\Habi$ satisfying
\[
 b_{\boldsymbol\alpha}
       =\delta_{\boldsymbol\alpha}c_{\boldsymbol\alpha}.
\]
If all the color values $f_{\mathbf r}$ lie in $\A$, then all
$c_{\boldsymbol\alpha}$ lie in $\A$.
\end{proposition}

\begin{proof}
At a primitive $m$th root $\zeta$, the Taylor series of the
$f_{\mathbf j}$ are the values of the single local component
$g_\zeta$ on the rectangle $t_i=(\zeta+x)^{j_i}$. Read in
$\mathbb Q(\zeta)((x))$, the quotient
\eqref{eq:rel-formal-coefficient} is therefore the mixed divided
difference of $g_\zeta$, which has no negative powers of $x$ by
Lemma~\ref{lem:rel-local-divided-difference}. Hence
\[
 \operatorname{ord}_{x=0}b_{\boldsymbol\alpha,\zeta}
       \ge v_{\Phi_m}(\delta_{\boldsymbol\alpha}),
\]
where $b_{\boldsymbol\alpha,\zeta}$ is the Taylor series of
$b_{\boldsymbol\alpha}$ at $\zeta$. Corollary~\ref{cor:rel-jet-divisibility}
now gives $b_{\boldsymbol\alpha}\in\delta_{\boldsymbol\alpha}\Habi$,
and the quotient is unique because $\delta_{\boldsymbol\alpha}$ is a
non-zero-divisor by Lemma~\ref{lem:rel-cyclotomic-quotient}.

If the color values lie in $\A$, then so does $b_{\boldsymbol\alpha}$,
and the identity $\A\cap\delta_{\boldsymbol\alpha}\Habi
=\delta_{\boldsymbol\alpha}\A$ of
Lemma~\ref{lem:rel-cyclotomic-quotient} gives
$b_{\boldsymbol\alpha}\in\delta_{\boldsymbol\alpha}\A$. The quotient in
$\A$ is then the unique quotient in $\Habi$, so
$c_{\boldsymbol\alpha}\in\A$.
\end{proof}

\begin{lemma}\label{lem:rel-color-uniqueness}
A family of local components of the form~\eqref{eq:rel-local}, with
no gluing condition imposed between different roots, is determined by
its values at all dominant colors.
\end{lemma}

\begin{proof}
It suffices to show that a family whose colored values all vanish is
zero. Fix $\zeta$. Color Taylor expansion at $\mathbf t=\mathbf1$ is
injective on $\mathbb Q(\zeta)[\mathbf t^{\pm1},D(\mathbf t^m)^{-1}]$,
as one sees by multiplying by the denominator and a monomial, so it
suffices to show that the joint germ of $g_\zeta$ at
$(\mathbf t,Q)=(\mathbf1,\zeta)$ vanishes. This germ vanishes under
the substitutions $t_i=Q^{ma_i}$, $a_i\in\mathbb N$, because these
give the Taylor expansions of the colored values at the colors
$m\mathbf a$. Lemma~\ref{lem:rel-color-rigidity}, applied with
$b_i=0$ and $Q-\zeta=\zeta\epsilon$, shows that the germ is zero.
\end{proof}

\begin{theorem}\label{thm:rel-main-embedding}
Color evaluation followed by tensor Newton interpolation defines
a canonical injective ring homomorphism
\[
 \iota_D:\mathcal H_{\rm rel}(D)\hookrightarrow\mathcal N_d(\Habi).
\]
Let $\mathcal H_{\rm rel}(D)_{\rm pol}$ denote the subring of
elements for which every dominant color value $f_{\mathbf r}$
belongs to $\A$.  Then
\[
 \iota_D\bigl(\mathcal H_{\rm rel}(D)_{\rm pol}\bigr)
                  \subseteq\mathcal N_d(\A).
\]
\end{theorem}

\begin{proof}
Using the coefficients of
Proposition~\ref{prop:rel-integral-coefficients}, set
\[
 \iota_D(g)=\sum_\alpha c_\alpha P_\alpha\in\mathcal N_d(\mathcal H).
\]
We first check that $\iota_D(g)$ has the values $f_{\mathbf r}$. Every
node difference $Q^a-Q^b$, $a\ne b$, is a Laurent unit times a
product of cyclotomic polynomials, hence a non-zero-divisor in
$\mathcal H$ by Lemma~\ref{lem:rel-cyclotomic-quotient}. Let $\Pi$ be
the product of the node differences on a given box. In
$\mathcal H[\Pi^{-1}]$ the $c_\alpha$ are the Newton coefficients of
the data $(f_{\mathbf j})$, so tensor interpolation there gives
\[
 f_{\mathbf r}
 =\sum_{\alpha_i\le r_i}c_\alpha
       P_\alpha(Q^{r_1},\ldots,Q^{r_d}).
\]
Both sides lie in $\mathcal H$, and $\mathcal H\to\mathcal H[\Pi^{-1}]$
is injective, so the identity holds in $\mathcal H$.

The evaluation matrix on a box is triangular with diagonal entries
that are non-zero-divisors in $\mathcal H$, so evaluation on the
dominant grid is injective on $\mathcal N_d(\mathcal H)$. Since color
specialization is a ring homomorphism on $\mathcal H_{\rm rel}(D)$,
the elements $\iota_D(g+g')$ and $\iota_D(g)+\iota_D(g')$ have the same
values, and likewise for products; hence $\iota_D$ is a ring
homomorphism. The assertion about $\mathcal H_{\rm rel}(D)_{\rm pol}$
follows from Proposition~\ref{prop:rel-integral-coefficients}.
Finally,
\[
 \iota_D(g)=0\ \Longrightarrow\ f_{\mathbf r}=0
       \text{ for every }\mathbf r
 \ \Longrightarrow\ g=0
\]
by Lemma~\ref{lem:rel-color-uniqueness}, which proves injectivity.
\end{proof}

\subsection{The image criterion}

Which integral Newton series come from relative collections? Every
Newton series already has integral local germs at the identity color,
and these germs are compatible under re-expansion
(Proposition~\ref{prop:rel-automatic-gluing}). The only question is
whether they are rational in the colors.

\begin{proposition}\label{prop:rel-automatic-gluing}
Every $F\in\mathcal N_d^+(A)$ has integral, Galois-equivariant joint
germs at all $(\mathbf t,Q)=(\zeta^{\mathbf b},\zeta)$, for $\zeta$ of
order $m$ and $0\le b_i<m$. At the identity color center, these germs
are compatible under re-expansion at every pair of $p$-adically close
roots of orders $m$ and $pm$. Their colored specializations are the
ordinary Taylor expansions of the Laurent polynomial grid values of
$F$.
\end{proposition}

\begin{proof}
Write $F=\sum_\alpha c_\alpha(Q)P_\alpha(\mathbf t;Q)$.
At $(\zeta^{\mathbf b},\zeta)$,
\[
 \ord P_\alpha
 =\sum_i\#\{0\le j<\alpha_i:j\equiv b_i\pmod m\}
 \longrightarrow\infty\qquad(|\alpha|\longrightarrow\infty),
\]
so the Newton sum converges to a joint germ with coefficients in
$\Z[\zeta]$. Galois equivariance holds termwise. On $t_i=Q^{r_i}$ the
sum truncates to $\alpha_i\le r_i$, which gives the ordinary Taylor
expansion of the Laurent polynomial grid value.

For the re-expansion, take $\zeta,\xi$ of orders $m,pm$ in a common
complete cyclotomic valuation ring $\mathcal O$ with uniformizer
$\pi$, such that $\xi\equiv\zeta\pmod\pi$. Put
\[
 x=Q-\xi,\qquad w_i=t_i-1,\qquad
 \mathfrak m=(\pi,x,w_1,\ldots,w_d).
\]
If $m\mid j$, then $\xi^j$ is a $p$th root of unity, hence congruent
to $1$ modulo $\pi$, and
\[
 1+w_i-(\xi+x)^j\in\mathfrak m.
\]
Consequently
\[
 c_\alpha(\xi+x)P_\alpha(\mathbf1+\mathbf w;\xi+x)
 \in\mathfrak m^{\sum_i\lfloor\alpha_i/m\rfloor},
\]
where $c_\alpha(\xi+x)$ has integral coefficients because $\xi$ is a
unit. The Newton sum therefore converges in the $\mathfrak m$-adically
complete ring $\mathcal O[[x,\mathbf w]]$. The substitution
$Q-\zeta=x+\xi-\zeta$ is continuous for this topology, and it carries
the germ at $\zeta$ to this same sum, which is also the germ at $\xi$.
This proves the compatibility.
\end{proof}

\begin{theorem}
\label{thm:rel-exact-image}
Let $D\in\Z[\mathbf T^{\pm1}]$ satisfy $D(\mathbf1)=1$, and let
$F\in\mathcal N_d^+(A)$. The following conditions are equivalent.
\begin{enumerate}
\item $F$ lies in the image of
$\mathcal H_{\rm rel}(D)_{\rm pol}\hookrightarrow\mathcal N_d^+(A)$.
\item At every primitive $m$th root $\zeta$, the joint germ of $F$ at
$\mathbf t=\mathbf1$, written as a series in $x=Q-\zeta$, is the
color Taylor expansion of a series
\begin{equation}\label{eq:rel-image-condition}
 g_\zeta(\mathbf t,x)
  \in\Q(\zeta)[\mathbf t^{\pm1},D(\mathbf t^m)^{-1}][[x]].
\end{equation}
\end{enumerate}
When these conditions hold, the representatives have integral
cyclotomic numerators, they form the unique relative lift, and they
agree with the Newton germs at every torsion-color center.
\end{theorem}

In \eqref{eq:rel-image-condition} the power of the denominator may
grow with the order in $x$.

\begin{proof}
Suppose first that $F$ has a relative preimage. At $\zeta$, its local
component and the Newton germ of $F$ have the same values on the
curves
\[
 t_i=(\zeta+x)^{ma_i},\qquad a_i\ge0,
\]
namely the Taylor expansions of the colored values at the colors
$m\mathbf a$. After the rescaling $x=\zeta\epsilon$,
Lemma~\ref{lem:rel-color-rigidity} identifies the two germs, which
proves condition~2.

Conversely, assume condition~2. Rational Taylor expansion at
$\mathbf1$ is injective, so each $g_\zeta$ is unique. Since the Newton
germ has integral coefficients, Lemma~\ref{lem:rel-integral-saturation},
applied to each coefficient of $x^h$ with $D_0=D(\mathbf t^m)$, gives
\[
 g_\zeta\in
 \Z[\zeta][\mathbf t^{\pm1},D(\mathbf t^m)^{-1}][[x]].
\]
Uniqueness and the Galois equivariance of the Newton germs give Galois
equivariance of the $g_\zeta$.

For the gluing, let $\zeta,\xi$ be $p$-adically close roots of orders
$m$ and $pm$, and work in
\[
 \widehat{\mathcal O[\mathbf t^{\pm1},D(\mathbf t^m)^{-1}]}_p[[x]].
\]
The denominator $D(\mathbf t^{pm})$ is invertible in this ring because
$D(\mathbf t^{pm})\equiv D(\mathbf t^m)^p\pmod p$. The coefficients of
$g_\zeta$ are integral and $\xi-\zeta$ has positive valuation, so the
re-expansion $g_\zeta(\mathbf t,x+\xi-\zeta)$ converges in this ring.
Its color germ is the re-expanded Newton germ at $\zeta$, which by
Proposition~\ref{prop:rel-automatic-gluing} is the Newton germ at
$\xi$, that is, the color germ of $g_\xi$. By
Lemma~\ref{lem:rel-completed-injection},
\[
 g_\zeta(\mathbf t,x+\xi-\zeta)=g_\xi(\mathbf t,x).
\]
Writing $t_i=T_i^{1/m}$ turns this into
\eqref{eq:rel-glue}, so the representatives form a
relative collection.

At a dominant color $\mathbf r$, let $h_{\mathbf r}\in\mathcal H$ be
the specialization of this collection. Its Taylor series at $Q=1$
agrees with that of the Laurent polynomial grid value of $F$, so the
injectivity of the Taylor map at $Q=1$
(Theorem~\ref{thm:rel-scalar-gluing}) gives
\[
 h_{\mathbf r}=F(Q^{r_1},\ldots,Q^{r_d};Q)\in A.
\]
The collection therefore belongs to the polynomial-valued subring, and
its Newton image is $F$ by injectivity of evaluation. The embedding
makes the preimage unique. Finally, applying
Lemma~\ref{lem:rel-color-rigidity} to the curves
$t_i=Q^{ma_i+b_i}$ identifies the germs at every torsion-color center
$(\zeta^{\mathbf b},\zeta)$.
\end{proof}

\subsection{Cyclotomic neighborhoods}\label{subsec:rel-blocks}

The criterion of Theorem~\ref{thm:rel-exact-image} is stated in terms
of germs. In this subsection and the next it is translated into a
condition on the Newton coefficients. The first step is to describe a
Newton series near $Q=\zeta$ simultaneously at all torsion colors.

\begin{theorem}
\label{thm:rel-integral-cyclotomic-neighborhood}
For $m,M\ge1$, put $B_{m,M}=A/(\Phi_m(Q)^M)$ and $z_i=t_i^m-1$.
There is a natural isomorphism
\begin{equation}\label{eq:rel-integral-cyclotomic-neighborhood}
 \frac{\mathcal N_d^+(A)}{\Phi_m(Q)^M\mathcal N_d^+(A)}
 \cong
 B_{m,M}[[z_1,\ldots,z_d]][t_1,\ldots,t_d]
             /(t_i^m-1-z_i)_{i=1}^d.
\end{equation}
The right side is finite free of rank $m^d$ over
$B_{m,M}[[\mathbf z]]$.
\end{theorem}

\begin{proof}
Reducing the Newton coefficients identifies the left side with the
Newton completion over $B_{m,M}$. In one variable we use the cofinal
subsequence of Newton polynomials $P_{mk}$. Its $j$th complete block
satisfies
\[
 \prod_{r=0}^{m-1}(t-Q^{mj+r})\equiv t^m-1=z\pmod{\Phi_m(Q)}.
\]
At a primitive $m$th root the nodes $Q^{mj+r}$, $0\le r<m$, run
through all $m$th roots of unity, so the difference of the two sides
vanishes at every primitive $m$th root; since $\Phi_m$ is monic,
division by it shows that the difference is divisible by $\Phi_m$ in
$A[t]$. Every block thus has the form $z+\Phi_m V_j(t)$, and in the
product of $k$ blocks every term containing $M$ factors
$\Phi_m V_j$ vanishes over $B_{m,M}$. Hence, for $k\ge M$,
\[
 P_{mk}(t)\in (z^{k-M+1}).
\]
Conversely, $P_{mk}=z^k+\Phi_m g$ for some $g\in B_{m,M}[t]$, so
$z^k\in\Phi_m B_{m,M}[t]$ modulo $P_{mk}$, and $z^{kM}=0$ in
$B_{m,M}[t]/(P_{mk})$. The ideals $(P_{mk})$ and the powers of $(z)$
are therefore cofinal, and the same holds in each coordinate.

The polynomial ring $B_{m,M}[\mathbf t]$ is finite free over
$B_{m,M}[\mathbf z]$ with basis $\mathbf t^{\mathbf b}$,
$0\le b_i<m$, and its $\mathbf z$-adic completion is the right side
of~\eqref{eq:rel-integral-cyclotomic-neighborhood}. In this completion
each $t_i$ is a unit, because $t_i^m=1+z_i$ is, so it is also the
completion of the Laurent polynomial ring. This proves the isomorphism
and the rank statement.
\end{proof}

Fix a primitive $m$th root $\zeta$, put $\Bbbk=\mathbb Q(\zeta)$ and
$x=Q-\zeta$, and introduce
\[
 z_i=t_i^m-1,
 \qquad
 H_b(t)=P_b(t;\zeta)\quad(0\le b<m),
 \qquad H_{\mathbf b}(\mathbf t)=\prod_iH_{b_i}(t_i).
\]
Consider the finite free algebra
\begin{equation}\label{eq:rel-block-simultaneous-ring}
 \mathscr C_\zeta
 =\Bbbk[[x,z_1,\ldots,z_d]]
       [t_1,\ldots,t_d]/(t_i^m-1-z_i)_{i=1}^d.
\end{equation}
The $t_i$ are units in it, and the $H_{\mathbf b}$ with $0\le b_i<m$
form a basis, since the change of basis from the
$\mathbf t^{\mathbf b}$ is unitriangular.

Every tensor Newton series has a well-defined image in
$\mathscr C_\zeta$. If $\alpha_i=mk_i+b_i$, then
\[
 P_{\alpha_i}(t_i;\zeta+x)
 \equiv z_i^{k_i}H_{b_i}(t_i)\pmod x,
\]
and since $P_{\alpha_i}$ is a product of $k_i$ complete blocks of
length $m$, each congruent to $z_i$ modulo $x$, it lies in
$(x,z_i)^{k_i}$. The series therefore converges in the
$(x,\mathbf z)$-adic topology.

Over $\Bbbk$ the ring~\eqref{eq:rel-block-simultaneous-ring} is the
product of the $m^d$ formal color charts
\[
 t_i=\tau_i(1+z_i)^{1/m},\qquad
 \tau_i\in\mu_m=\{\tau\in\Bbbk:\tau^m=1\}.
\]
This follows by factoring $t_i^m-(1+z_i)$ over $\Bbbk[[z_i]]$: the
differences of its roots have nonzero constant term, so the Chinese
remainder theorem applies. Equality in $\mathscr C_\zeta$ therefore
means equality of the germs at all torsion colors, whereas the
projection to the single chart $\boldsymbol\tau=\mathbf1$ is not
injective.

For
\[
 F=\sum_{\boldsymbol\alpha}c_{\boldsymbol\alpha}(Q)
                              P_{\boldsymbol\alpha}(\mathbf t;Q),
\]
write the Taylor expansion of the coefficients as
\[
 c_{\boldsymbol\alpha}(\zeta+x)
      =\sum_{h\ge0}c^{\zeta,h}_{\boldsymbol\alpha}x^h,
 \qquad c^{\zeta,h}_{\boldsymbol\alpha}\in\mathbb Z[\zeta],
\]
and define the raw block series
\begin{equation}\label{eq:rel-raw-block-series}
 C_{\zeta,\mathbf b,h}(\mathbf z)
    =\sum_{\mathbf k\in\mathbb N^d}
       c^{\zeta,h}_{m\mathbf k+\mathbf b}\mathbf z^{\mathbf k}.
\end{equation}
The image of $F$ in $\mathscr C_\zeta$ has unique coordinates
$B_{\zeta,\mathbf b,h}(\mathbf z)\in\Bbbk[[\mathbf z]]$, the transported
blocks, defined by
\begin{equation}\label{eq:rel-transported-blocks}
 F=\sum_{h\ge0}x^h
          \sum_{0\le b_i<m}H_{\mathbf b}(\mathbf t)
                           B_{\zeta,\mathbf b,h}(\mathbf z).
\end{equation}
In quantum order zero the two kinds of blocks coincide. In higher
orders they differ, because the Newton polynomials themselves depend
on $x$. The next theorem describes the difference. We write
$\mathbf B_{\zeta,h}$ and $\mathbf C_{\zeta,h}$ for the vectors of
blocks indexed by $\mathbf b$.

\begin{theorem}
\label{thm:rel-block-transport}
For every $\ell\ge0$ there is a matrix $\mathcal L_{\zeta,\ell}$ of
polynomial-coefficient differential operators in $\mathbf z$, with
coefficients in $\Bbbk$, such that
\begin{equation}\label{eq:rel-block-triangular-transport}
 \mathbf B_{\zeta,h}
 =\sum_{j=0}^h\mathcal L_{\zeta,h-j}\mathbf C_{\zeta,j}.
\end{equation}
Here $\mathcal L_{\zeta,0}$ is the identity, and every entry of
$\mathcal L_{\zeta,\ell}$ has differential order at most $2\ell$.
These matrices depend only on $m,\zeta,d$ and the residue basis, not
on $F$.
\end{theorem}

\begin{proof}
In one variable fix $b<m$ and factor
\[
 P_{mk+b}(t;\zeta+x)
 =\prod_{j=0}^{k-1}G_j(t,x)\,H_{k,b}(t,x),
\]
where
\[
 G_j=\prod_{r=0}^{m-1}\bigl(t-(\zeta+x)^{mj+r}\bigr),
 \qquad
 H_{k,b}=\prod_{r=0}^{b-1}\bigl(t-(\zeta+x)^{mk+r}\bigr).
\]
At $x=0$ these are $z=t^m-1$ and $H_b(t)$, respectively.
Since
\[
 [x^a](\zeta+x)^{mj+r}
       =\binom{mj+r}{a}\zeta^{r-a},
\]
the coefficient $[x^a]G_j$ is a polynomial in $j$ of degree at
most $a$, and $[x^a]H_{k,b}$ is a polynomial in $k$ of degree at
most $a$.

The coefficient of $x^\ell$ in the product is a sum over the
choices of $s\le\ell$ complete blocks $j_1<\cdots<j_s<k$ of positive
order, with the remaining order taken by the final partial block.
Summing a polynomial in $j_1,\ldots,j_s$ over
$0\le j_1<\cdots<j_s<k$ raises its degree by at most $s$, as one sees
by applying
\[
 \sum_{j=0}^{k-1}(j)_{\underline r}
       =\frac{(k)_{\underline{r+1}}}{r+1}
\]
once for each index. This gives
\begin{equation}\label{eq:rel-one-block-derivative}
 [x^\ell]P_{mk+b}(t;\zeta+x)
 =\sum_{s=0}^{\ell}z^{k-s}A_{\ell,s,b}(t,k),\qquad
 \deg_k A_{\ell,s,b}\le\ell+s,
\end{equation}
where terms with $s>k$ are zero.

For $s>0$ the sum defining $A_{\ell,s,b}$ is empty when
$k=0,\ldots,s-1$, so $A_{\ell,s,b}$ vanishes there, and its expansion
in falling factorials starts in degree $s$:
\[
 A_{\ell,s,b}(t,k)
 =\sum_{r=s}^{\ell+s}a_{\ell,s,b,r}(t)(k)_{\underline r}.
\]
For $C(z)=\sum_{k\ge0}u_kz^k$ this yields
\[
 \begin{split}
 \sum_{k\ge0}u_kz^{k-s}A_{\ell,s,b}(t,k)
 &=\sum_{r=s}^{\ell+s}a_{\ell,s,b,r}(t)
                  z^{r-s}\partial_z^r C(z),\\
 r-s&\ge0,\qquad r\le\ell+s\le2\ell.
 \end{split}
\]
The same formula holds for $s=0$. Reducing each coefficient in $t$ to
the basis $H_0,\ldots,H_{m-1}$ by means of $t^m=1+z$ introduces only
polynomial coefficients in $z$. This constructs the one-variable
matrix of differential order at most $2\ell$, and in order zero it is
the identity.

In $d$ variables the order $\ell$ is distributed among the tensor
factors. Differential orders add, so the resulting matrix still has
order at most $2\ell$. Finally, distributing the total order $h$
between the coefficient $c_\alpha$ and the Newton polynomial gives
\eqref{eq:rel-block-triangular-transport}. Differential operators of
finite order are continuous in the formal topology, so these
manipulations are compatible with the convergent expansion in
$\mathscr C_\zeta$.
\end{proof}

\begin{example}
\label{ex:rel-first-block-transport}
In one variable write $G(t,Q)=P_m(t;Q)$ and $H_b(t,Q)=P_b(t;Q)$.
In the next formula a subscript $Q$ or $t$ denotes the corresponding
partial derivative, evaluated at $Q=\zeta$. For block series
$C_b(z)$, the contribution in quantum order one of the Newton
coefficients of order zero is
\begin{align}\label{eq:rel-first-block-operator}
 \sum_{b=0}^{m-1}\biggl\{
 &H_b(t)G_Q(t)\,C_b'(z)
 -\frac{m^2}{2\zeta}H_b(t)z\,C_b''(z)
 +H_{b,Q}(t)C_b(z)\\
 &\hspace{14mm}
 -\frac m\zeta\bigl(tH_{b,t}(t)-bH_b(t)\bigr)z\,C_b'(z)
 \biggr\}.\nonumber
\end{align}
Reducing this expression in the residue basis gives
$\mathcal L_{\zeta,1}\mathbf C$. The formula follows from
\[
 [x]G_j=G_Q-\frac{m^2j}{\zeta},\qquad
 [x]H_{k,b}=H_{b,Q}
       -\frac{mk}{\zeta}(tH_{b,t}-bH_b),
\]
together with $\sum_{j<k}j=k(k-1)/2$. For $m=1$ it reduces to the
scalar operator $\mathcal L_{1,1}=-\tfrac12 z\partial_z^2$.
\end{example}

\subsection{The coefficient criterion}

Fix $D\in\mathbb Z[\mathbf T^{\pm1}]$ with $D(\mathbf1)=1$, and put
\[
 S_D(\Bbbk)=\Bbbk[(1+z_1)^{\pm1},\ldots,(1+z_d)^{\pm1},
                         D(1+\mathbf z)^{-1}].
\]
This is a differential subring of $\Bbbk[[\mathbf z]]$: differentiation
only produces further powers of the inverted elements $1+z_i$ and
$D(1+\mathbf z)$.

\begin{proposition}
\label{prop:rel-fixed-root-block-recognition}
The image of $F$ in $\mathscr C_\zeta$ is the expansion of an element
of
\[
 \Bbbk[\mathbf t^{\pm1},D(\mathbf t^m)^{-1}][[x]]
\]
if and only if $C_{\zeta,\mathbf b,h}\in S_D(\Bbbk)$ for every
$\mathbf b$ and $h$. Here the expansions are compared simultaneously
at all $\boldsymbol\tau\in\mu_m^d$.
\end{proposition}

\begin{proof}
With $T_i=t_i^m$, the localized Laurent ring
$\Bbbk[\mathbf t^{\pm1},D(\mathbf t^m)^{-1}]$ is free over
$\Bbbk[\mathbf T^{\pm1},D(\mathbf T)^{-1}]$ with basis the
$H_{\mathbf b}$. Its image in $\mathscr C_\zeta$ therefore consists
precisely of the series whose transported blocks satisfy
\[
 \mathbf B_{\zeta,h}\in S_D(\Bbbk)^{m^d}
       \qquad\text{for every }h.
\]
Each transport matrix $\mathcal L_{\zeta,\ell}$ preserves
$S_D(\Bbbk)^{m^d}$, since $S_D(\Bbbk)$ is closed under multiplication
by polynomials and under differentiation. The identity diagonal in
\eqref{eq:rel-block-triangular-transport} gives the inverse recursion
\[
 \mathbf C_{\zeta,h}
 =\mathbf B_{\zeta,h}
   -\sum_{\ell=1}^{h}\mathcal L_{\zeta,\ell}
                          \mathbf C_{\zeta,h-\ell},
\]
and induction on $h$ shows that the conditions on the transported and
on the raw blocks are equivalent.
\end{proof}

\begin{theorem}
\label{thm:rel-cyclotomic-block-criterion}
For $F\in\mathcal N_d^+(A)$, the following conditions are equivalent.
\begin{enumerate}
\item $F$ belongs to the image of
$\mathcal H_{\mathrm{rel}}(D)_{\mathrm{pol}}
\hookrightarrow\mathcal N_d^+(A)$.
\item For every primitive $m$th root $\zeta$, every
$\mathbf b\in\{0,\ldots,m-1\}^d$, and every $h\ge0$, the raw
coefficient series~\eqref{eq:rel-raw-block-series} belongs to
$S_D(\Bbbk)$.
\end{enumerate}
When they hold, the block series in condition 2 lie in
\[
 \mathbb Z[\zeta][(1+z_1)^{\pm1},\ldots,(1+z_d)^{\pm1},
                                  D(1+\mathbf z)^{-1}].
\]
\end{theorem}

\begin{proof}
If $F$ has a relative preimage, its local component and the Newton
series of $F$ have the same values on every curve
$t_i=Q^{ma_i+b_i}$. By Lemma~\ref{lem:rel-color-rigidity} they have
the same germ in every torsion-color chart, hence the same image in
$\mathscr C_\zeta$, and
Proposition~\ref{prop:rel-fixed-root-block-recognition} gives
condition~2.

Conversely, condition~2 and the same proposition give a rational
representative in all charts, in particular in the chart of the
identity color. Theorem~\ref{thm:rel-exact-image} then
supplies the unique relative preimage.

For integrality, write a raw block as $C=N/E$, where $N$ and $E$ are
polynomials over $\Bbbk$, $E$ is a product of powers of the $1+z_i$
and of the polynomial obtained from $D(1+\mathbf z)$ by clearing the
negative powers of the units $1+z_i$, and $E(\mathbf0)=1$. Then
$E\in\Z[\mathbf z]$, and by definition
\[
 C\in\Z[\zeta][[\mathbf z]],\qquad
 N=EC\in\Bbbk[\mathbf z]\cap\Z[\zeta][[\mathbf z]]
       =\Z[\zeta][\mathbf z],
\]
which proves the integral assertion.
\end{proof}

Clearing the permitted denominators turns the criterion into linear
recurrences for the blocks of Newton coefficients.

\begin{corollary}
\label{cor:rel-block-recurrence}
Choose a Laurent monomial $\mathbf T^{\boldsymbol\ell}$ so that
$\widetilde D(\mathbf T)=\mathbf T^{\boldsymbol\ell}D(\mathbf T)$
is a polynomial, and write $d_D(\mathbf z)=\widetilde D(1+\mathbf z)$.
Condition 2 of Theorem~\ref{thm:rel-cyclotomic-block-criterion} holds
if and only if for every $\zeta,\mathbf b,h$ there are $e\ge0$,
$a_i\ge0$, and a polynomial
$N(\mathbf z)\in\mathbb Z[\zeta][\mathbf z]$ such that
\begin{equation}\label{eq:rel-block-recurrence-symbol}
 \prod_i(1+z_i)^{a_i}d_D(\mathbf z)^e
       C_{\zeta,\mathbf b,h}(\mathbf z)=N(\mathbf z).
\end{equation}
Equivalently, if
$E(\mathbf z)=\sum_{\boldsymbol\nu}e_{\boldsymbol\nu}
\mathbf z^{\boldsymbol\nu}$ denotes the polynomial multiplying $C$ on
the left and $u_{\mathbf k}=c^{\zeta,h}_{m\mathbf k+\mathbf b}$, with
$u_{\mathbf k}=0$ when a coordinate of $\mathbf k$ is negative, then
\[
 \sum_{\boldsymbol\nu}e_{\boldsymbol\nu}
                  u_{\mathbf k-\boldsymbol\nu}=0
\]
for all but finitely many $\mathbf k\in\mathbb N^d$.
The admissible recurrence polynomials do not depend on the order $m$
of the root; $m$ only determines the subsequences of coefficients to
which they are applied.
\end{corollary}

\begin{proof}
Every element of $S_D(\Bbbk)$ can be written with a denominator that
is a product of powers of the $1+z_i$ and of $d_D$, and by the last
part of the theorem the numerator can be taken integral. Conversely,
\eqref{eq:rel-block-recurrence-symbol} places the series in
$S_D(\Bbbk)$, because $E(\mathbf0)=1$. Comparing coefficients gives
the formulation in terms of recurrences.
\end{proof}

In type $A_{n-1}$ the relevant denominator is
\[
 D_{K,n}(\mathbf T)=\prod_{1\le a<b\le n}
                    \Delta_K(T_aT_{a+1}\cdots T_{b-1}).
\]
The type-$A$ lift constructed below places every coefficient block, at
every root and in every quantum order, in this Alexander localization
in the variables $T_i=1+z_i$ (Corollary~\ref{cor:trace-all-rank-block-recurrences}).

\subsection{Examples}

Two families of Newton series with constant or polynomial coefficients
illustrate the criterion.

\begin{proposition}
\label{prop:rel-polynomial-quantum-phases}
Let $p\in\mathbb Z[X_1,\ldots,X_d]$ and put
\[
 F_p(\mathbf t;Q)=\sum_{\boldsymbol\alpha\in\mathbb N^d}
      Q^{p(\boldsymbol\alpha)}P_{\boldsymbol\alpha}(\mathbf t;Q),
 \qquad D_0(\mathbf T)=\prod_{i=1}^d(2-T_i).
\]
Then $F_p$ has a unique relative Habiro lift for $D_0$. If
$e_i=\deg_{X_i}p$, with $e_i=0$ when $p$ does not involve $X_i$, then
the coefficient of $x^h$ in its local component at a primitive
$m$th root lies in
\begin{equation}\label{eq:rel-phase-pole-bound}
 \frac{\mathbb Z[\zeta][\mathbf t]}
 {\displaystyle\prod_i(2-t_i^m)^{\,1+h\max(e_i,2)}}.
\end{equation}
This bound is uniform in the order of the root.
\end{proposition}

\begin{proof}
On the residue block $\boldsymbol\alpha=m\mathbf k+\mathbf b$ we have
$p(m\mathbf k+\mathbf b)\equiv p(\mathbf b)\pmod m$, so
\[
 [x^h](\zeta+x)^{p(m\mathbf k+\mathbf b)}
 =\zeta^{p(\mathbf b)-h}
       \binom{p(m\mathbf k+\mathbf b)}h,
\]
where for negative exponents the binomial coefficient is the usual
polynomial. As a function of $\mathbf k$ this is a polynomial of
degree at most $he_i$ in $k_i$. Expanding in falling factorials and
using
\[
 \sum_{k\ge0}(k)_{\underline r}z^k
       =\frac{r!z^r}{(1-z)^{r+1}},
\]
we see that its generating series has denominator dividing
$\prod_i(1-z_i)^{he_i+1}$. Since $D_0(1+\mathbf z)=\prod_i(1-z_i)$,
the raw blocks satisfy the condition of
Theorem~\ref{thm:rel-cyclotomic-block-criterion}.

For the pole bound, raw order $j$ contributes to quantum order $h$
through a differential operator of order at most $2(h-j)$ by
Theorem~\ref{thm:rel-block-transport}, so the exponent of the
denominator in $z_i$ is at most
\[
 je_i+1+2(h-j)\le1+h\max(e_i,2).
\]
The transport coefficients are polynomials, and neither the reduction
to the residue basis nor the substitution $z_i=t_i^m-1$ introduces
denominators, so the numerator is a polynomial over
$\mathbb Q(\zeta)$. The Newton germ is integral at $\mathbf t=\mathbf1$
and the denominators take the value one there, so
Lemma~\ref{lem:rel-integral-saturation} shows that the numerator
is integral. This proves~\eqref{eq:rel-phase-pole-bound}.
\end{proof}

For the one-variable quadratic phase $p(k)=k^2$, the component at the
identity root begins with
\[
 F_p(t;1+x)
 =\frac1{2-t}+x\frac{(t-1)^2}{(2-t)^3}
       +O(x^2).
\]
Here the first two raw blocks are $C_0=(1-z)^{-1}$ and
$C_1=z(1+z)/(1-z)^3$, and the correction
$-\tfrac12zC_0''=-z/(1-z)^3$ comes from
$\mathcal L_{1,1}=-\tfrac12z\partial_z^2$. At a general primitive
root the representative in order zero is
\[
 \frac{\sum_{b=0}^{m-1}\zeta^{b^2}P_b(t;\zeta)}{2-t^m}.
\]

\begin{proposition}
\label{prop:rel-geometric-classification}
For an integer $a$, set
\[
 F_a(t;Q)=\sum_{k\ge0}a^kP_k(t;Q),
\]
so that $F_0=1$. There is a Laurent polynomial $D(T)$ with $D(1)=1$
for which $F_a$ has a relative Habiro lift if and only if
$a\in\{-1,0,1\}$.  For $a=\pm1$ one can take $D(T)=2-T$.
\end{proposition}

\begin{proof}
For $a\ne0$ the Newton coefficients do not depend on $Q$, so the raw
blocks are
\[
 C_{\zeta,b,0}(z)=\frac{a^b}{1-a^mz},\qquad
 C_{\zeta,b,h}=0\quad(h>0).
\]
For $a=1$ their denominator is $1-z=2-T$. For $a=-1$ it is $1-z=2-T$
when the order of the root is even, and $1+z=T$, a unit, when it is
odd. In both cases the block criterion gives a lift, and the case
$a=0$ is trivial.

If $|a|>1$, the block with $b=0$ at a root of order $m$ has a pole at
$T=1+a^{-m}$. These points are nonzero and pairwise distinct for
different $m$. Membership in $\mathbb Q(\zeta)[T^{\pm1},D(T)^{-1}]$ for
all $m$ would make all of them zeros of the single nonzero Laurent
polynomial $D$, which is impossible. By the block criterion there is no
lift.
\end{proof}

\section{Braid operators at roots of unity}
\label{sec:braid}

In this section and the next we prove that the local expansions of the
knot element at a nontrivial root of unity are rational. The strategy
is as follows. At a primitive $m$th root of unity we split each PBW
exponent into its quotient and remainder modulo $m$. The remainders
span a finite-dimensional module, while the quotients become exponents
of polynomial variables. We first express each fixed root transfer as
a polynomial differential operator in these variables, and then sum
the transfers to obtain the local expansion of a braid crossing
(Theorem~\ref{thm:braid-all-rank-insertion}). The traces of the resulting
operators are computed in Section~\ref{sec:trace}.

Write the positive roots of $\mathfrak{sl}_n$ as
$\alpha_{ab}=\varepsilon_a-\varepsilon_b$, $1\le a<b\le n$, in the
lexicographic order on $(a,b)$. Here $\varepsilon_1,\ldots,\varepsilon_n$
is the standard orthonormal basis, so
$(\alpha_i,\alpha_j)=2\delta_{ij}-\delta_{i,j+1}-\delta_{i+1,j}$
and $\operatorname{ht}(\alpha_{ab})=b-a$.
We use the Drinfeld--Jimbo algebra $U_q(\mathfrak{sl}_n)$ over
$\mathbb Q(q)$, with
\[
 K_iE_jK_i^{-1}=q^{(\alpha_i,\alpha_j)}E_j,\quad
 K_iF_jK_i^{-1}=q^{-(\alpha_i,\alpha_j)}F_j,\quad
 [E_i,F_j]=\delta_{ij}\frac{K_i-K_i^{-1}}{q-q^{-1}},
\]
and the usual quantum Serre relations. Set $d=q-q^{-1}$ and fix
root vectors recursively by
\[
 F_{a,a+1}=F_a,\qquad
 F_{ab}=F_{b-1}F_{a,b-1}-qF_{a,b-1}F_{b-1},
\]
\[
 E_{a,a+1}=E_a,\qquad
 E_{ab}=E_{a,b-1}E_{b-1}-q^{-1}E_{b-1}E_{a,b-1}.
\]
The coproduct is
$\Delta(E_i)=E_i\otimes K_i+1\otimes E_i$ and
$\Delta(F_i)=F_i\otimes1+K_i^{-1}\otimes F_i$.
Write $[k]_q=(q^k-q^{-k})/d$,
$[k]_q!=\prod_{j=1}^k[j]_q$, and
$E_\alpha^{(k)}=E_\alpha^k/[k]_q!$.
For independent invertible parameters $u_1,\ldots,u_{n-1}$, let
$M(\mathbf u)$ be the Verma module over
$\mathbb Q(q,u_1,\ldots,u_{n-1})$ generated by $v_0$ with
$K_iv_0=u_iv_0$ and $E_iv_0=0$.
For vectors $w_\gamma,w_\delta$ with root deficits $\gamma,\delta$,
so that $K_iw_\gamma=u_iq^{-(\alpha_i,\gamma)}w_\gamma$, write
$\mathbf u^\gamma=\prod_i u_i^{\gamma_i}$ if
$\gamma=\sum_i\gamma_i\alpha_i$. The highest-normalized Cartan
factor and the checked crossing are
\begin{align}
 C(w_\gamma\otimes w_\delta)
   &=
    \mathbf u^{-\gamma-\delta}q^{(\gamma,\delta)}
                       w_\gamma\otimes w_\delta,\notag\\
 \check R
   &=P C\prod_{\alpha>0}^{\longrightarrow}
       \left(\sum_{k\ge0}q^{k(k-1)/2}d^k
                          E_\alpha^{(k)}\otimes F_\alpha^k\right).
 \label{eq:braid-normalized-crossing}
\end{align}
Here $P$ is the tensor flip, the root factors are taken in the
lexicographic convex order, and products act from right to left.
This is the convex-order factorization of the universal $R$-matrix
with the Cartan scalar of the highest vector removed; see
\cite[Theorem~4.1.2 and Corollary~40.2.2]{Lusztig} for the
quasi-$R$-matrix and the PBW construction, and \cite{RT} for the tangle
invariant. The removed scalar is restored in the zero-framed
normalization of the cut. On any fixed PBW vector in the first factor
only finitely many terms of the root sums act nontrivially. For a
signed braid $\beta$ we denote the resulting operator by
$\mathcal B_\beta$. The interval straightening relations used below
agree with \cite[Section~1, equations (1.1)--(1.12)]{FHW} once the
root vectors of height $h$ there are multiplied by $(-1)^{h-1}$.

By the quantum PBW theorem \cite[Corollary~40.2.2]{Lusztig},
$M(\mathbf u)$ has the ordinary PBW basis
\[
 w_{\mathbf N}=\prod_{(a,b)}F_{ab}^{N_{ab}}v_0.
\]
We first compute matrix coefficients in this basis at generic $q$ and
then specialize them at roots of unity.

For an interval $\beta$ put
\[
 \Omega_i(\beta;\mathbf N)
   =\sum_{\gamma>\beta}(\alpha_i,\gamma)N_\gamma,
 \qquad
 P_i(\mathbf N)=\sum_{a<i}(N_{a,i}-N_{a,i+1}).
\]
We write $\mathbf e_{ab}$ for the unit exponent vectors, and we use
the conventions that a basis vector with a negative exponent is zero
and that empty sums vanish.

\begin{lemma}
\label{lem:braid-simple-generators}
One has
\begin{align}
 F_iw_{\mathbf N}
 ={}&q^{P_i(\mathbf N)}w_{\mathbf N+\mathbf e_{i,i+1}}
 +\sum_{a<i}q^{\sum_{c<a}(N_{c,i}-N_{c,i+1})}[N_{a,i}]_q
              w_{\mathbf N-\mathbf e_{a,i}+\mathbf e_{a,i+1}},
 \label{eq:braid-F-simple}\\
 E_iw_{\mathbf N}
 ={}&\frac{[N_{i,i+1}]_q}{d}
 \left(u_iq^{1-N_{i,i+1}-\Omega_i(\alpha_i;\mathbf N)}
       -u_i^{-1}q^{N_{i,i+1}-1+\Omega_i(\alpha_i;\mathbf N)}\right)
                 w_{\mathbf N-\mathbf e_{i,i+1}}\notag\\
 &+\sum_{a<i}u_i^{-1}q^{\Omega_i(\alpha_{a,i+1};\mathbf N)}
              [N_{a,i+1}]_q
                 w_{\mathbf N+\mathbf e_{a,i}-\mathbf e_{a,i+1}}\notag\\
 &-\sum_{b>i+1}u_i
  q^{2-N_{i,b}-\sum_{c>b}N_{i,c}+\sum_{c\ge b}N_{i+1,c}}
              [N_{i,b}]_q
                 w_{\mathbf N-\mathbf e_{i,b}+\mathbf e_{i+1,b}}.
 \label{eq:braid-E-simple}
\end{align}
Also $K_iw_{\mathbf N}=u_iq^{-(\alpha_i,\sum_\beta N_\beta\beta)}
w_{\mathbf N}$.
\end{lemma}

\begin{proof}
Induction on the length of the interval, using the defining recursions
and the quantum Serre relations, gives
\begin{equation}
 [E_i,F_{i,b}]=-qF_{i+1,b}K_i\quad(b>i+1),\qquad
 [E_i,F_{a,i+1}]=F_{a,i}K_i^{-1}\quad(a<i).
 \label{eq:braid-interval-endpoints}
\end{equation}
The commutators of $E_i$ with all other nonsimple interval root
vectors vanish. The only possible interior contribution in the
induction occurs for $i=b-2$ in
$F_{a,b}=F_{b-1}F_{a,b-1}-qF_{a,b-1}F_{b-1}$. There the two
resulting terms $F_{b-1}F_{a,b-2}K_i^{-1}$ and
$qF_{a,b-2}K_i^{-1}F_{b-1}$ coincide and cancel, because the two
negative interval vectors commute and
$K_i^{-1}F_{b-1}=q^{-1}F_{b-1}K_i^{-1}$. The outer endpoints give the
two identities above, and in the disjoint cases the operators commute
at every step of the induction.

To compute $F_iw_{\mathbf N}$, we move $F_i$ through the preceding rows
of the PBW monomial. The only relation that produces a new term is
$F_iF_{a,i}=qF_{a,i}F_i+F_{a,i+1}$, whose version for powers has
transfer coefficient $[N_{a,i}]_q$. The new interval $(a,i+1)$ follows
$(a,i)$ in its row, and moving it past $F_{a,i+1}^{N_{a,i+1}}$ in the
main term contributes $q^{-N_{a,i+1}}$; this accounts for the prefix
exponents in~\eqref{eq:braid-F-simple}.

For $E_i$, the contraction with a power of the simple root vector
follows from
\[
 E_iF_i^N=F_i^NE_i+
 \frac{[N]_q}{d}F_i^{N-1}
       \bigl(q^{1-N}K_i-q^{N-1}K_i^{-1}\bigr).
\]
At a right endpoint, commuting $E_i$ through a power first gives
$q^{N-1}[N]_qF_{a,i+1}^{N-1}F_{a,i}K_i^{-1}$.
Moving $F_{a,i}$ to the left cancels $q^{N-1}$, and moving $K_i^{-1}$
through the remaining factors gives the suffix $\Omega_i$. At a left
endpoint the commutator with a power is
$-q^{2-N}[N]_qF_{i,b}^{N-1}F_{i+1,b}K_i$; moving $K_i$ to the right
contributes $q^{-\Omega_i(\alpha_{i,b};\mathbf N)}$, and moving
$F_{i+1,b}$ into its row contributes
$q^{-\sum_{i+1<c<b}N_{i+1,c}}$. Since
\[
 \Omega_i(\alpha_{i,b};\mathbf N)
 =\sum_{c>b}N_{i,c}-\sum_{c>i+1}N_{i+1,c},
\]
the product of these factors is the last exponent
in~\eqref{eq:braid-E-simple}. These are all the possible contractions,
so there are no further terms.
\end{proof}

\subsection{Divided powers at roots of unity}

We write
$\left[\begin{smallmatrix}N\\t\end{smallmatrix}\right]_q
=[N]_q!/([t]_q![N-t]_q!)$
for the symmetric Gaussian coefficient, which is zero unless
$0\le t\le N$.

\begin{lemma}
\label{lem:braid-normal-order}
For $k,N\ge0$ and the indicated intervals,
\begin{align}
 E_i^{(k)}F_{i,b}^{N}
 ={}&\sum_{t=0}^{\min(k,N)}(-1)^tq^{t(2t+1-N-k)}
 \begin{bmatrix}N\\t\end{bmatrix}_q
 F_{i,b}^{N-t}F_{i+1,b}^{t}K_i^tE_i^{(k-t)},
 \quad b>i+1,\label{eq:braid-normal-left}\\
 E_i^{(k)}F_{a,i+1}^{N}
 ={}&\sum_{t=0}^{\min(k,N)}q^{t(k-t)}
 \begin{bmatrix}N\\t\end{bmatrix}_q
 F_{a,i}^{t}F_{a,i+1}^{N-t}K_i^{-t}E_i^{(k-t)},
 \quad a<i.\label{eq:braid-normal-right}
\end{align}
At the simple root,
\begin{equation}\label{eq:braid-normal-simple}
 E_i^{(k)}F_i^N
 =\sum_{t=0}^{\min(k,N)}\begin{bmatrix}N\\t\end{bmatrix}_q
 F_i^{N-t}
 \prod_{j=1}^{t}
 \frac{K_iq^{2t-k-N-j+1}-K_i^{-1}q^{-2t+k+N+j-1}}d
 E_i^{(k-t)}.
\end{equation}
\end{lemma}

\begin{proof}
For the left endpoint put $A=E_i$, $B=F_{i,b}$ and
$C=F_{i+1,b}K_i$, so that
$AB=BA-qC$, $AC=q^{-2}CA$ and $CB=q^{-2}BC$.
Induction with the Gaussian Pascal rule shows that the coefficient of
$B^{N-t}C^tA^{k-t}$ in $A^kB^N$ is
\[
 (-q)^t \binom{k}{t}_{q^{-2}}
   \prod_{j=0}^{t-1}(N-j)_{q^{-2}},
\]
where $(N)_v=(1-v^N)/(1-v)$ and $\binom{k}{t}_v$ is the corresponding
Gaussian polynomial; these coefficients satisfy
\[
 a_{k+1,N,t}=q^{-2t}a_{k,N,t}
              -q(N-t+1)_{q^{-2}}a_{k,N,t-1},
 \qquad a_{0,N,t}=\delta_{t0}.
\]
Converting to symmetric brackets, using
$C^t=q^{t(t-1)/2}F_{i+1,b}^tK_i^t$ and dividing by $[k]_q!$ gives
\eqref{eq:braid-normal-left}. For the right endpoint take
$C=F_{a,i}K_i^{-1}$. The relations become $AB=BA+C$, $AC=q^2CA$ and
$CB=q^2BC$, and the same induction, followed by moving $F_{a,i}^t$ to
the left, gives \eqref{eq:braid-normal-right}. Formula
\eqref{eq:braid-normal-simple} follows by induction from
$[E_i,F_i]=(K_i-K_i^{-1})/d$ and $K_iE_i=q^2E_iK_i$; the two
contributions with $t$ contractions are combined by the symmetric
Gaussian Pascal rule.
\end{proof}

Fix $m>1$ and a primitive $m$th root $\zeta$, and choose
\[
 \xi^2=\zeta,\quad \xi^m=\chi=(-1)^{m+1},\quad
 Q=\zeta(1+\epsilon),\qquad q=\xi(1+\epsilon)^{1/2}.
\]
Such a $\xi$ exists: for odd $m$ one of the two square roots of $\zeta$
has $\xi^m=1$, and for even $m$ both satisfy $\xi^m=\zeta^{m/2}=-1$.
Then $d(\xi)\ne0$ and $\xi^{m^2}=1$. Put
$S=\mathbb Q(\xi)[u_1^{\pm1},\ldots,u_{n-1}^{\pm1}]$ and
$N_\alpha=mA_\alpha+r_\alpha$, $0\le r_\alpha<m$.

For the rest of the paper, $\epsilon$ denotes the multiplicative local
coordinate defined by $Q=\zeta(1+\epsilon)$. It differs from the additive
coordinate $x=Q-\zeta=\zeta\epsilon$ of Sections~\ref{sec:newton}
and~\ref{sec:relative} by the unit factor $\zeta$, which changes
neither orders of vanishing nor integrality.

\begin{proposition}
\label{prop:braid-regular-jets}
For every positive root $\alpha$ and every fixed $k\ge0$, all ordinary
PBW matrix coefficients of $E_\alpha^{(k)}$ belong to
\[
 \mathbb Z[q^{\pm1},\mathbf u^{\pm1},(q-q^{-1})^{-1}].
\]
In particular they are regular at the roots of unity fixed above.
Only finitely many exponent shifts occur, independently of the input
exponents. On a fixed residue class, each quantum coefficient of order
$h$ is a finite sum of sign characters in $\mathbf A$ times
polynomials of total degree at most
\[
 h+\left\lfloor\frac{k}{m}\right\rfloor.
\]
The same holds for a fixed ordinary power $F_\alpha^k$ of a negative
root vector, without the quantum factorial in the denominator and with
degree bound $h$.
\end{proposition}

\begin{proof}
Put
\[
 \mathscr R=\Z[q^{\pm1},\mathbf u^{\pm1},d^{-1}],
 \qquad d=q-q^{-1}.
\]
We first prove regularity for a fixed PBW input. Move $E_i^{(k)}$
through the root factors of the PBW monomial using
Lemma~\ref{lem:braid-normal-order}. Every contraction has coefficient in
$\mathscr R$: its Gaussian factor is a Laurent polynomial, and its
toral factors introduce only powers of $d^{-1}$. The negative words
produced in this way are put back into PBW order by expanding the root
vectors into simple generators and applying \eqref{eq:braid-F-simple},
which introduces no denominators.

For composite roots we use the divided-power identity
\begin{equation}\label{eq:braid-root-divided}
 E_{a,b}^{(k)}
 =\sum_{t=0}^k(-1)^tq^{-t}
 E_{a,b-1}^{(k-t)}E_{b-1}^{(k)}E_{a,b-1}^{(t)}.
\end{equation}
For the proof, put $X=E_{a,b-1}$ and $Y=E_{b-1}$. These satisfy the
rank-two Serre relations, which for $Z=XY-q^{-1}YX$ give
\[
 XZ=qZX,\qquad ZY=qYZ;
\]
induction on $k$, using the Gaussian Pascal identity, then gives
\eqref{eq:braid-root-divided}. The required Serre relations follow by
applying $T_a\cdots T_{b-3}$ to the adjacent simple roots, with
$T_i(E_j)=E_iE_j-q^{-1}E_jE_i$ for $|i-j|=1$.
Induction on $b-a$ now proves regularity for all positive roots.
The assertion for negative roots follows directly from
\eqref{eq:braid-F-simple}.

For bounds that are uniform in the input, write
\[
 E_\alpha^{(k)}=\frac{E_\alpha^k}{[k]_q!}.
\]
The numerator is a finite linear combination of words of length
$k\hgt(\alpha)$ in the simple positive generators. By
\eqref{eq:braid-E-simple}, each generator has a fixed finite set of
exponent shifts of $\ell^1$-norm at most two. Consequently every
resulting shift $s$ satisfies
\[
 \lVert s\rVert_1\le2k\hgt(\alpha).
\]
After expanding the quantum integers, the coefficient of a path of
shifts is a finite sum of terms $b(q,\mathbf u)q^{\ell(\mathbf N)}$,
where $b\in\mathscr R$ and $\ell$ is affine linear in the input
exponents $\mathbf N$; the intermediate shifts change only the
constant term of $\ell$. Write $\ell_0$ for the linear part of $\ell$.
On the residue class $\mathbf N=m\mathbf A+\mathbf r$,
\[
 q^{\ell(m\mathbf A+\mathbf r)}
 =\xi^{\ell(\mathbf r)}\chi^{\ell_0(\mathbf A)}
 (1+\epsilon)^{\ell(m\mathbf A+\mathbf r)/2},
\]
so, apart from the sign character $\chi^{\ell_0(\mathbf A)}$, the
coefficient of $\epsilon^j$ is a polynomial in $\mathbf A$ of degree at
most $j$.

The factors of $[k]_q!$ that vanish at $q=\xi$ are exactly those whose
index is divisible by $m$, and each has a simple zero. Thus
\[
 [k]_q!=\epsilon^v u(\epsilon),\qquad
 v=\left\lfloor\frac{k}{m}\right\rfloor,
 \qquad u(0)\ne0.
\]
By regularity the quotient has no pole, and its coefficient of
$\epsilon^h$ is determined by the numerator coefficients of order at
most $h+v$. This proves the degree bound $h+v$. For $F_\alpha^k$
there is no factorial denominator, and the same argument gives degree
at most $h$.
\end{proof}

\subsection{Polynomial differential operators}

Define the coarse and residue root weights
\[
 \gamma=\sum_\alpha A_\alpha\alpha,
 \qquad \varrho=\sum_\alpha r_\alpha\alpha,
\]
and the local regrouping exponent
\[
 L(\mathbf A,\mathbf r)
  =\sum_{\beta<\alpha}(\alpha,\beta)A_\alpha r_\beta.
\]
On $b$ strands we make the diagonal change of basis
\begin{equation}\label{eq:braid-global-gauge}
 \mathcal G
 =\chi^{\sum_jL(\mathbf A_j,\mathbf r_j)
       +\sum_j(\gamma_j,\sum_{i<j}\varrho_i)}.
\end{equation}
This is the identity for odd $m$ and an involution for even $m$.
We encode the quotient indices by the monomials
$\prod_{\alpha,j}X_{\alpha,j}^{A_{\alpha,j}}$ in independent
polynomial variables. The residue factor has the finite basis
$e_{\mathbf r_1}\otimes\cdots\otimes e_{\mathbf r_b}$.
At $\epsilon=0$, the first sum in the exponent of $\mathcal G$
accounts for moving the $m$th powers of root vectors ahead of the
residue factors. The second sum cancels the signs that remain when a
root transfer acts on two adjacent strands.

\begin{theorem}
\label{thm:braid-fixed-transfer-jets}
For every positive root $\alpha$, fixed $k\ge0$, adjacent pair of
strands, and quantum order $h\ge0$, the operator
\[
 [\epsilon^h]\,\mathcal G
       (E_\alpha^{(k)}\otimes F_\alpha^k)\mathcal G^{-1}
\]
is a finite matrix of polynomial differential operators in the coarse
variables, with coefficients in $S$, of differential order at most
\[
 h+\left\lfloor\frac{k}{m}\right\rfloor.
\]
The same statement holds after multiplication by any scalar in
$S[[\epsilon]]$. The change of basis does not depend on the crossing,
on $k$ or on $h$.
\end{theorem}

\begin{proof}
The proof has two steps. First we show that after conjugation by
$\mathcal G$ no sign characters remain, so that for each fixed
exponent shift the coefficient on a residue class is a polynomial in
the coarse indices. Then we show that such a polynomial coefficient is
realized by a polynomial differential operator.

Since $\chi=1$ for odd $m$, only even $m$ requires a sign calculation.
For an ordinary exponent shift $s$ with coarse shift $b$, the new
residue vector is $\mathbf r+s-mb$. Hence
\begin{align}
 &L(\mathbf A+b,\mathbf r+s-mb)-L(\mathbf A,\mathbf r)\notag\\
 &\qquad\equiv
 \sum_{\beta<\alpha}(\alpha,\beta)
 \bigl(A_\alpha s_\beta+b_\alpha(r_\beta+s_\beta)\bigr)
 \pmod2.
 \label{eq:braid-local-carry}
\end{align}
The part depending on $\mathbf A$ is the linear form
\[
 \ell_s(\mathbf A)
 =\sum_\gamma A_\gamma
          \sum_{\beta<\gamma}(\gamma,\beta)s_\beta.
\]
The three types of shifts in \eqref{eq:braid-E-simple} satisfy
\begin{align*}
 \ell_{-\mathbf e_{i,i+1}}(\mathbf A)
 &\equiv\sum_{\gamma>\alpha_i}(\alpha_i,\gamma)A_\gamma,\\
 \ell_{\mathbf e_{a,i}-\mathbf e_{a,i+1}}(\mathbf A)
 &\equiv A_{a,i+1}
       +\sum_{\gamma>(a,i+1)}(\alpha_i,\gamma)A_\gamma,\\
 \ell_{-\mathbf e_{i,b}+\mathbf e_{i+1,b}}(\mathbf A)
 &\equiv\sum_{c>b}A_{i,c}+\sum_{c\ge b}A_{i+1,c}
 \pmod2.
\end{align*}
The right sides are exactly the sign characters
$\chi^{\ell_0(\mathbf A)}$ that the coefficients of the simple-root
contraction, the right-endpoint transfer and the left-endpoint transfer
in \eqref{eq:braid-E-simple} acquire on a residue class (see the proof of
Proposition~\ref{prop:braid-regular-jets}). The identities follow by
substituting $(\alpha_{a,b},\alpha_{c,d})=
\delta_{ac}-\delta_{ad}-\delta_{bc}+\delta_{bd}$ in the fixed PBW
order. In the second identity the two intervals are consecutive in
their row. In the third, the intermediate roots contribute the entries
of row $i$ after column $b$ together with $(i+1,b)$, and the later
roots contribute the entries of row $i+1$ after column $b$. Thus the
local part of the gauge cancels every coarse character of a positive
generator.

For $F_i$, the prefix in \eqref{eq:braid-F-simple} and the local gauge
together leave the character $\chi^{(\alpha_i,\gamma)}$. The same
character occurs for the main term and for every transfer term,
because the additional $A_{a,i}$ in the local gauge supplies the
coefficient that is missing at the transferred interval. Under
composition, an intermediate coarse shift changes a character only by
a constant factor. Therefore
\[
 \begin{array}{c|c}
 \text{homogeneous word}&\text{remaining coarse character}\\ \hline
 \text{positive}&1\\
 \text{negative of root weight }\nu&\chi^{(\nu,\gamma)}.
 \end{array}
\]
Scalar quantum factorials do not affect these characters, and
Proposition~\ref{prop:braid-regular-jets} allows us to specialize after
division.

For a paired transfer of weight $\nu=k\alpha$, let $b_1,b_2$ be the
coarse root-weight shifts on the two strands. Then
\[
 \Delta\varrho_1=-\nu-mb_1,\qquad
 \Delta\varrho_2=\nu-mb_2.
\]
If $R$ is the sum of the residue weights on the preceding strands,
the change of the prefix exponent in \eqref{eq:braid-global-gauge} is
\begin{align*}
 \Delta P={}&(b_1+b_2,R)+(\gamma_2,\Delta\varrho_1)
       +(b_2,\varrho_1+\Delta\varrho_1)\\
 &\quad+\sum_{j>2}(\gamma_j,\Delta\varrho_1+\Delta\varrho_2)\\
 \equiv{}&-(\gamma_2,\nu)+(b_1+b_2,R)
                 +(b_2,\varrho_1-\nu)\pmod2.
\end{align*}
The first term cancels the character $\chi^{(\nu,\gamma_2)}$ of the
negative word on the second strand. The other terms depend only on the
residue transition and on the fixed shift, not on $\mathbf A$. This
proves the cancellation for every adjacent pair, including carries and
borrows between residue and coarse exponents.

After the cancellation, each fixed shift has a coefficient
$p(\mathbf A)$ that is a polynomial of degree at most
$h+\lfloor k/m\rfloor$, by Proposition~\ref{prop:braid-regular-jets}.
To realize this shift by a differential operator we use the vanishing
of $p$ at the boundary of the nonnegative exponent lattice.

Suppose the shift lowers $A_i$ by $d_i>0$. For
$A_i=0,\ldots,d_i-1$ the resulting ordinary PBW exponent would be
negative. Every path to such an exponent lowers an exponent from $0$
at the step where it first leaves the lattice, and the coefficient of
that step contains $[0]_q=0$. The coefficient is therefore
identically zero before specialization, and it remains zero after
division by a scalar quantum factorial and specialization. Fixing one
of these values of $A_i$ and letting the other coordinates run through
sufficiently large integers shows that $p|_{A_i=s}=0$ as a polynomial
for each $0\le s<d_i$. The factor theorem, applied successively in
the different variables, gives
\[
 \prod_i(A_i)_{d_i}\mid p(\mathbf A),
 \qquad (A_i)_{d_i}=A_i(A_i-1)\cdots(A_i-d_i+1).
\]

For a shift $b=b_+-b_-$ with $b_\pm\ge0$, write
$p=(\mathbf A)_{b_-}\widetilde p$. On the monomial $X^{\mathbf A}$,
the operator
\[
 X^{b_+}\partial^{b_-}\widetilde p(\boldsymbol\Theta),
 \qquad \Theta_i=X_i\partial_{X_i},
\]
acts by $p(\mathbf A)X^{\mathbf A+b}$. Its differential order is at
most
\[
 |b_-|+\deg\widetilde p\le\deg p
             \le h+\left\lfloor\frac{k}{m}\right\rfloor.
\]
There are only finitely many shifts. Multiplication by a scalar
series only combines coefficients of orders at most $h$, so it
preserves this bound.
\end{proof}

\subsection{Root powers and polynomial flows}
\label{subsec:braid-principal-fields}

We keep $m,\zeta,\xi,\chi$ and the interval root vectors fixed above,
and at $q=\xi$ we put
\[
 c=\xi^{m(m-1)/2}(\xi-\xi^{-1})^m,\qquad
 X_{a,b}=F_{a,b}^m,\qquad
 U_{a,b}=\prod_{i=a}^{b-1}u_i^m,\quad
 S_{a,b}=U_{a,b}-U_{a,b}^{-1}.
\]
By Proposition~\ref{prop:braid-regular-jets}, every divided-power
operator on the ordinary PBW module can be specialized at $q=\xi$.

\begin{lemma}\label{lem:braid-endpoint-root-power}
At $q=\xi$, the variables $X_{a,b}$ commute and are algebraically
independent. For every positive root $\gamma$,
\[
 F_\gamma X_\alpha
 =\chi^{(\gamma,\alpha)}X_\alpha F_\gamma,
 \qquad [E_i,X_\alpha]=0.
\]
The only nonzero commutators of the $cE_i^{(m)}$ with the $X_\alpha$
are
\begin{align}
 [cE_i^{(m)},X_{i,i+1}]&=K_i^m-K_i^{-m},
 \label{eq:braid-simple-center}\\
 [cE_i^{(m)},X_{i,b}]&=-cX_{i+1,b}K_i^m\quad(b>i+1),\notag\\
 [cE_i^{(m)},X_{a,i+1}]&=cX_{a,i}K_i^{-m}\quad(a<i).\notag
\end{align}
\end{lemma}

\begin{proof}
For a negative simple generator and an interval root vector, the
straightening relations give commutation, $q$-commutation at a shared
endpoint, or a relation at adjacent endpoints. The last is a rank-two
relation $AB=q^{\pm1}BA+C$ in which $B$ and $C$ $q$-commute. When $B$
is replaced by $B^m$, the correction term is divisible by $[m]_q$, so
it vanishes at $q=\xi$ and only the sign $\chi^{(\alpha_i,\alpha)}$
remains. The defining recursion shows that these are all the cases: a
simple interval in the interior of $\alpha$ is nested and commutes, a
disjoint nonadjacent one commutes, and the four endpoint positions
give the two shared-endpoint and the two adjacent-endpoint relations.
Passing through a homogeneous word in the $F_i$ gives the formula for
$F_\gamma$, and hence
$X_\gamma X_\alpha=\chi^{m(\gamma,\alpha)}X_\alpha X_\gamma
=X_\alpha X_\gamma$. Algebraic independence follows from the ordinary
PBW basis. The positive endpoint identities
\eqref{eq:braid-interval-endpoints} and the rank-one formula produce a
factor $[m]_\xi$ in every commutator of $E_i$ with an $m$th power, so
$[E_i,X_\alpha]=0$.

For the divided powers, apply the endpoint normal-order formulas with
$N=k=m$. Each intermediate contraction contains the factor
$\begin{bmatrix}m\\t\end{bmatrix}_\xi=0$, $0<t<m$, while the remaining
operator is regular, so these terms vanish after specialization. The
full contraction at the left endpoint has scalar $(-1)^m\xi^m=-1$, and
the one at the right endpoint has scalar one, which gives the last two
lines of~\eqref{eq:braid-simple-center}. At the simple root the full
contraction is
\[
 c\prod_{j=0}^{m-1}
 \frac{K_i\xi^{-j}-K_i^{-1}\xi^j}{\xi-\xi^{-1}}
 =K_i^m\prod_{j=0}^{m-1}(1-K_i^{-2}\zeta^j)
 =K_i^m-K_i^{-m}.
\]
\end{proof}

The commutators in Lemma~\ref{lem:braid-endpoint-root-power} lead to the
following derivations of the polynomial ring
$\mathbb Q(\xi)[u^{\pm1}][X_{a,b}:a<b]$. For each positive root
$\alpha_{a,b}$ put
\begin{equation}\label{eq:braid-affine-root-field}
 D_{a,b}=S_{a,b}\partial_{X_{a,b}}
 -cU_{a,b}\sum_{j>b}X_{b,j}\partial_{X_{a,j}}
 +cU_{a,b}^{-1}\sum_{k<a}X_{k,a}\partial_{X_{k,b}}.
\end{equation}
On two adjacent strands, with superscripts indicating the strand, we
also put $A_{a,b}=X_{a,b}^{(2)}D_{a,b}^{(1)}$. The fields $D_{a,b}$
are linear in a convenient matrix coordinate, which makes their flows
explicit.

\begin{lemma}
\label{lem:braid-matrix-field}
Let $\mathsf X$ be the strictly upper triangular matrix with entries
$X_{a,b}$, and put $Y=I-c\mathsf X$. For the ordinary matrix unit
$\mathsf e_{a,b}$,
\begin{equation}\label{eq:braid-matrix-field}
 D_{a,b}Y=-cU_{a,b}\mathsf e_{a,b}Y
             +cU_{a,b}^{-1}Y\mathsf e_{a,b},
\end{equation}
where $D_{a,b}$ acts on $Y$ entrywise.
Consequently, whenever both index pairs are positive roots,
\[
 [D_{a,b},D_{a',b'}]
 =c\bigl(\delta_{b,a'}D_{a,b'}-\delta_{a,b'}D_{a',b}\bigr).
\]
Both $e^{sD_{a,b}}$ and its inverse are polynomial in $s$ and the
coordinates jointly, and
\begin{equation}\label{eq:braid-matrix-flow}
 e^{sD_{a,b}}Y
 =(I-scU_{a,b}\mathsf e_{a,b})Y
                    (I+scU_{a,b}^{-1}\mathsf e_{a,b}).
\end{equation}
The same holds for the paired fields $A_{a,b}$.
\end{lemma}

\begin{proof}
Compare entries. The diagonal entries of $Y$ are one, the entry
$(a,b)$ of \eqref{eq:braid-matrix-field} gives
$-c(U_{a,b}-U_{a,b}^{-1})$, the entries $(a,j)$ with $j>b$ and $(k,b)$
with $k<a$ give the other two sums in \eqref{eq:braid-affine-root-field},
and all other entries vanish. Since $D_{a,b}$ is a derivation and the
matrices $\mathsf e_{a,b}$ are constant, applying the identity twice
gives
\[
 [D_{a,b},D_{a',b'}]Y
 =-c^2U_{a,b}U_{a',b'}[\mathsf e_{a,b},\mathsf e_{a',b'}]Y
  +c^2U_{a,b}^{-1}U_{a',b'}^{-1}
             Y[\mathsf e_{a,b},\mathsf e_{a',b'}];
\]
the mixed terms cancel. Using
$[\mathsf e_{a,b},\mathsf e_{a',b'}]
=\delta_{b,a'}\mathsf e_{a,b'}-\delta_{a,b'}\mathsf e_{a',b}$
and, for $b=a'$, the identity $U_{a,b}U_{b,b'}=U_{a,b'}$ (the other
case is analogous), we obtain $c$ times the right side of
\eqref{eq:braid-matrix-field} for $D_{a,b'}$, respectively for $D_{a',b}$.
A derivation of the polynomial ring is determined by its values on the
entries of $Y$, which proves the bracket formula.

The matrix $Y_s=e^{sD_{a,b}}Y$ solves the linear differential
equation $\frac{d}{ds}Y_s=-cU_{a,b}\mathsf e_{a,b}Y_s+cU_{a,b}^{-1}
Y_s\mathsf e_{a,b}$ with initial value $Y$, and the right side of
\eqref{eq:braid-matrix-flow} is its solution. Since
$\mathsf e_{a,b}^2=0$ and
$\mathsf e_{a,b}Y\mathsf e_{a,b}=Y_{b,a}\mathsf e_{a,b}=0$, its entries
are linear in $s$. Replacing $s$ by $-s$ gives the inverse. For the
paired field, $X_{a,b}^{(2)}$ is a constant for the derivation
$D^{(1)}_{a,b}$ on the first strand, so replacing $s$ by
$sX_{a,b}^{(2)}$ gives the joint polynomiality of $e^{\pm sA_{a,b}}$.
\end{proof}

From now on we use the grouped basis $X^A e_{\mathbf r}$, where
$0\le r_\alpha<m$ and $e_{\mathbf r}$ is the ordered PBW monomial in
the residues, and we set
\[
 H(\mathbf r)=\sum_\alpha\operatorname{ht}(\alpha)r_\alpha.
\]
Let $V$ be the free residue module with basis $e_{\mathbf r}$,
filtered by
$V_{\le h}=\operatorname{span}_S\{e_{\mathbf r}:H(\mathbf r)\le h\}$.
On $S[\mathbf X]\otimes V$ the filtration only involves the residue
factor, and on tensor products we use the sum of the residue heights.
The change from the ordinary PBW basis to the grouped basis involves
the local sign
$\chi^{\sum_{\beta<\alpha}(\alpha,\beta)A_\alpha r_\beta}$.

\begin{proposition}\label{prop:braid-all-root-fields}
On the graded module associated with the residue-height filtration,
$cE_{a,b}^{(m)}$ acts on residue $\mathbf r$ by
$\chi^{(\alpha_{a,b},\varrho)}D_{a,b}$. The remaining part lowers the
residue height by at least $m$. After conjugation by $\mathcal G$,
every paired root operator has the form
\begin{equation}\label{eq:braid-paired-principal}
 cE_{a,b}^{(m)}\otimes F_{a,b}^{m}
 =A_{a,b}I+N_{a,b},\qquad
 A_{a,b}=X_{a,b}^{(2)}D_{a,b}^{(1)},
\end{equation}
where $N_{a,b}$ lowers the total residue height by at least $m$.
Each small paired power $E_\alpha^{(r)}\otimes F_\alpha^r$,
$0\le r<m$, preserves the filtration and acts on the associated graded
residue space by a constant matrix.
\end{proposition}

\begin{proof}
Put $P=S[\mathbf X]$. For a simple root,
Lemma~\ref{lem:braid-endpoint-root-power} and the Leibniz rule give
\[
 cE_i^{(m)}(p e_{\mathbf r})
 \equiv\chi^{(\alpha_i,\varrho)}D_{i,i+1}(p)e_{\mathbf r}
 \pmod{P\otimes V_{\le H(\mathbf r)-m}},
\]
because $K_i^{\pm m}$ acts on $e_{\mathbf r}$ by
$\chi^{(\alpha_i,\varrho)}u_i^{\pm m}$. The omitted term is
$p\,cE_i^{(m)}e_{\mathbf r}$. If it contains $X^Be_{\mathbf s}$, then
conservation of ordinary root height gives
\[
 H(\mathbf s)
 =H(\mathbf r)-m
       -m\sum_\alpha\hgt(\alpha)B_\alpha
 \le H(\mathbf r)-m.
\]
This proves both the formula and the filtration estimate for simple
roots.

For $b>a+1$, apply \eqref{eq:braid-root-divided} with $k=m$:
\[
 E_{a,b}^{(m)}
 =\sum_{r=0}^{m}(-1)^rq^{-r}
 E_{a,b-1}^{(m-r)}E_{b-1}^{(m)}E_{a,b-1}^{(r)}.
\]
When $0<r<m$, the outer divided powers have scalar factorials that are
nonzero at $\xi$, so they commute with the coarse algebra. Together
they lower the residue height by at least $m\hgt(\alpha_{a,b-1})$;
factoring out newly created coarse powers can only lower it further.
The middle factor preserves the filtration. Hence every intermediate
term lowers the residue height by at least $m$.

The two endpoint terms preserve each residue label in their principal
parts. Since $(-1)^m\xi^{-m}=-1$, their sum, multiplied by $c$, has
principal part
\[
 \chi^{(\alpha_{a,b},\varrho)}c^{-1}
       [D_{a,b-1},D_{b-1,b}]
 =\chi^{(\alpha_{a,b},\varrho)}D_{a,b},
\]
where the equality is the bracket formula of
Lemma~\ref{lem:braid-matrix-field}. Every term involving a remainder
still lowers the filtration by at least $m$, which completes the
induction on $b-a$.

A principal paired transfer changes the coarse root weights of the
two strands by $(-\alpha,+\alpha)$ and leaves the residues fixed. The
prefix exponent of $\mathcal G$ therefore changes by
$(\alpha,\varrho_1)$, which cancels the toral sign, while the negative
$m$th power acts by multiplication by $X_\alpha^{(2)}$. Hence the
principal paired action is $X_\alpha^{(2)}D_\alpha^{(1)}$, as in
\eqref{eq:braid-paired-principal}.

For $0\le r<m$, move the input coarse monomials past the small
paired power. The positive factor commutes with them, and the gauge
cancels the sign from the negative factor. On bare residues, write an
output term as $X^{B_1}e_{\mathbf s_1}\otimes X^{B_2}e_{\mathbf s_2}$.
The coarse exponents are nonnegative, and
\[
 H(\mathbf s_1)+H(\mathbf s_2)
 =H(\mathbf r_1)+H(\mathbf r_2)
   -m\sum_\alpha\hgt(\alpha)(B_{1,\alpha}+B_{2,\alpha}).
\]
A term that preserves the height therefore has $B_1=B_2=0$ and acts
only on the residue labels, with coefficients in $S$, while every other
term lowers the height by at least $m$. For $r>0$, a nonzero graded
term lowers the residue height of the first strand by
$r\hgt(\alpha)$. Summing the graded small powers with their
coefficients in \eqref{eq:braid-normalized-crossing} therefore gives the
identity plus a nilpotent residue matrix. This sum is unipotent, and
its inverse is a finite geometric sum.
\end{proof}

\begin{corollary}\label{cor:braid-uniform-residue-bound}
Let
\[
 H_n=\sum_{1\le a<b\le n}(b-a)=\frac{n(n^2-1)}6.
\]
On two strands the residue-height filtration has width
$2(m-1)H_n$. A product of more than
\[
 L_{n,m}=\left\lfloor\frac{2(m-1)H_n}{m}\right\rfloor
 \le 2H_n-1
\]
strictly lowering remainders is zero, and the same holds for products
of their conjugates by the scalar principal flows, with independent
flow parameters.
\end{corollary}

\begin{proof}
Each remainder lowers the integer-valued filtration by at least $m$.
The scalar flows act on the coarse coordinates and trivially on the
residue factor, so their conjugates do the same. The bound follows
from the width of the filtration, and $2H_n-1$ bounds it uniformly in
$m$ because $2(m-1)H_n/m<2H_n$.
\end{proof}

\subsection{The crossing at a root of unity}
\label{subsec:braid-principal-crossing}

To analyze a crossing we group the terms of each root factor of
\eqref{eq:braid-normalized-crossing} according to the residue of $k$
modulo $m$. For a fixed positive root put
\[
 c_k(q)=q^{k(k-1)/2}(q-q^{-1})^k,\qquad
 V_r=c_r(q)E_\alpha^{(r)}\otimes F_\alpha^r\quad(0\le r<m),
 \qquad Z=c_m(q)E_\alpha^{(m)}\otimes F_\alpha^m.
\]
We always conjugate these operators by $\mathcal G$.

For $\kappa\ge0$ let $M_\kappa$ be the span of the tensor PBW monomials
of total ordinary root height $\kappa$. It is finite-dimensional and is
preserved by root transfers, Cartan factors, flips and inverse
crossings. A transfer with index $k$ vanishes on $M_\kappa$ when
$k\operatorname{ht}(\alpha)>\kappa$; in particular $Z(\epsilon)^a=0$ on
$M_\kappa$ when $ma\operatorname{ht}(\alpha)>\kappa$. All operator
identities below are first established on these finite-dimensional
spaces.

\begin{lemma}
\label{lem:braid-scalar-grouping}
On every block $M_\kappa$,
\[
 \sum_{k\ge0}c_k(q)E_\alpha^{(k)}\otimes F_\alpha^k
 =\sum_{r=0}^{m-1}V_r
     \sum_{a\ge0}A_{m,r,a}(Q)\frac{Z^a}{a!},
\]
where
\[
 A_{m,r,a}(Q)=a!Q^{m^2a(a-1)/2+mar}
 \frac{(Q;Q)_m^a(Q;Q)_r}{(Q;Q)_{ma+r}}.
\]
This scalar is regular at $Q=\zeta$ with value one there, and for
every $h\ge0$ the coefficient of $\epsilon^h$ in
$A_{m,r,a}(\zeta(1+\epsilon))$ and in its reciprocal is a polynomial
in $a$ of degree at most $2h$ over $\mathbb Q(\zeta)$.
\end{lemma}

\begin{proof}
The identity for the scalar follows from
\[
 (E^{(m)})^a=\frac{[ma]_q!}{[m]_q!^a}E^{(ma)},\qquad
 E^{(r)}E^{(ma)}=\begin{bmatrix}ma+r\\r\end{bmatrix}_q E^{(ma+r)}
\]
and $[k]_q!=q^{-k(k-1)/2}(1-Q)^{-k}(Q;Q)_k$; the computation does not
depend on the rank or on the positive root.

For the remaining assertions put $z=\log(1+\epsilon)$, so that
$Q=\zeta e^z$, and
\[
 G(z)=\frac{e^z-1}{z},\qquad
 H_j(z)=\frac{1-\zeta^je^z}{1-\zeta^j}\quad(1\le j<m).
\]
All these series have constant term one. A factor $1-Q^N$ equals
$-NzG(Nz)$ if $m\mid N$, and $(1-\zeta^j)H_j(Nz)$ if $N\equiv j$
modulo $m$ with $1\le j<m$. In the denominator the factors with index
divisible by $m$ give $(-mz)^a a!\prod_{b=1}^aG(mbz)$, which cancels
against the factor $a!$ and the corresponding zeros of the numerator,
and the remaining nonzero constants cancel in complete residue blocks.
Hence $A_{m,r,a}(\zeta)=1$, and
\begin{align*}
 \log A_{m,r,a}(\zeta e^z)
 ={}&\left(\frac{m^2a(a-1)}2+mar\right)z
       +a\log G(mz)-\sum_{b=1}^{a}\log G(mbz)\\
 &+\sum_{j=1}^{m-1}
   \left(a\log H_j(jz)-\sum_{b=0}^{a-1}\log H_j((mb+j)z)\right)\\
 &+\sum_{j=1}^{r}
   \left(\log H_j(jz)-\log H_j((ma+j)z)\right).
\end{align*}
In order $\ell$ in $z$, each of the finite sums is a sum of $\ell$th
powers of affine functions of $b$, hence a polynomial in $a$ of degree
at most $\ell+1$; for $\ell=1$ this includes the quadratic term
displayed first. After exponentiation, a product contributing to order
$h$ has $\sum_\ell\ell n_\ell=h$, and therefore degree at most
$\sum_\ell(\ell+1)n_\ell\le2h$. Changing the sign of the logarithm
gives the statement for the reciprocal, and returning from $z$ to
$\epsilon$ only mixes orders at most $h$.
\end{proof}

For a polynomial map $\Phi$ we write $C_\Phi f=f\circ\Phi$ for the
substitution operator. On a polynomial module with residue
coefficients, $C_\Phi$ acts on the polynomial factor and fixes the
residue labels.

\begin{proposition}\label{prop:braid-triangular-crossing}
After conjugation by $\mathcal G$, the associated graded of the
normalized checked crossing has the form $R^{\rm res}C_\Phi$. In the
coordinates $X_\alpha=(X_\alpha^{(1)},X_\alpha^{(2)})^T$ of a root
$\alpha$, the map $\Phi$ is given by
\[
 \Phi_\alpha(\mathbf X)
 =B(T_\alpha)X_\alpha
     +P_\alpha(X_\beta:\operatorname{ht}(\beta)
                                      <\operatorname{ht}(\alpha)),
\]
where
\[
 T_{a,b}=U_{a,b}^2,\qquad
 B(T)=\begin{pmatrix}1-T^{-1}&T^{-1/2}\\T^{-1/2}&0\end{pmatrix}.
\]
Here $T_{a,b}^{1/2}$ denotes the Laurent monomial $U_{a,b}$. The
coordinate of the root $\alpha_{a,b}$ is homogeneous of root weight
$\alpha_{a,b}$ and of ordinary degree at most $b-a$, and both $\Phi$
and $\Phi^{-1}$ are polynomial over
$\mathbb Q(\xi)[u_1^{\pm1},\ldots,u_{n-1}^{\pm1}]$. The residue
matrix of the associated graded is invertible over the same ring.
Compositions along signed braid words have the same structure, with
the corresponding Burau matrices as linear blocks.
\end{proposition}

\begin{proof}
At $\epsilon=0$ the scalars $A_{m,r,a}$ of
Lemma~\ref{lem:braid-scalar-grouping} equal one, so a root factor becomes
$\sum_rV_r\,e^{Z}$. By Proposition~\ref{prop:braid-all-root-fields}, its
associated graded is a constant residue matrix composed with
$e^{A_{a,b}}$. The matrix flow \eqref{eq:braid-matrix-flow}, with $s$
replaced by $X^{(2)}_{a,b}$, gives the coordinate changes
\begin{align*}
 X_{a,b}^{(1)}&\longmapsto
       X_{a,b}^{(1)}+(U_{a,b}-U_{a,b}^{-1})X_{a,b}^{(2)},\\
 X_{a,j}^{(1)}&\longmapsto
       X_{a,j}^{(1)}-cU_{a,b}X_{b,j}^{(1)}X_{a,b}^{(2)}
                       &&(j>b),\\
 X_{k,b}^{(1)}&\longmapsto
       X_{k,b}^{(1)}+cU_{a,b}^{-1}X_{k,a}^{(1)}X_{a,b}^{(2)}
                       &&(k<a).
\end{align*}
Every nonlinear term replaces an interval by a product of two shorter
ones. It therefore preserves root weight and involves only variables
of lower height, and a monomial of root weight $\alpha_{a,b}$ has
ordinary degree at most $b-a$.

In the normalized Cartan factor, the phase between two coarse parts is
one because $\xi^{m^2}=1$. The mixed phase is
\[
 \chi^{(\gamma_1,\varrho_2)+(\varrho_1,\gamma_2)},
\]
and it is cancelled by the change of the prefix gauge under the flip.
The remaining color dilation and shear on a root block give, with
$U=U_\alpha$,
\[
 U^{-1}
 \begin{pmatrix}1&U-U^{-1}\\0&1\end{pmatrix}
 \begin{pmatrix}0&1\\1&0\end{pmatrix}
 =\begin{pmatrix}1-U^{-2}&U^{-1}\\U^{-1}&0\end{pmatrix}
 =B(U^2).
\]
Its determinant is $-U^{-2}$, a Laurent unit. The triangular map can
therefore be inverted successively in increasing root height, and the
inverse has Laurent polynomial coefficients and preserves root
weight.

Each small-root residue matrix is unipotent: its nonidentity part
strictly lowers the residue height of the first strand while
preserving the total residue height. Its inverse is thus a finite
polynomial in that part. The residual Cartan matrix and the flip are
invertible monomial matrices, so $R^{\rm res}$ is invertible over $S$.
Composition and inversion preserve the triangular, root-homogeneous
form, and on each root block they give the corresponding product or
inverse of Burau matrices.
\end{proof}

\subsection{Expansion to all quantum orders}
\label{subsec:braid-all-order-resummation}

The associated graded of a crossing was computed in
Proposition~\ref{prop:braid-triangular-crossing}. To reach all quantum
orders we expand the exponential $e^{Z(\epsilon)}$ around the principal
flow $e^{A}$. The following lemma controls this expansion.

\begin{lemma}
\label{lem:braid-filtered-insertion}
Let $P$ be a polynomial algebra over a commutative $\mathbb Q$-algebra
$S$, and let $V$ be a finite free $S$-module with a finite split
integer filtration of width $W$.
Suppose
\[
 Z(\epsilon)=AI+N+\sum_{j\ge1}\epsilon^j Z_j
\]
is a matrix of polynomial differential operators on $P\otimes V$,
$A$ is a scalar locally nilpotent derivation, and $e^{sA}$ and
$e^{-sA}$ are polynomial jointly in $s$ and the coordinates.
If $N$ lowers the filtration by at least $\delta>0$, then
\[
 [\epsilon^h]e^{Z(\epsilon)}=\mathcal D_h e^A
\]
for a finite matrix $\mathcal D_h$ of polynomial differential operators,
and with $L=\lfloor W/\delta\rfloor$ the Duhamel expansion in order $h$
involves at most $h+(h+1)L$ insertions.
\end{lemma}

\begin{proof}
Write
\[
 B(\epsilon)=N+\sum_{j\ge1}\epsilon^jZ_j,
 \qquad B_s(\epsilon)=e^{sA}B(\epsilon)e^{-sA}.
\]
The function $U(s)=e^{s(A+B)}e^{-sA}$ satisfies $U'(s)=U(s)B_s$ and
$U(0)=I$, and iterating this equation gives the ordered-insertion
(Duhamel) expansion
\[
 e^{Z(\epsilon)}=
 \left(I+\sum_{\ell\ge1}
  \int_{0\le s_1\le\cdots\le s_\ell\le1}
       B_{s_1}(\epsilon)\cdots B_{s_\ell}(\epsilon)
       \,ds_1\cdots ds_\ell\right)e^A.
\]
We show that only finitely many terms contribute to each coefficient
of $\epsilon$. First, each $B_s$ is a polynomial differential operator
whose coefficients are polynomial in $s$: if $e^{sA}=C_{\phi_s}$, then
\[
 C_{\phi_s}M_fC_{\phi_s^{-1}}=M_{f\circ\phi_s},\qquad
 C_{\phi_s}\partial_iC_{\phi_s^{-1}}
 =\sum_j(\partial_i(\phi_s^{-1})_j)(\phi_s(\mathbf X))\partial_j,
\]
where $M_f$ denotes multiplication by $f$.

The scalar flow fixes the residue factor, so every conjugate of $N$
still lowers the filtration by $\delta$, and any product of $L+1$ such
conjugates is zero. A term of order $h$ contains at most $h$
insertions of positive order, and between two consecutive ones, as
well as before the first and after the last, at most $L$ insertions of
$N$. Its total length is therefore at most
\[
 h+(h+1)L.
\]
Only finitely many words contribute. Their coefficients are
polynomial on the ordered integration simplices, and their integrals
lie in $S$ because $S$ is a $\Q$-algebra. The resulting finite sum is
a polynomial differential operator $\mathcal D_h$ with
$[\epsilon^h]e^{Z(\epsilon)}=\mathcal D_he^A$.
\end{proof}

\begin{theorem}
\label{thm:braid-all-rank-insertion}
Fix $n\ge2$, $m>1$, and a signed braid $\beta$ acting on tensor powers
of $M(\mathbf u)$. With the change of basis
\eqref{eq:braid-global-gauge}, for every quantum order $h$
\[
 [\epsilon^h](\mathcal G\mathcal B_\beta\mathcal G^{-1})
       =\mathcal D_{\beta,h}C_{\Phi_\beta}.
\]
Here $\mathcal D_{\beta,h}$ is a finite matrix of polynomial
differential operators over $\mathbb Q(\xi)[\mathbf u^{\pm1}]$, and
$\Phi_\beta$ is the composite, along the braid word, of the triangular
polynomial substitutions of
Proposition~\ref{prop:braid-triangular-crossing}; its linear block at the
root $\alpha$ is the Burau matrix of the braid at $T_\alpha=U_\alpha^2$.
In order zero the residue matrix of the associated graded is constant
and invertible, and all other terms strictly lower the total residue
height.
\end{theorem}

\begin{proof}
We treat a single root factor, then a positive crossing, then inverse
crossings and composition.

For a positive root factor, put
\[
 P_{r,j}(a)=[\epsilon^j]A_{m,r,a}(\zeta(1+\epsilon)).
\]
Lemma~\ref{lem:braid-scalar-grouping} gives $\deg P_{r,j}\le2j$.
Expanding this polynomial in falling factorials,
$P_{r,j}(a)=\sum_\ell p_{r,j,\ell}(a)_\ell$, gives
\[
 \sum_{a\ge0}P_{r,j}(a)\frac{Z^a}{a!}
 =T_{r,j}(Z)e^Z,
 \qquad T_{r,j}(Z)=\sum_\ell p_{r,j,\ell}Z^\ell,
\]
an identity in the single operator $Z=Z(\epsilon)$. Hence,
coefficientwise on each finite height block $M_\kappa$,
\[
 \mathcal R_\alpha(\epsilon)
 =\sum_{r=0}^{m-1}\sum_{j\ge0}
      \epsilon^jV_r(\epsilon)T_{r,j}(Z(\epsilon))
                   e^{Z(\epsilon)},
\]
and in order $h$ only $j\le h$ occurs. The Taylor coefficients of
$V_r$ and $Z$ are polynomial differential operators by
Theorem~\ref{thm:braid-fixed-transfer-jets}. Moreover,
\[
 Z(0)=A_\alpha I+N_\alpha,
\]
where $A_\alpha$ has a polynomial flow and $N_\alpha$ lowers the
residue height by at least $m$ (Proposition~\ref{prop:braid-all-root-fields}
and Lemma~\ref{lem:braid-matrix-field}). Lemma~\ref{lem:braid-filtered-insertion}
and Corollary~\ref{cor:braid-uniform-residue-bound} therefore give
\[
 [\epsilon^h]\mathcal R_\alpha(\epsilon)
       =\mathcal D_{\alpha,h}e^{A_\alpha}.
\]

This form is preserved under composition, since
\[
 (D_1C_{\Phi_1})(D_2C_{\Phi_2})
 =D_1\bigl(C_{\Phi_1}D_2C_{\Phi_1}^{-1}\bigr)
           C_{\Phi_1}C_{\Phi_2},
\]
and conjugation by a polynomial map with polynomial inverse preserves
the polynomial Weyl algebra. The product over the positive roots
therefore has the required form. The normalized Cartan factor adds the
common Laurent dilation of
Proposition~\ref{prop:braid-triangular-crossing}: its mixed phase is
cancelled by the prefix gauge, and its correction of order $h$ is a
polynomial of degree at most $2h$ in the coarse indices, hence an
Euler differential operator. This proves the assertion for a positive
checked crossing.

For inversion, write
$\mathcal R(\epsilon)=\mathcal D(\epsilon)C_\Phi$ and
$\mathcal D_0=R_0+N_0$, where $R_0$ is the invertible constant residue
matrix and $N_0$ lowers the total residue height by at least $m$.
Since $R_0^{-1}$ preserves this height,
Corollary~\ref{cor:braid-uniform-residue-bound} gives
\[
 \mathcal D_0^{-1}
 =\left(\sum_{j=0}^{L_{n,m}}(-R_0^{-1}N_0)^j\right)R_0^{-1}.
\]
The inverse series $\mathcal E(\epsilon)=\mathcal D(\epsilon)^{-1}$
is then determined in the polynomial Weyl algebra by
\[
 \mathcal E_0=\mathcal D_0^{-1},\qquad
 \mathcal E_h=-\mathcal D_0^{-1}
          \sum_{j=1}^h\mathcal D_j\mathcal E_{h-j}.
\]
Consequently
\[
 \mathcal R(\epsilon)^{-1}
 =C_{\Phi^{-1}}\mathcal E(\epsilon)
 =\bigl(C_{\Phi^{-1}}\mathcal E(\epsilon)C_\Phi\bigr)
                  C_{\Phi^{-1}}
\]
has the same form over $S$.

The gauge is common to all strands, so composition along a signed
braid word proves the asserted formula. In order zero the
height-preserving residue matrices multiply, and every other term
lowers the same total filtration.

Finally, the operators obtained in this way do not depend on the
height block. In the Duhamel expansion of each root factor, the number
of insertions in order $h$ is at most $h+(h+1)L_{n,m}$ for every
$M_\kappa$, and the remaining factors use only $j\le h$ and
$\deg T_{r,j}\le2j$. The differential operators above are therefore
single elements of the polynomial Weyl algebra acting on every block.
Since the PBW module is the direct sum of the $M_\kappa$, the
identities proved on the finite blocks hold on the whole module.
\end{proof}

\section{Graded traces and knot invariants}
\label{sec:trace}

The braid substitutions of Section~\ref{sec:braid} preserve root
weight but not ordinary polynomial degree: their nonlinear terms raise
the degree, whereas a fixed differential insertion can lower it only
by a bounded amount. Only finitely many nonlinear terms can therefore
contribute to a diagonal matrix coefficient. We use this to reduce the
braid trace to traces of linear substitutions, whose determinants give
the Alexander denominators. We then remove the regulator, compare the
trace with the colored invariants, and treat the root $Q=1$ by the
rational Kontsevich invariant.

\subsection{Traces of polynomial substitutions}

Let $S$ be a commutative $\mathbb Q$-algebra and $I$ a finite set of
blocks, with a positive integer $w_\alpha$ attached to each
$\alpha\in I$, and write
\[
 \mathbf X_\alpha=(X_{\alpha,1},\ldots,X_{\alpha,d_\alpha}),
 \qquad P=S[\mathbf X].
\]
The variables in block $\alpha$ have regulating weight $w_\alpha$, so
the regulating degree of a monomial $\mathbf X^a$ is
$\operatorname{wt}(a)=\sum_\alpha w_\alpha|a_\alpha|$, while its
ordinary degree is $|a|=\sum_\alpha|a_\alpha|$.
Let $P_e$ be the homogeneous part of regulating degree $e$ and
$\pi_e:P\to P_e$ the projection; since all weights are positive,
$P_e$ is finite free over $S$. For an $S$-linear operator $T:P\to P$
define
\[
 \operatorname{Tr}_s(T)=
 \sum_{e\ge0}s^e\operatorname{Tr}(\pi_eT|_{P_e}).
\]
This is a formal power series in $s$ whose coefficients are traces on
finite free modules.

Consider a polynomial map of the form
\begin{equation}\label{eq:trace-weighted-map}
 \Phi_\alpha(\mathbf X)=B_\alpha\mathbf X_\alpha+
                \mathbf P_\alpha(\mathbf X),
 \qquad L_\alpha=B_\alpha\mathbf X_\alpha,
\end{equation}
where the coordinates of $\mathbf P_\alpha$ are homogeneous of
regulating weight $w_\alpha$ and of ordinary degree at least two; the
matrices $B_\alpha$ may be singular. This covers the root-homogeneous
maps of Section~\ref{sec:braid}: there each variable
$X_{\alpha,j}$ has root weight $\alpha$, and a nonlinear monomial of
weight $\alpha$ only involves roots of smaller height, so
$w_\alpha=\hgt(\alpha)$ is a grading of the required kind.
Multiplying all weights by $m$ gives the grading on the coarse PBW
variables. We write $C_\Phi f=f\circ\Phi$, $M_f$ for multiplication by
$f$, and $[f]_r$ for the component of ordinary degree $r$ of a
polynomial $f$.

\begin{lemma}
\label{lem:trace-finite-replacement}
Let $\mathcal D=\mathbf X^a\partial^b$ be a normal-ordered differential
monomial and put $N=|b|-|a|$. The trace of $\mathcal D C_\Phi$ vanishes
if $N<0$ or $\operatorname{wt}(a)\ne\operatorname{wt}(b)$; otherwise
\begin{equation}\label{eq:trace-finite-replacement}
 \operatorname{Tr}_s(\mathcal D C_\Phi)=
 \sum_{|\nu|\le N}\frac{s^{\operatorname{wt}(\nu)}}{\nu!}
 \operatorname{Tr}_s\left(
   \partial^\nu\mathcal D\,
   M_{[\mathbf P^\nu]_{N+|\nu|}} C_L\right).
\end{equation}
Here $\nu$ runs over multi-indices in the polynomial coordinates.
\end{lemma}

\begin{proof}
Taylor expansion of $f(L(\mathbf X)+\mathbf P(\mathbf X))$ around
$L(\mathbf X)$ gives the locally finite operator identity
\[
 C_\Phi=\sum_{\nu\ge0}\frac1{\nu!}
           M_{\mathbf P^\nu}C_L\partial^\nu.
\]
Consider the term of $\mathcal DC_\Phi$ coming from the component of
ordinary degree $r$ of $\mathbf P^\nu$. It changes ordinary degree by
\[
 |a|-|b|+r-|\nu|=r-|\nu|-N,
\]
so it contributes to a diagonal coefficient only if $r=N+|\nu|$. On
the other hand $r\ge2|\nu|$ whenever this component is nonzero, because
every coordinate of $\mathbf P$ has degree at least two. Together these
give
\[
 |\nu|\le N,
\]
which proves finiteness and the vanishing for $N<0$. Since $C_\Phi$
preserves regulating degree, a nonzero trace also requires
$\operatorname{wt}(a)=\operatorname{wt}(b)$.

For a remaining index $\nu$, set $r_\nu=\operatorname{wt}(\nu)$ and
\[
 A_\nu=\mathcal D M_{[\mathbf P^\nu]_{N+|\nu|}}C_L.
\]
The operators $A_\nu$ and $\partial^\nu$ shift regulating degree by
$r_\nu$ and $-r_\nu$, respectively, so they give maps between the
finite free modules $P_e$ and $P_{e-r_\nu}$. The identity
$\operatorname{Tr}(AB)=\operatorname{Tr}(BA)$ for rectangular matrices
gives
\[
 \Tr_{P_e}(A_\nu\partial^\nu)
 =\Tr_{P_{e-r_\nu}}(\partial^\nu A_\nu),
\]
where $P_j=0$ for $j<0$. Multiplying by $s^e$ and changing the index
to $e'=e-r_\nu$ produces the factor $s^{r_\nu}$ in
\eqref{eq:trace-finite-replacement}.
\end{proof}

Put $\delta_\alpha(s)=\det(I-s^{w_\alpha}B_\alpha)$. Its constant term
is one, so it is invertible in $S[[s]]$.

\begin{lemma}
\label{lem:trace-linear-kernel}
For a block of weight $w$, write $C_Bf(X)=f(BX)$.
For auxiliary vectors $a,b$ of commuting indeterminates,
\begin{equation}\label{eq:trace-linear-kernel}
 \operatorname{Tr}_s
   \left(e^{a\cdot X}e^{b\cdot\partial}C_B\right)
 =\frac{\exp\left(a^T s^wB(I-s^wB)^{-1}b\right)}
        {\det(I-s^wB)}.
\end{equation}
In particular, the trace of $X^c\partial^d C_B$ vanishes unless
$|c|=|d|$, and its denominator divides $\det(I-s^wB)^{1+|d|}$.
\end{lemma}

\begin{proof}
In one variable, expand both exponentials and use
\[
 \sum_{n\ge0}(n)_r z^n=\frac{r!z^r}{(1-z)^{r+1}}.
\]
The formula for a diagonal matrix follows by multiplying the
one-variable formulas. A generic matrix is diagonalizable over an
algebraic extension of the field of rational functions in its entries.
Under the change of variables $X=PZ$ the data become
\[
 B'=P^{-1}BP,\qquad a'=P^Ta,\qquad b'=P^{-1}b,
\]
and both sides of \eqref{eq:trace-linear-kernel} are unchanged, so the
formula holds for the generic matrix. For each coefficient in $a,b$,
clearing the determinant gives a polynomial identity in the matrix
entries. It therefore holds over $\mathbb Q$ and, by base change, over
every $S$.

The coefficient of $a^cb^d$ involves at most $|d|$ entries of
$(I-s^wB)^{-1}$, each with denominator $\det(I-s^wB)$, and the
determinant outside the exponential contributes one more factor; this
gives the power $1+|d|$. The trace vanishes when $|c|\ne|d|$ because
$C_B$ preserves ordinary degree.
\end{proof}

\begin{theorem}
\label{thm:trace-general-trace}
For every finite polynomial differential operator $\mathcal D$,
\begin{equation}\label{eq:trace-general-trace-ring}
 \operatorname{Tr}_s(\mathcal D C_\Phi)
 \in S[s,\delta_\alpha(s)^{-1}:\alpha\in I].
\end{equation}
For a normal monomial $\mathcal D=X^a\partial^b$ with $N=|b|-|a|\ge0$,
a common denominator is
\begin{equation}\label{eq:trace-coarse-pole-bound}
 \prod_{\alpha\in I}
       \delta_\alpha(s)^{1+|b_\alpha|+N}.
\end{equation}
Without insertions,
\begin{equation}\label{eq:trace-uninserted-trace}
 \operatorname{Tr}_s(C_\Phi)
       =\prod_{\alpha\in I}\delta_\alpha(s)^{-1}.
\end{equation}
\end{theorem}

\begin{proof}
By linearity it suffices to treat $\mathcal D=X^a\partial^b$. In
Lemma~\ref{lem:trace-finite-replacement} every remaining index satisfies
$|\nu|\le N$. Normal ordering
\[
 \partial^\nu\mathcal D
      M_{[\mathbf P^\nu]_{N+|\nu|}}
\]
can only remove derivatives, so the number of derivatives in block
$\alpha$ is at most $|b_\alpha|+|\nu_\alpha|$. Since $L$ is block
diagonal, the trace of a normal monomial composed with $C_L$ is the
product of the traces over the blocks, and
Lemma~\ref{lem:trace-linear-kernel} gives
\[
 \text{denominator}\ \mid\
 \prod_\alpha\delta_\alpha(s)^{1+|b_\alpha|+|\nu_\alpha|}
 \ \mid\
 \prod_\alpha\delta_\alpha(s)^{1+|b_\alpha|+N}.
\]
The replacement formula involves only nonnegative powers of $s$,
which proves \eqref{eq:trace-general-trace-ring} and
\eqref{eq:trace-coarse-pole-bound}. For $\mathcal D=1$ one has $N=0$, only
$\nu=0$ remains, and Lemma~\ref{lem:trace-linear-kernel} with $a=b=0$
gives \eqref{eq:trace-uninserted-trace}.
\end{proof}

The theorem will also be used for operators on $P\otimes_S E$, where
$E$ is a finite free module with a homogeneous basis
$e_1,\ldots,e_v$ of nonnegative degrees $d_1,\ldots,d_v$. The graded
trace of a matrix operator $T=(T_{ij})$ is then
$\sum_i s^{d_i}\operatorname{Tr}_s(T_{ii})$, and the same determinant
denominators apply entry by entry.

\subsection{Partial traces}

Cutting a strand produces operators with derivatives on the right of
the substitution. The following corollary moves them to the left.

\begin{corollary}\label{cor:trace-two-sided-trace}
The conclusion \eqref{eq:trace-general-trace-ring} also holds for finite
sums of operators $\mathcal D_\eta C_\Phi\partial^\eta$; more precisely,
\begin{equation}\label{eq:trace-right-cyclicity}
 \operatorname{Tr}_s(\mathcal D_\eta C_\Phi\partial^\eta)
 =s^{\operatorname{wt}(\eta)}
       \operatorname{Tr}_s(\partial^\eta\mathcal D_\eta C_\Phi).
\end{equation}
For a cut normal monomial $M_{X^a}C_\Phi\partial^b$ with
$N=|b|-|a|$, the trace vanishes if $N<0$, and otherwise
\eqref{eq:trace-coarse-pole-bound} is a common denominator.
\end{corollary}

\begin{proof}
Decompose $\mathcal D_\eta$ according to its finitely many shifts of
regulating degree; only the shift $\operatorname{wt}(\eta)$
contributes to the trace. On each finite free homogeneous piece the
identity $\operatorname{Tr}(AB)=\operatorname{Tr}(BA)$, used as in
Lemma~\ref{lem:trace-finite-replacement}, applies, and reindexing by the
degree of the right derivative gives \eqref{eq:trace-right-cyclicity}.
Since $\partial^\eta\mathcal D_\eta$ is a finite polynomial differential
operator, Theorem~\ref{thm:trace-general-trace} applies to it. For the cut
monomial, normal ordering $\partial^bM_{X^a}$ lowers the multiplication
and derivative exponents together in each coordinate, so $N$ is
unchanged and the number of derivatives in block $\alpha$ is at most
$|b_\alpha|$; this gives the bound.
\end{proof}

\begin{lemma}
\label{lem:trace-general-cut}
Suppose that the coordinates of each block are distributed over
several strands. Let $\iota$ be the inclusion of the polynomials that
do not depend on the coordinates of the cut strand, and let
$\mathrm{ev}_0$ set those coordinates to zero. For a polynomial
differential operator $\mathcal D$ of order at most $d$,
\[
 \mathrm{ev}_0\mathcal D C_\Phi\iota
 =\sum_\eta M_{p_\eta}C_\Psi\partial^\eta,
 \qquad
 \Psi=\Phi_{\mathrm{noncut}}|_{\mathrm{cut}=0},
\]
with $|\eta|\le d$ and polynomials $p_\eta$. The map $\Psi$ satisfies
the same weighted hypotheses, and its linear block at each root is the
corresponding principal submatrix of $B_\alpha$.
\end{lemma}

\begin{proof}
Apply the chain rule to $\mathcal D(f(\Phi_{\mathrm{noncut}}))$ and set
the cut coordinates to zero afterwards. The coefficients are
polynomials, and the remaining derivatives act on $f$ at $\Psi$.
Setting variables to zero preserves regulating homogeneity, and a
monomial of ordinary degree at least two either keeps its degree or
vanishes; so no linear terms are created, and the linear part of $\Psi$
is the principal submatrix. Each derivative of $f$ comes from one
application of the chain rule, and there are at most $d$ of them,
which gives $|\eta|\le d$.
\end{proof}

The following example shows how the replacement formula selects the
nonlinear terms in a case with three heights.

\begin{example}
\label{ex:trace-height-three-trace}
Give $(x,z,w,u,v,y)$ the weights $(1,1,1,2,2,3)$ and consider the map
\[
 \begin{split}
 &(x,z,w)\longmapsto(ax,bz,cw),\\
 &u\longmapsto du+p xz,\qquad v\longmapsto ev+q zw,\\
 &y\longmapsto fy+r xv+\ell uw+t xzw,
 \end{split}
\]
with coefficients in $S$; the six variables can be thought of as the
coordinates of the six positive roots of $A_3$. For
$\mathcal D=y\partial_x\partial_z\partial_w$ we have $N=2$. The
nonlinear terms contain no $y$, so equality of the multiplication and
differentiation degrees in the $y$ block forces $\nu_y=1$. Since
$|\nu|\le2$ and only the coordinates $u$, $v$ and $y$ have nonlinear
terms, the possible indices are $e_y$, $e_u+e_y$ and $e_v+e_y$.
Balancing the remaining blocks selects the polynomial terms $t xzw$,
$p\ell xzuw$ and $qr xzwv$, respectively. The resulting traces are
evaluated with the scalar identity
$\operatorname{Tr}_s(\partial_X M_X C_a)=(1-as^w)^{-2}$. Writing
$\delta_a=1-as$, $\delta_d=1-ds^2$ and $\delta_f=1-fs^3$, with the
analogous notation for $b,c,e$, one obtains
\[
 \operatorname{Tr}_s(\mathcal D C_\Phi)
 =\frac{s^3}{\delta_a^2\delta_b^2\delta_c^2
                  \delta_d\delta_e\delta_f^2}
       \left(t+\frac{p\ell s^2}{\delta_d}
                    +\frac{qr s^2}{\delta_e}\right).
\]
The two nested contributions have poles only along divisors that are
already present in the linear part.
\end{example}

\subsection{Alexander polynomials and normalization}

We use the blocks
\[
 B(T)=\begin{pmatrix}
       1-T^{-1}&T^{-1/2}\\ T^{-1/2}&0
      \end{pmatrix},\qquad
 B(T)^{-1}=\begin{pmatrix}
       0&T^{1/2}\\T^{1/2}&1-T
      \end{pmatrix}.
\]
Placing them in adjacent positions defines $B_\beta(T)$ for a signed
braid $\beta$ on $b$ strands. The superscript $\mathrm{tr}$ denotes
deletion of the first row and column, which corresponds to cutting
the first strand. If the closure of $\beta$ is a knot $K$ and
$w=w(\beta)$ is the writhe, then
\begin{equation}\label{eq:trace-exact-burau-alexander}
 \det(I-B_\beta^{\mathrm{tr}}(T))
    =T^{(b-1-w)/2}\Delta_K(T).
\end{equation}
For these blocks this is the Burau--Alexander identity of
\cite[Section~2.2, equations~(53) and~(61)]{GW}. The exponent of $T$ is
an integer: the permutation of the braid is a $b$-cycle, of sign
$(-1)^{b-1}$, and it is a product of as many transpositions as there
are crossings, so its sign is also $(-1)^w$.

It can also be derived from the classical formula for the reduced
Burau representation, as follows.
With $\Gamma=\operatorname{diag}(1,T^{1/2},\ldots,T^{(b-1)/2})$,
conjugation $B\mapsto\Gamma^{-1}B\Gamma$ takes the standard unreduced
blocks
$\left(\begin{smallmatrix}1-T^{-1}&T^{-1}\\1&0\end{smallmatrix}\right)$
to the symmetric blocks, and these fix
$v=(1,T^{-1/2},\ldots,T^{-(b-1)/2})^T$ from both sides. At $T=1$ the
matrix $B_\beta$ is the permutation matrix of a single cycle, so
$I-B_\beta$ has rank $b-1$ over the field of rational functions, and
its adjugate is $c\,vv^T$ for some $c$. The first principal cofactor
is then $c$, the determinant of $I$ minus the reduced Burau matrix is
the trace of the adjugate, $c\sum_{j=0}^{b-1}T^{-j}$, and the usual
normalized formula for the reduced Burau representation gives
\eqref{eq:trace-exact-burau-alexander}.

\subsection{Dominant colors}

Let $\mathfrak g=\mathfrak{sl}_n$ and let $\lambda$ be dominant. We
work over $\mathbb Q(q)$, adjoining a root of $q$ where the
unnormalized Cartan scalar requires it. Let $M_\lambda$ be the Verma
module, $W_\lambda$ its maximal proper weight submodule,
$L_\lambda=M_\lambda/W_\lambda$, and
$\rho=\frac12\sum_{\alpha>0}\alpha$ the Weyl vector. A vector of root
deficit $\gamma$ is weighted by the regulator
$s^{\operatorname{ht}(\gamma)}$. Explicitly, with $w_{\mathbf N}$ the
ordinary PBW vector and
$|\mathbf N|_{\rm ht}=\sum_{\alpha>0}\operatorname{ht}(\alpha)N_\alpha$,
the normalized regulated cut scalar is
\begin{equation}\label{eq:trace-regulated-source}
 \begin{split}
 Z_\beta(\mathbf u,q;s)
 ={}&\mathbf u^{(b-1-w)2\rho}
 \sum_{\mathbf N_2,\ldots,\mathbf N_b}
       (s q^{-2})^{\sum_{j=2}^b|\mathbf N_j|_{\rm ht}}\\
 &\quad{}\cdot
 [w_{\mathbf0}\otimes w_{\mathbf N_2}\otimes\cdots\otimes w_{\mathbf N_b}]
 \mathcal B_\beta
 (w_{\mathbf0}\otimes w_{\mathbf N_2}\otimes\cdots\otimes w_{\mathbf N_b}).
 \end{split}
\end{equation}
Here the brackets extract the indicated diagonal coefficient, and the
coefficient of each power of $s$ is a finite sum. The following
comparison, at generic $q$, is a regulated form of the Verma-module
cut-trace argument in \cite[Lemma~7.7]{GaroufalidisLeQHolonomic}.

\begin{theorem}
\label{thm:trace-quotient-selection}
For a signed braid on $b$ strands whose permutation is a single
cycle, the regulated highest-vector partial trace of the normalized
braid operator on $M_\lambda^{\otimes b}$ equals, coefficientwise in
the regulator, the corresponding trace on $L_\lambda^{\otimes b}$. The
latter is a polynomial in $s$, and its value at $s=1$, together with
the pivot and the zero-framing correction, is the reduced colored
invariant.
\end{theorem}

\begin{proof}
Choose a complement $V_\lambda$ to $W_\lambda$, spanned by weight
vectors and containing the highest vector. The quotient map identifies
$V_\lambda$ with $L_\lambda$ as a graded vector space. For a subset
$J$ of tensor positions, let $W^{(J)}$ be the span of the tensors with
a factor in $W_\lambda$ at each position in $J$ and arbitrary factors
elsewhere. Since the braiding and its inverse are module maps,
\[
 \mathcal B_\beta(W^{(J)})\subseteq W^{(\pi(J))},
\]
where $\pi$ is the permutation of the braid. Use a tensor basis adapted
to $W_\lambda\oplus V_\lambda$. A basis vector whose factors in
$W_\lambda$ sit exactly at the positions in $J$ is sent to a
combination of basis vectors with a factor in $W_\lambda$ at every
position of $\pi(J)$. Its diagonal coefficient can therefore be
nonzero only if $\pi(J)\subseteq J$, and since the two sets have the
same cardinality, only if $\pi(J)=J$.

The permutation $\pi$ is a single cycle, so its only invariant subsets
are the empty set and the full set. The cut factor is the highest
vector, which lies in $V_\lambda$, and this excludes the full set.
Only the block with every tensor factor in $V_\lambda$ contributes,
and modulo $W_\lambda$ this block is the braid operator on
$L_\lambda^{\otimes b}$.

Both the pivot and the regulator are scalar on each weight space.
Each regulating degree contains only finitely many weights and PBW
vectors, so the adapted change of basis preserves every coefficient
of the partial trace. The quotient is finite-dimensional, so its trace
is a polynomial in $s$, and at $s=1$ it is the reduced ribbon partial
trace with the stated normalization.
\end{proof}

Put $u_i=q^{r_i}$ for
$\lambda=\sum_{i=1}^{n-1}r_i\omega_i$, and let $w$ be the writhe.
The product of the Cartan scalar of the highest vector, the ribbon
factor and the pivotal factor of the highest vector is
\begin{equation}\label{eq:trace-general-prefactor}
 \mathbf u^{(b-1-w)2\rho}
 =\prod_{i=1}^{n-1}u_i^{(b-1-w)i(n-i)}.
\end{equation}
The three factors contribute the exponents $w(\lambda,\lambda)$,
$-w(\lambda,\lambda+2\rho)$ and $(b-1)(2\rho,\lambda)$ of $q$. In their
sum $(\lambda,\lambda)$ cancels, which leaves
$(b-1-w)(2\rho,\lambda)$, and $2\rho=\sum_ii(n-i)\alpha_i$. On a
vector of root deficit $\gamma$ the remaining pivot is
$q^{-(2\rho,\gamma)}=q^{-2\operatorname{ht}(\gamma)}$.

\begin{lemma}
\label{lem:trace-square-color-descent}
The normalized regulated cut scalar of a knot is invariant under the
sign changes $u_i\mapsto\sigma_i u_i$, $\sigma_i\in\{\pm1\}$, chosen
independently. Consequently, a rational germ at a root of unity whose
denominator is a polynomial in the squares of the root characters
descends to the variables $t_i=u_i^2$.
\end{lemma}

\begin{proof}
Identify the two Verma modules through their ordinary negative PBW
bases. Then a positive root operator of root weight $\nu$ is
multiplied by the character $\sigma^\nu$, negative root operators are
unchanged, and the highest-normalized checked crossing with incoming
deficits $\gamma,\delta$ and root transfer $\nu$, whose outgoing
deficits are $(\delta+\nu,\gamma-\nu)$, acquires the factor
$\sigma^{\gamma+\delta+\nu}$. The diagonal gauge transformation
\[
 v_{\gamma_1}\otimes\cdots\otimes v_{\gamma_b}
 \longmapsto
 \sigma^{\sum_{j=1}^b j\gamma_j}
 v_{\gamma_1}\otimes\cdots\otimes v_{\gamma_b}
\]
contributes at an adjacent crossing the ratio
$\sigma^{\gamma-\delta-\nu}
 =\sigma^{\gamma+\delta+\nu}$, so it conjugates each crossing and each
inverse crossing to its sign-changed version. It commutes with the
pivot and the regulator, and it cancels in a diagonal entry of the cut.

It remains to consider the prefactor \eqref{eq:trace-general-prefactor}.
The permutation of the braid has sign $(-1)^w$ and a $b$-cycle has
sign $(-1)^{b-1}$, so $b-1-w$ is even, all exponents in
\eqref{eq:trace-general-prefactor} are even, and the prefactor is
invariant under the sign changes. Since the denominator is invariant
as well, taking invariants in the localized Laurent ring over a field
of characteristic zero gives the descent to the variables
$t_i=u_i^2$.
\end{proof}

\subsection{Rationality at nontrivial roots}

\begin{proposition}\label{prop:trace-all-rank-trace-closure}
Let $\beta$ be a signed braid whose closure is a knot $K$, and let
$\zeta$ be a primitive $m$th root with $m>1$.
Each coefficient of the expansion of the normalized cut scalar at
$Q=\zeta$ has a rational value at $s=1$, and these values form a
series
\[
 g_{K,n,\zeta}(\mathbf t,x)\in
 \mathbb Q(\zeta)
 [t_1^{\pm1},\ldots,t_{n-1}^{\pm1},
      D_{K,n}(t_1^m,\ldots,t_{n-1}^m)^{-1}][[x]],
 \qquad x=Q-\zeta.
\]
Its germ at $\mathbf t=\mathbf1$ is the Newton knot germ, its
specializations are the dominant colored invariants, and its constant
coefficient has denominator dividing $D_{K,n}(\mathbf t^m)$.
\end{proposition}

\begin{proof}
Choose $\xi^2=\zeta$ with $\xi^m=(-1)^{m+1}$ and work over
$S=\Q(\xi)[\mathbf u^{\pm1}]$. Theorem~\ref{thm:braid-all-rank-insertion}
expresses each quantum coefficient as a finite differential insertion
followed by a common root-homogeneous substitution, whose block at
the root $\alpha$ is $B_\beta(T_\alpha)$ with
\[
 T_\alpha=(u^\alpha)^{2m},\qquad
 w_\alpha=m\hgt(\alpha).
\]
The coarse part of the pivot is one in order zero; its higher
coefficients are polynomials in weighted Euler operators and can be
included in the insertion. The prefactor
\eqref{eq:trace-general-prefactor} is a Laurent monomial.

After the cut, Lemma~\ref{lem:trace-general-cut},
Corollary~\ref{cor:trace-two-sided-trace} and
Theorem~\ref{thm:trace-general-trace} allow only the denominators
\[
 \delta_{\alpha,s}
 =\det\!\left(I-s^{m\hgt(\alpha)}
                  B_\beta^{\mathrm{tr}}(T_\alpha)\right).
\]
By the Burau--Alexander identity \eqref{eq:trace-exact-burau-alexander},
\[
 \delta_{\alpha,1}
 =T_\alpha^{(b-1-w)/2}\Delta_K(T_\alpha),
\]
and the monomial factor is a unit since $b-1-w$ is even. Each quantum
coefficient therefore has a rational value at $s=1$ in
\[
 \Q(\xi)[\mathbf u^{\pm1},D_{K,n}(\mathbf u^{2m})^{-1}].
\]

To compare with a fixed dominant color $\mathbf r$, expand at
$u_i=\xi^{r_i}$ and put $y_i=u_i/\xi^{r_i}-1$. At this center
$U_\alpha=\pm1$, so each elementary root block is a signed
transposition. The braid block is therefore a signed permutation
matrix whose permutation is a single cycle. Deleting the cut vertex
breaks this cycle into a path, so the remaining matrix is nilpotent.
Consequently
\[
 \delta_{\alpha,s}(\mathbf0)=1,\qquad
 \delta_{\alpha,s}=1+E_\alpha(s,\mathbf y),\quad
 E_\alpha\in(\mathbf y)\Q(\xi)[s][[\mathbf y]].
\]
For every fixed $a\ge1$ and $L\ge0$,
\[
 \delta_{\alpha,s}^{-a}
 \equiv\sum_{j=0}^{L}(-1)^j\binom{a+j-1}{j}E_\alpha^j
                     \pmod{(\mathbf y)^{L+1}}.
\]
The numerator of each quantum coefficient is polynomial in $s$, and
this finite sum replaces each denominator to the required color
precision. For each quantum order $h$ and each bound $L$ on the color
degree, the coefficient of order $h$ modulo $(\mathbf y)^{L+1}$ is
therefore a polynomial in $s$.

At generic $q$, Theorem~\ref{thm:trace-quotient-selection} identifies
the regulated cut scalar with the finite-dimensional colored trace.
Every coefficient of $s$ is a finite sum of regular PBW coefficients,
so this identity can be expanded at $q=\xi$. On the dominant-color
curve, where $q=\xi(1+\epsilon)^{1/2}$, we have
\[
 y_i=(1+\epsilon)^{r_i/2}-1\in\epsilon\Q(\xi)[[\epsilon]].
\]
The coefficient of $\epsilon^H$ after substitution only involves
quantum orders $h\le H$ and color degrees at most $H-h$. Each of these
finitely many jets is polynomial in $s$, so evaluation at $s=1$
commutes with substitution and coefficient extraction. It follows that
\[
 g_{K,n,\zeta}(Q^{r_1},\ldots,Q^{r_{n-1}},Q-\zeta)
       =J_{K,\sum r_i\omega_i}^{(n)}(Q).
\]

Lemma~\ref{lem:trace-square-color-descent} first descends the series to
the variables $t_i=u_i^2$. If $m$ is odd, $\Q(\xi)=\Q(\zeta)$. If $m$
is even, compare the series with its conjugate under $\xi\mapsto-\xi$
on the curves
\[
 t_i=(1+x/\zeta)^{2ma_i},\qquad a_i\in\N.
\]
These curves are fixed by the conjugation, since
$(1+x/\zeta)^{2ma_i}=Q^{2ma_i}$, and on them both series give the same
colored Laurent polynomial. The proof of
Lemma~\ref{lem:rel-color-rigidity} applies to this arithmetic
progression and identifies the two germs at $\mathbf t=\mathbf1$.
Since $D_{K,n}(\mathbf1)=1$, rational Taylor expansion is injective,
which gives the descent to $\Q(\zeta)$. The same comparison with the
Newton knot element identifies their germs at the identity color for
every $m$.

In order zero, the strictly residue-lowering operators have zero
diagonal trace. Only the constant residue matrix and the common
substitution remain, and \eqref{eq:trace-uninserted-trace} contributes
each $\delta_{\alpha,s}$ once. Thus
\[
 D_{K,n}(\mathbf t^m)g_{K,n,\zeta}(\mathbf t,0)
          \in\Q(\zeta)[\mathbf t^{\pm1}].
\]
Its Taylor germ is integral because it equals the corresponding
Newton germ, and Lemma~\ref{lem:rel-integral-saturation} puts
this Laurent numerator in $\Z[\zeta][\mathbf t^{\pm1}]$.
\end{proof}

\subsection{The identity root and the general type-\texorpdfstring{$A$}{A} lift}
\label{subsec:trace-identity-lift}

At $Q=1$ we use the rational Kontsevich invariant \cite{GK},
the Lie group weight systems of \cite{Beads}, and the comparison
with quantum invariants in \cite{LM}.

Normalize the invariant form by $(\alpha_i,\alpha_i)=2$, and put
$q=e^{h/2}$, $Q=e^h$. Let $\mathfrak h$ be the Cartan subalgebra,
$H_i$ its simple coroots, and $C=((\alpha_i,\alpha_j))$ the Cartan
matrix. In the next display $E_\alpha,F_\alpha$ are the classical
root vectors, normalized to be dual for this form. Put
\[
 \Omega_{\mathfrak h}=\sum_{i,j}(C^{-1})_{ij}H_i\otimes H_j,
 \qquad
 \Omega=\Omega_{\mathfrak h}+\sum_{\alpha>0}
       (E_\alpha\otimes F_\alpha+F_\alpha\otimes E_\alpha).
\]
The unnormalized universal $R$-matrix underlying our crossing
convention satisfies
\[
 R=1+h\left(\tfrac12\Omega_{\mathfrak h}
       +\sum_{\alpha>0}E_\alpha\otimes F_\alpha\right)+O(h^2),
 \qquad R_{21}R=1+h\Omega+O(h^2).
\]
We apply the comparison of the Kontsevich and Reshetikhin--Turaev
invariants in \cite[Section~10, Theorem~10]{LM} with the invariant
form normalized as above. The braiding on the associator side is then
$\exp(h\Omega/2)$. The parameters match: on the highest summand of
$V_\lambda\otimes V_\mu$ the two monodromies have eigenvalues
$q^{2(\lambda,\mu)}$ and $e^{h(\lambda,\mu)}$, so the comparison uses
$h=\log Q$. The Cartan scalar cancels in the reduced zero-framed
invariant by \eqref{eq:trace-general-prefactor}.

The rationality theorem of Garoufalidis and Kricker
\cite[Section~7.4, Theorem~14]{GK} and the Lie group weight system
\cite[Definitions~1.6--1.7 and Theorem~1]{Beads} give a reduced series
\begin{equation}\label{eq:trace-identity-rational-input}
 \mathscr R_{K,n}(\mathbf z,h)
 \in\mathbb Q[z_1^{\pm1},\ldots,z_{n-1}^{\pm1},
                       D_{K,n}(\mathbf z)^{-1}][[h]]
\end{equation}
with
\begin{equation}\label{eq:trace-identity-colored-input}
 \mathscr R_{K,n}(e^{h(r_1+1)},\ldots,e^{h(r_{n-1}+1)},h)
       =J_{K,\sum_i r_i\omega_i}^{(n)}(e^h).
\end{equation}
The shift by $\rho$ comes from \cite[Theorem~1(b)]{Beads}: the Lie
group weight system is evaluated at $\lambda+\rho$, and the unreduced
invariant carries the common quantum-dimension factor, which is
removed by dividing by the invariant of the unknot. The rational
invariant of the unknot is the empty matrix with the empty graph
\cite[Remark~7.4]{GK}, and its Lie group weight is one.

Beads act through the adjoint representation, whose torus characters
are roots, so every coefficient is a rational function of the
simple-root characters. At $\lambda+\rho$ these characters are
$z_i=e^{h(r_i+1)}=Qt_i$, which gives
\eqref{eq:trace-identity-colored-input}.

The poles can be read off from the weight system. For a fixed knot
each bead lies in $\Delta_K(T)^{-1}\mathbb Q[T^{\pm1}]$
\cite[Section~1.2]{Beads}. Under the weight system, its
denominator becomes $\Delta_K(\operatorname{Ad}(\mathbf z))$. This
operator acts on the Cartan subalgebra by $\Delta_K(1)=1$, and on the
root spaces of $\pm\alpha_{a,b}$ by $\Delta_K(z_{a,b}^{\pm1})$, where
$z_{a,b}=\prod_{i=a}^{b-1}z_i$. The matrix part, whose normalized
determinant is $\Delta_K$, contributes
\[
 \det\bigl(\Delta_K(\operatorname{Ad}(\mathbf z))\bigr)^{-1/2}
      =D_{K,n}(\mathbf z)^{-1}.
\]
Here Alexander symmetry gives equal factors for opposite roots, and
the branch of the square root is fixed by its value one at the
identity. Graphs of positive degree contribute positive powers of
$h$, so
\[
 [h^0]\mathscr R_{K,n}=D_{K,n}^{-1}.
\]
All tensors and contractions can be taken in a rational root basis,
so the coefficient field in \eqref{eq:trace-identity-rational-input} is
$\mathbb Q$.

\begin{proposition}\label{prop:trace-identity-root}
For every $n\ge2$ and every knot $K$, there is a unique
\[
 g_{K,n,1}(\mathbf t,\eta)
 \in\mathbb Z[t_1^{\pm1},\ldots,t_{n-1}^{\pm1},
                    D_{K,n}(\mathbf t)^{-1}][[\eta]],
 \qquad \eta=Q-1,
\]
whose joint color germ is the integral Newton knot germ. Its values
at dominant colors are the colored invariants, and
$g_{K,n,1}(\mathbf t,0)=D_{K,n}(\mathbf t)^{-1}$.
\end{proposition}

\begin{proof}
Define
\[
 g_{K,n,1}(\mathbf t,\eta)=
 \mathscr R_{K,n}((1+\eta)t_1,\ldots,(1+\eta)t_{n-1},
                           \log(1+\eta)).
\]
For each interval root, the denominator satisfies
\[
 \Delta_K((1+\eta)^{b-a}t_{a,b})
       =\Delta_K(t_{a,b})+O(\eta),
\]
so its inverse expands in powers of $\eta$ with coefficients in the
localization at $D_{K,n}$, and \eqref{eq:trace-identity-colored-input}
gives the required colored values.

Since $D_{K,n}(\mathbf1)=1$, this rational series has a joint Taylor
germ at $(\mathbf t,\eta)=(\mathbf1,0)$. Its difference from the
Newton germ vanishes on every curve
\[
 t_i=(1+\eta)^{r_i},\qquad r_i\in\N,
\]
so it is zero by Lemma~\ref{lem:rel-color-rigidity}. Applying
Lemma~\ref{lem:rel-integral-saturation} to each quantum
coefficient gives the integral coefficient ring, and uniqueness
follows from the injectivity of rational Taylor expansion. Finally,
setting $\eta=0$ in the definition gives
$g_{K,n,1}(\mathbf t,0)=[h^0]\mathscr R_{K,n}(\mathbf t)
=D_{K,n}(\mathbf t)^{-1}$.
\end{proof}

\begin{proof}[Proof of Theorem~\ref{main:type-a}]
Fix $n\ge2$ and a zero-framed knot $K$.
Theorem~\ref{thm:sym-dominant-knot-newton} gives
\[
 \mathcal J_{K,n}\in\mathcal N_{n-1}^+(A),\qquad
 \mathcal J_{K,n}(Q^{\mathbf r};Q)
       =J_{K,\sum r_i\omega_i}^{(n)}(Q).
\]
By Propositions~\ref{prop:trace-all-rank-trace-closure}
and~\ref{prop:trace-identity-root}, at each root $\zeta$ of order
$m\ge1$ its germ at the identity color is the germ of a series in
\[
 \Q(\zeta)[\mathbf t^{\pm1},
                     D_{K,n}(\mathbf t^m)^{-1}][[Q-\zeta]].
\]
The rational-germ condition of
Theorem~\ref{thm:rel-exact-image} therefore holds, and
that theorem supplies the integral local coefficients, the Frobenius
gluing and the unique relative lift.

For $m>1$, Proposition~\ref{prop:trace-all-rank-trace-closure} gives
\[
 D_{K,n}(\mathbf t^m)g_{K,n,\zeta}(\mathbf t,0)
                  \in\Z[\zeta][\mathbf t^{\pm1}],
\]
and for $m=1$ this numerator is one by
Proposition~\ref{prop:trace-identity-root}. Finally, the lift is
determined by the colored values, so it does not depend on the braid
presentation used in the proof.
\end{proof}

\begin{corollary}\label{cor:trace-all-rank-block-recurrences}
Every raw cyclotomic coefficient block of $\mathcal J_{K,n}$ lies in
\[
 \mathbb Z[\zeta]
 [(1+z_i)^{\pm1},D_{K,n}(1+\mathbf z)^{-1}].
\]
Equivalently, it satisfies the recurrence criterion of
Corollary~\ref{cor:rel-block-recurrence}. If $\Delta_K=1$, then in
every rank each quantum coefficient of each local component is a
Laurent polynomial in the independent colors.
\end{corollary}

\begin{proof}
Apply Theorem~\ref{thm:rel-cyclotomic-block-criterion} to the lift.
When $\Delta_K=1$ the denominator $D_{K,n}$ is one, so the local
coefficients are Laurent polynomials.
\end{proof}

\section{Consequences for knot invariants}
\label{sec:consequences}

\subsection{A uniform bound on pole order}
\label{subsec:cons-uniform-poles}

The operator estimates of Section~\ref{sec:braid} give a bound
on the order of the Alexander poles that does not depend on the order
of the root of unity. The differential order of a matrix of polynomial
differential operators is the maximum of the orders of its entries;
we write it $\ord_{\rm diff}$.

\begin{theorem}\label{thm:cons-uniform-pole-order}
Fix $n\ge2$ and a signed braid $\beta$ whose closure is the knot $K$.
There is an integer $C_{n,\beta}$, independent of $m>1$, of the
primitive $m$th root $\zeta$, and of $h\ge0$, such that
\begin{equation}\label{eq:cons-uniform-pole-order}
 D_{K,n}(t_1^m,\ldots,t_{n-1}^m)^{\,1+C_{n,\beta}h}
 [x^h]g_{K,n,\zeta}(\mathbf t,x)
 \in\mathbb Z[\zeta][t_1^{\pm1},\ldots,t_{n-1}^{\pm1}],
\end{equation}
where $x=Q-\zeta$.
\end{theorem}

The constant can be made explicit. Put
\[
 H=\frac{n(n^2-1)}6,\quad L=2H-1,\quad R=\frac{n(n-1)}2,
\]
\[
 a_n=L+2,\qquad b_n=RL,\qquad c_n=Lb_n=RL^2,
 \qquad A_n=a_n+b_n+c_n.
\]
Here $R$ is the number of positive roots and $L$ bounds the constants
$L_{n,m}$ of Corollary~\ref{cor:braid-uniform-residue-bound}. If the
chosen word for $\beta$ has $\ell_+$ positive and $\ell_-$ negative
crossings, set $B_{n,\beta}=\ell_+b_n+\ell_-c_n$. Then one may take
\begin{equation}\label{eq:cons-explicit-uniform-constant}
 C_{n,\beta}=2(A_n+B_{n,\beta}).
\end{equation}

\begin{proof}
We bound the differential order of the operators
$\mathcal D_{\beta,h}$ of Theorem~\ref{thm:braid-all-rank-insertion}, and
then convert this into a bound on the exponent of the trace
denominator.

For the root operators $V_r$ and $Z$ of
Section~\ref{subsec:braid-principal-crossing},
Theorem~\ref{thm:braid-fixed-transfer-jets} gives
\begin{equation}\label{eq:cons-VZ-order-bounds}
 \ord_{\rm diff}[\epsilon^j]V_r\le j,\qquad
 \ord_{\rm diff}[\epsilon^j]Z\le j+1.
\end{equation}
In particular, $Z_0=A_\alpha+N_\alpha$ has order at most one, and by
Corollary~\ref{cor:braid-uniform-residue-bound} a product of more than
$L_{n,m}\le L$ conjugates of $N_\alpha$ is zero.

Consider a Duhamel term of quantum order $e$ in
Lemma~\ref{lem:braid-filtered-insertion}, and suppose it contains $t$
insertions of positive order, of quantum orders $j_1,\ldots,j_t$. Then
$\sum j_i=e$ and $t\le e$, and by \eqref{eq:cons-VZ-order-bounds} these
insertions have total differential order at most
$\sum(j_i+1)=e+t$. Each of the $t+1$ gaps before, between and after
them contains at most $L$ conjugates of $N_\alpha$, each of order at
most one, and conjugation by the polynomial flow preserves
differential order. This gives
\begin{equation}\label{eq:cons-exp-order-bound}
 \ord_{\rm diff}\bigl(( [\epsilon^e]e^{Z(\epsilon)})
                                      e^{-A_\alpha}\bigr)
 \le e+t+(t+1)L\le a_ne+L.
\end{equation}

By Lemma~\ref{lem:braid-scalar-grouping}, for each residue $r$ the
grouping polynomial of quantum order $j$ has degree at most $2j$, so
the associated polynomial $T_{r,j}(Z)$ satisfies
\[
 \ord_{\rm diff}[\epsilon^\ell]T_{r,j}(Z(\epsilon))
       \le\ell+2j.
\]
A grouped root term of total order $h=j+v+\ell+e$, built from the
coefficient of $\epsilon^v$ in $V_r$, the coefficient of
$\epsilon^\ell$ in $T_{r,j}(Z)$ and the coefficient of $\epsilon^e$ in
$e^{Z}$, therefore has order at most
\[
 v+\ell+2j+a_ne+L\le a_nh+L,
\]
since $a_n\ge2$. Differential orders add under composition, and
conjugation by polynomial automorphisms preserves them, so the product
of the $R$ root factors adds $RL=b_n$ to the constant term of this
bound. The Cartan coefficient of order $j$ is an Euler operator of
order at most $2j\le a_nj$. Thus a positive checked crossing satisfies
\begin{equation}\label{eq:cons-positive-crossing-order}
 \ord_{\rm diff}\mathcal D_h^+\le a_nh+b_n.
\end{equation}

For an inverse crossing, write $\mathcal D_0=R_0+N_0$ as in the proof
of Theorem~\ref{thm:braid-all-rank-insertion}. The constant matrix
$R_0^{-1}$ has order zero, and $N_0$ has order at most $b_n$ and
lowers the residue height by at least $m$. Therefore
\[
 \mathcal D_0^{-1}
 =\sum_{j=0}^{L}(-R_0^{-1}N_0)^jR_0^{-1},
 \qquad \ord_{\rm diff}\mathcal D_0^{-1}\le Lb_n=c_n.
\]
For $h>0$, unwinding the recursion for $\mathcal E_h$ writes the
inverse coefficient as a sum over ordered compositions
$j_1+\cdots+j_s=h$ with $j_i>0$. Each summand contains $s\le h$
factors $\mathcal D_{j_i}$ of positive order and $s+1$ copies of
$\mathcal D_0^{-1}$, so by \eqref{eq:cons-positive-crossing-order}
\[
 \ord_{\rm diff}\le a_nh+sb_n+(s+1)c_n
                 \le A_nh+c_n.
\]
Conjugation by the inverse substitution preserves this order, and
together with the case $h=0$ treated above we obtain
\begin{equation}\label{eq:cons-negative-crossing-order}
 \ord_{\rm diff}\mathcal D_h^-\le A_nh+c_n
 \qquad(h\ge0).
\end{equation}

Composition along the braid word now gives
\[
 \ord_{\rm diff}\mathcal D_{\beta,h}
       \le A_nh+\ell_+b_n+\ell_-c_n
       =A_nh+B_{n,\beta}=:d_h.
\]
The pivot contributes, in order $j$, an Euler operator of order at
most $j$, so it preserves this bound, and the framing prefactor has
order zero.

We now pass to the trace. Lemma~\ref{lem:trace-general-cut} gives cut
terms $M_pC_\Psi\partial^\eta$ with $|\eta|\le d_h$. By
Corollary~\ref{cor:trace-two-sided-trace} and normal ordering, their
traces reduce to traces of $X^\mu\partial^\nu C_\Psi$ with
$|\nu|\le d_h$. For such a term, Lemma~\ref{lem:trace-finite-replacement}
allows only replacement indices $\kappa$ with
\[
 |\kappa|\le|\nu|-|\mu|\le d_h,
\]
and by the proof of Theorem~\ref{thm:trace-general-trace} the exponent of
the determinant in each root block is at most
\[
 1+|\nu_\alpha|+|\kappa_\alpha|\le1+2d_h.
\]
A common regulated denominator is therefore
\[
 \prod_{\alpha>0}
 \det\!\left(I-s^{m\hgt(\alpha)}
              B_\beta^{\mathrm{tr}}(T_\alpha)\right)^{1+2d_h},
\]
and summing over the finite residue basis does not increase the
exponent. At $s=1$ this becomes, up to Laurent units,
$D_{K,n}(\mathbf t^m)^{1+2A_nh+2B_{n,\beta}}$.

This denominator depends only on $\mathbf t=\mathbf u^2$ and has
integer coefficients, so it is unaffected by the two descents in the
proof of Proposition~\ref{prop:trace-all-rank-trace-closure}. Multiplying
the coefficient of order $h$ by it gives a Laurent polynomial over
$\Q(\zeta)$ whose Taylor germ is integral, because the local series is
the Newton knot germ, and Lemma~\ref{lem:rel-integral-saturation}
puts this numerator in $\Z[\zeta][\mathbf t^{\pm1}]$. The change of
coordinates $x=\zeta\epsilon$ only multiplies the coefficient of order
$h$ by a cyclotomic unit. For $h\ge1$,
\[
 1+2A_nh+2B_{n,\beta}
       \le1+2(A_n+B_{n,\beta})h,
\]
and for $h=0$ the denominator divides one Alexander product by
Proposition~\ref{prop:trace-all-rank-trace-closure}. This proves the
theorem with the constant \eqref{eq:cons-explicit-uniform-constant}.
\end{proof}

\subsection{Alexander recovery and the inverse Cartan matrix}

Write
\[
 \mathcal J_{K,n}(\mathbf t;Q)
   =\sum_{\alpha\in\mathbb N^{n-1}}c_\alpha(Q)
                         P_\alpha(\mathbf t;Q).
\]
Let $C_{n-1}$ be the Cartan matrix of $A_{n-1}$, and let $\mathbf e_i$
be the coordinate unit vectors in $\mathbb Z^{n-1}$.

\begin{corollary}\label{cor:cons-alexander-hessian}
For each $n\ge2$, the Newton knot element of type $A_{n-1}$
determines the normalized Alexander polynomial, and its Newton
coefficients of degree one and two satisfy
\begin{align*}
 c_{\mathbf e_i}(1)&=0,\\
 c_{2\mathbf e_i}(1)&=-\frac{\Delta_K''(1)}2\,i(n-i),\\
 c_{\mathbf e_i+\mathbf e_j}(1)&=-\Delta_K''(1)\,
        \min(i,j)\bigl(n-\max(i,j)\bigr)\qquad(i\ne j).
\end{align*}
Equivalently, the color Hessian of the constant quantum coefficient is
\begin{equation}\label{eq:cons-inverse-cartan-hessian}
 \left.
 \left(\frac{\partial^2 g_{K,n,1}(\mathbf t,0)}
                  {\partial t_i\partial t_j}\right)_{i,j=1}^{n-1}
 \right|_{\mathbf t=\mathbf1}
   =-n\Delta_K''(1)\,C_{n-1}^{-1}.
\end{equation}
In particular, for $i<j$,
\begin{align}
 J_{K,\omega_i+\omega_j}^{(n)}
 -J_{K,\omega_i}^{(n)}-J_{K,\omega_j}^{(n)}+1
 &\equiv-\Delta_K''(1)i(n-j)(Q-1)^2,
 \label{eq:cons-mixed-color-congruence}\\
 J_{K,2\omega_i}^{(n)}
 -(Q+1)J_{K,\omega_i}^{(n)}+Q
 &\equiv-\Delta_K''(1)i(n-i)(Q-1)^2,
 \label{eq:cons-diagonal-color-congruence}
\end{align}
where both congruences are modulo $(Q-1)^3A$.
\end{corollary}

\begin{proof}
By Theorem~\ref{main:type-a},
$g_{K,n,1}(\mathbf t,0)=D_{K,n}(\mathbf t)^{-1}$. On the $i$th
coordinate axis, with $t_i=t$ and the other coordinates equal to one,
this is $\Delta_K(t)^{-i(n-i)}$, because exactly $i(n-i)$ interval
roots contain $\alpha_i$. Since $\Delta_K(1)=1$, the polynomial
$\Delta_K$ is recovered as the unique power series root with value one
at $t=1$.

Alexander symmetry gives $\Delta_K'(1)=0$, so
$\Delta_K(1+y)=1+\tfrac12\Delta_K''(1)y^2+O(y^3)$. For an interval
root, $t_a\cdots t_{b-1}-1=w_a+\cdots+w_{b-1}+O((\mathbf w)^2)$ with
$t_i=1+w_i$, and hence, up to total degree two in the colors,
\[
 D_{K,n}(1+\mathbf w)^{-1}
 =1-\frac{\Delta_K''(1)}2
    \sum_{1\le a<b\le n}
         (w_a+\cdots+w_{b-1})^2+O((\mathbf w)^3).
\]
The number of interval roots containing both positions $i$ and $j$ is
$\min(i,j)(n-\max(i,j))$. At $Q=1$ the Newton basis consists of the
monomials $\mathbf w^\alpha$, so reading off the coefficients of
$w_i^2$ and $w_iw_j$ gives the formulas for $c_{2\mathbf e_i}(1)$ and
$c_{\mathbf e_i+\mathbf e_j}(1)$, and the absence of linear terms
gives $c_{\mathbf e_i}(1)=0$. The inverse Cartan matrix is
\[
 (C_{n-1}^{-1})_{ij}
   =\frac{\min(i,j)(n-\max(i,j))}{n},
\]
as one checks by taking second differences in $i$, and this gives
\eqref{eq:cons-inverse-cartan-hessian}. The diagonal coefficients are
integers: writing $\Delta_K(t)=a_0+\sum_{k>0}a_k(t^k+t^{-k})$ gives
$\Delta_K''(1)=2\sum_{k>0}a_k k^2$.

For the congruences we use $J^{(n)}_{K,0}=1$. By
Proposition~\ref{prop:newton-divided-differences}, the left side of
\eqref{eq:cons-mixed-color-congruence} equals
$(Q-1)^2c_{\mathbf e_i+\mathbf e_j}(Q)$, and that of
\eqref{eq:cons-diagonal-color-congruence} equals
$(Q^2-1)(Q^2-Q)c_{2\mathbf e_i}(Q)$. Their terms of order two in
$Q-1$ are $c_{\mathbf e_i+\mathbf e_j}(1)$ and $2c_{2\mathbf e_i}(1)$,
which gives the congruences.
\end{proof}

\subsection{Knot operations and difference equations}

\begin{proposition}\label{prop:cons-knot-operations}
For knots $K,L$ and every $n\ge2$,
\[
 \mathscr J_{K\mathbin{\#}L,n}
     =\mathscr J_{K,n}\mathscr J_{L,n}
 \quad\text{in}\quad
 \mathcal H_{\rm rel}(D_{K,n}D_{L,n})_{\rm pol}.
\]
The substitution
\[
 \iota(F)(\mathbf t;Q)=F(\mathbf t^{-1};Q^{-1})
\]
defines an involution of the dominant Newton ring, which preserves
the relative image for $D_{K,n}$. If $\overline K$ is the mirror
knot, then
\[
 \mathcal J_{\overline K,n}=\iota(\mathcal J_{K,n}),
 \qquad
 g_{\overline K,n,\zeta}(\mathbf t,x)
  =g_{K,n,\zeta^{-1}}
     \left(\mathbf t^{-1},(\zeta+x)^{-1}-\zeta^{-1}\right).
\]
Moreover, $\mathcal J_{K\#L,n}=\mathcal J_{K\#M,n}$ implies
$\mathcal J_{L,n}=\mathcal J_{M,n}$.
\end{proposition}

\begin{proof}
The invariant of a reduced $(1,1)$-tangle is multiplicative under
composition, and Seifert matrices of a connected sum are block sums.
Hence
\[
 J_{K\#L,\lambda}^{(n)}=J_{K,\lambda}^{(n)}J_{L,\lambda}^{(n)},
 \qquad D_{K\#L,n}=D_{K,n}D_{L,n}.
\]
The componentwise product $\mathscr J_{K,n}\mathscr J_{L,n}$ is a
relative collection for $D_{K,n}D_{L,n}$, since the coefficient rings
and the gluing conditions are compatible with products, and its color
values are those of $K\#L$. Uniqueness in
Theorem~\ref{main:recognition} gives the connected-sum identity.

The identity
\[
 P_k(t^{-1};Q^{-1})
      =(-1)^kt^{-k}Q^{-\binom k2}P_k(t;Q)
\]
shows that simultaneous inversion preserves every box ideal, so it
defines a continuous involution (Proposition~\ref{prop:newton-operations}).
Under dominant color evaluation it sends $J_{K,\lambda}^{(n)}(Q)$ to
$J_{K,\lambda}^{(n)}(Q^{-1})$, which is the invariant of the mirror
knot by the Reshetikhin--Turaev construction. Injectivity of Newton
evaluation gives $\mathcal J_{\overline K,n}=\iota(\mathcal J_{K,n})$.

At $Q=\zeta+x$ we have $Q^{-1}=\zeta^{-1}+((\zeta+x)^{-1}-\zeta^{-1})$,
so the germ of $\iota(F)$ at $\zeta$ is obtained from the germ of $F$
at $\zeta^{-1}$ by the change of variables
\[
 (\mathbf t,x)\longmapsto
 \left(\mathbf t^{-1},(\zeta+x)^{-1}-\zeta^{-1}\right).
\]
Alexander symmetry gives $D_{K,n}(\mathbf t^{-m})=D_{K,n}(\mathbf t^m)$,
so this change of variables preserves the permitted denominator. This
gives the mirror formula, and
Theorem~\ref{thm:rel-exact-image} shows that $\iota$
preserves the relative image. Finally, the Newton ring is a domain and
$\mathcal J_{K,n}$ has value one at the trivial color, so this element
can be cancelled.
\end{proof}

Quantum knot invariants in independent highest-weight coordinates
satisfy $q$-holonomic systems \cite[Theorem~6]{GaroufalidisLeQHolonomic}.
To relate colored recurrences to the local rational expansions, write
$E_i f(\mathbf t,Q)=f(t_1,\ldots,Qt_i,\ldots,t_{n-1},Q)$.

\begin{proposition}\label{prop:cons-recurrence-transfer}
Fix $n\ge2$. Let
$P=\sum_{\mathbf v\in S}a_{\mathbf v}(\mathbf t,Q)E^{\mathbf v}$,
where $S\subset\mathbb Z^{n-1}$ is finite,
$a_{\mathbf v}\in A[t_1^{\pm1},\ldots,t_{n-1}^{\pm1}]$, and
$E^{\mathbf v}f(\mathbf t,Q)=f(Q^{\mathbf v}\mathbf t,Q)$.
Then $P$ annihilates every local component of $\mathscr J_{K,n}$ if
and only if the corresponding recurrence holds for the colored
polynomials at every color $\mathbf r$ for which all the colors
$\mathbf r+\mathbf v$ are dominant. Moreover, it suffices that the
colored recurrence hold for all sufficiently large dominant colors.
\end{proposition}

\begin{proof}
At a root $\zeta$ of order $m$ consider the shifted local term
\[
 g_{K,n,\zeta}((\zeta+x)^{\mathbf v}\mathbf t,x).
\]
For shifts of either sign it lies in the same localization of quantum
series, because $\zeta^{m\mathbf v}=\mathbf1$ and so the shifted
denominator has constant term $D_{K,n}(\mathbf t^m)$. The shifts use
the same variable $Q$ in all charts, so they are compatible with
re-expansion and act on the whole relative ring. They need not
preserve the subring with Laurent polynomial color values, since a
negative shift moves boundary colors out of the dominant range.

If the colored recurrence holds for large colors, the local expression
$Pg_{K,n,\zeta}$ vanishes on every curve $t_i=Q^{ma_i}$ with all $a_i$
sufficiently large. The proof of Lemma~\ref{lem:rel-color-rigidity}
only uses that the $a_i$ range over infinite sets, so it applies to
this grid and shows that the germ of $Pg_{K,n,\zeta}$ at the identity
color vanishes. Injectivity of rational Taylor expansion then gives
$Pg_{K,n,\zeta}=0$. Conversely, a local identity can be specialized at
every admissible dominant color, and a Laurent polynomial with zero
Taylor expansion is zero.
\end{proof}

\subsection{The first-order term of a difference equation}

Let $d\ge1$, let $\mathcal R=\mathbb Q[t_1^{\pm1},\ldots,t_d^{\pm1}]$,
and put $\theta_i=t_i\partial_{t_i}$. The following algebraic
statement applies to the knot element of type $A_{n-1}$, with $d=n-1$
and $D=D_{K,n}$, since Proposition~\ref{prop:trace-identity-root}
provides the required leading term.

\begin{theorem}\label{thm:cons-logarithmic-annihilator}
Let $0\ne D\in\mathcal R$. Suppose
\[
 F(\mathbf t,Q)=D(\mathbf t)^{-1}+(Q-1)F_1(\mathbf t)+\cdots
 \in\operatorname{Frac}(\mathcal R)[[Q-1]]
\]
is annihilated by the left-normal-ordered difference operator
\[
 P=\sum_{v\in S}a_v(\mathbf t,Q)E^v,
 \qquad E^vF(\mathbf t,Q)
 =F(Q^{v_1}t_1,\ldots,Q^{v_d}t_d,Q),
\]
where $S\subset\Z^d$ is finite and
$a_v\in\operatorname{Frac}(\mathcal R)(Q)$. Assume that each $a_v$
is regular at $Q=1$, with $a_v(\mathbf t,1)$ and
$\partial_Qa_v(\mathbf t,1)$ in $\mathcal R$, and put
\[
 p_0(\mathbf t,\mathbf L)=\sum_v a_v(\mathbf t,1)\mathbf L^v,
 \quad A_i=\left.\partial_{L_i}p_0\right|_{\mathbf L=\mathbf1},
 \quad B=\sum_v\partial_Qa_v(\mathbf t,1).
\]
Then
\begin{equation}\label{eq:cons-logarithmic-constraint}
 p_0(\mathbf t,\mathbf1)=0,
 \qquad BD=\sum_iA_i\theta_iD.
\end{equation}
In particular $V=\sum_iA_i\theta_i$ preserves the ideal $(D)$, that
is, it is a logarithmic derivation along $D=0$, and
$\delta\mid V(\delta)$ in $\mathcal R$ for every irreducible factor
$\delta$ of $D$, whatever its multiplicity.
\end{theorem}

\begin{proof}
The constant term of $PF$ is $p_0(\mathbf t,\mathbf1)D^{-1}$, which
gives the first identity. With $x=Q-1$ we have
$Q^{v_i}t_i=t_i(1+v_ix+O(x^2))$ for shifts of either sign, and hence
\[
 E^vF=F_0+x\left(F_1+\sum_i v_i\theta_iF_0\right)+O(x^2),
 \qquad F_0=D^{-1}.
\]
Since $\sum_va_v(\mathbf t,1)v_i=A_i$, the coefficient of $x$ in $PF$
is $p_0(\mathbf t,\mathbf1)F_1+V(F_0)+BF_0$. The first term vanishes,
and substituting $F_0=D^{-1}$ and multiplying by $D^2$ gives
$-V(D)+BD=0$.

For the last assertion write $D=\delta^e U$ with $\delta\nmid U$.
Then $V(D)=e\delta^{e-1}UV(\delta)+\delta^eV(U)$, and $D\mid V(D)$
gives $\delta\mid eUV(\delta)$. As $\mathcal R$ is a unique
factorization domain and $e$ is invertible in $\Q$, it follows that
$\delta\mid V(\delta)$.
\end{proof}

Two remarks on the hypotheses. If the coefficients $a_v$ are Laurent
polynomials, dividing all of them by the largest power of $Q-1$ that
divides them does not affect the annihilation property and ensures
that not all $a_v(\mathbf t,1)$ vanish. At a nontrivial root $\zeta$ of
order $m$, the same expansion applies to operators in the shifts
$E_i^m$, which fix the center of the expansion: in the coordinate
$Q=\zeta(1+\epsilon)$ the first-order term of $E_i^m$ is $m\theta_i$.

\section{Explicit knot computations and cyclotomic expansions}
\label{sec:examples}

This section collects explicit calculations. The figure-eight knot
gives integral Newton coefficients and rational expansions at every
root of unity, torus knots compare independent Dynkin directions, and
connected sums and the Kinoshita--Terasaka knot show what the
classical color germ does and does not detect. Throughout we use the
reduced, zero-framed normalization and the Newton basis
$P_k(t;Q)=\prod_{j<k}(t-Q^j)$.

\subsection{The figure-eight knot}
\label{subsec:ex-figure-eight}

For the figure-eight knot $4_1$, take $n=2$ and $t=Q^r$ for the
highest weight $r\omega_1$. Its normalized Alexander polynomial is
\[
 \Delta(T)=3-T-T^{-1}.
\]
Habiro's formula \cite[Sections~3.3--3.4, equations~(5) and~(9)]{Habiro2002}
gives
\begin{equation}\label{eq:ex-figure-eight-series}
 \mathcal J_{4_1,2}(t;Q)=\sum_{k\ge0}K_k(t,Q),\qquad
 K_k(t,Q)=\prod_{j=1}^{k}
       \bigl(Qt+Q^{-1}t^{-1}-Q^j-Q^{-j}\bigr),
\end{equation}
where the dimension color of that formula is $Q^{r+1}=Qt$. Each factor
satisfies
\[
 Qt+Q^{-1}t^{-1}-Q^j-Q^{-j}
   =Q(t-Q^{j-1})(1-Q^{-j-1}t^{-1}),
\]
so $K_k\in P_k(t;Q)A[t^{\pm1}]$, the series converges in
$\mathcal N_1^+(A)$, and at $t=Q^r$ its summands with $k>r$ vanish.
For example,
\[
 J_{4_1,\omega_1}^{(2)}
     =Q^2-Q+1-Q^{-1}+Q^{-2}.
\]
Writing $\mathcal J_{4_1,2}=\sum_{k\ge0}c_k(Q)P_k(t;Q)$ and
interpolating the colors $r=0,1,2$ gives
\begin{equation}\label{eq:ex-figure-eight-newton}
 c_0=1,\qquad c_1=Q-Q^{-2},\qquad
 c_2=Q^2-Q^{-2}-Q^{-3}+Q^{-4}+Q^{-7}.
\end{equation}
Thus $c_1(1)=0$ and $c_2(1)=1$, in agreement with
$\Delta''(1)=-2$ and Corollary~\ref{cor:cons-alexander-hessian}.

At the singular color $t=Q^{-1}$, the boundary value of
Section~\ref{subsec:newton-walls} is
\begin{equation}\label{eq:ex-figure-eight-singular}
 \left.\mathcal J_{4_1,2}(t;Q)\right|_{t=Q^{-1}}
   =\sum_{k\ge0}(-1)^kQ^{-k(k+1)/2}(Q;Q)_k^2
     \in\mathcal H.
\end{equation}
This follows from $2-Q^j-Q^{-j}=-Q^{-j}(1-Q^j)^2$, and the right side
converges in the ordinary Habiro topology. At $Q=\zeta$ of order $m$
it terminates at $k=m-1$ and equals
$J_{4_1,(m-1)\omega_1}^{(2)}(\zeta)$, since that color also has
$Qt=1$. For $m=1,2,3,4$ (with $\zeta=i$ in the last case) the values
are $1,5,13,27$.

\subsection{Torus knots in type \texorpdfstring{$A_2$}{A2}}
\label{subsec:ex-torus-a2}

For $\mathfrak{sl}_3$, the color $r\omega_1+s\omega_2$ corresponds
to the partition $(r+s,s,0)$. The values at $\omega_1$, $2\omega_1$ and
$\omega_1+\omega_2$ determine the linear, the diagonal quadratic and
the mixed quadratic Newton coefficients, respectively.

We use the Rosso--Jones formula \cite{RJ}, in the zero-framed
normalization of \cite[equations~(2), (4), (5), and (9)]{BM}, with
its parameter $q$ replaced by $q^{-1}$.  Put $Q=q^2$ and
$[a]=(q^a-q^{-a})/(q-q^{-1})$.  If
$\psi_m(s_\lambda)=\sum_\mu c^{(m)}_{\lambda\mu}s_\mu$, then
\begin{equation}\label{eq:ex-rosso-jones}
 J^{T(m,b)}_\lambda(q)=
 \frac{\sum_\mu c^{(m)}_{\lambda\mu}
 q^{\frac bm C_\mu-mbC_\lambda}\dim_q\mu}{\dim_q\lambda},
 \qquad
 C_\lambda=\sum_i\lambda_i(\lambda_i-2i+1)+n|\lambda|,
 \qquad \gcd(m,b)=1.
\end{equation}
The correction $-|\lambda|^2/n$ from $\mathfrak{gl}_n$ to
$\mathfrak{sl}_n$ cancels in this exponent, since
$|\mu|=m|\lambda|$. In this convention the fundamental
$\mathfrak{sl}_2$ invariant of the trefoil is
$Q^{-1}+Q^{-3}-Q^{-4}$.

For an odd integer $\ell\ge1$, set
\[
 F_\ell=J^{T(2,\ell)}_{1,0}=J^{T(2,\ell)}_{0,1},\qquad
 B_\ell=J^{T(2,\ell)}_{2,0}=J^{T(2,\ell)}_{0,2},\qquad
 A_\ell=J^{T(2,\ell)}_{1,1}.
\]
The equalities under interchange of the two Dynkin coordinates follow
from duality. The three Adams decompositions in three variables are
\[
 \psi_2(s_1)=s_2-s_{11},\qquad
 \psi_2(s_2)=s_4-s_{31}+s_{22},\qquad
 \psi_2(s_{21})=s_{42}-s_{411}-s_{33}+s_{222};
\]
for the second Adams operation in arbitrary $\mathfrak{sl}_3$ colors,
see \cite[Theorem~2.4]{GMV}. With
$\dim_q(a,b,c)=[a-b+1][b-c+1][a-c+2]/[2]$ this gives
\begin{align}
 F_\ell&=Q^{-\ell}(Q+Q^{-1})-Q^{-2\ell},
 \label{eq:ex-torus-family-fundamental}\\
 B_\ell&=
 \frac{q^{-4\ell}[5][6]/[2]
       -q^{-8\ell}[3][5]+q^{-10\ell}[3][4]/[2]}
      {[3][4]/[2]},
 \label{eq:ex-torus-family-symmetric}\\
 A_\ell&=
 \frac{q^{-4\ell}[3]^2[6]/[2]
       -2q^{-6\ell}[4][5]/[2]+q^{-12\ell}}
      {[2][4]}.
 \label{eq:ex-torus-family-adjoint}
\end{align}
These formulas hold for every odd $\ell$, including
$5_1=T(2,5)$ and $7_1=T(2,7)$.

Their integrality can be seen directly from the expressions in $Q$.
Write $U_k(Q)=1+Q+\cdots+Q^{k-1}$ for $k\ge1$ and $U_0=0$. In terms of
$Q$ the last two formulas become
\begin{align}
 B_\ell&=Q^{-5\ell}
 \frac{U_5(Q)\bigl(Q^{3\ell-2}(Q^2-Q+1)-Q^{\ell-1}\bigr)
                      +Q^2+1}{Q^2+1},
 \label{eq:ex-torus-family-symmetric-Q}\\
 A_\ell&=Q^{2-6\ell}
 \frac{Q^{4\ell-4}U_3(Q)^2(Q^4+Q^2+1)
       -2Q^{3\ell-3}(Q^2+1)U_5(Q)+1}
      {(Q+1)^2(Q^2+1)}.
 \label{eq:ex-torus-family-adjoint-Q}
\end{align}
Modulo $Q^2+1$ the numerator of $B_\ell$ vanishes because
$Q^{2\ell}=-1$, and the numerator of $A_\ell$ reduces to $-1+1=0$.
At $Q=-1$ the numerator of $A_\ell$ and its derivative both vanish.
The displayed denominators therefore cancel in $\Z[Q]$, and
$F_\ell,B_\ell,A_\ell\in A$.

The corresponding Newton coefficients are
\begin{align}
 c_{1,0}^{(\ell)}=c_{0,1}^{(\ell)}
   &=-Q^{-2\ell}(Q-1)U_{\ell-1}(Q)U_{\ell+1}(Q),
 \label{eq:ex-torus-family-linear-newton}\\
 c_{1,1}^{(\ell)}
   &=\frac{A_\ell-2F_\ell+1}{(Q-1)^2},
 \label{eq:ex-torus-family-mixed-newton}\\
 c_{2,0}^{(\ell)}=c_{0,2}^{(\ell)}
   &=\frac{B_\ell-(Q+1)F_\ell+Q}{(Q^2-1)(Q^2-Q)}.
 \label{eq:ex-torus-family-diagonal-newton}
\end{align}
The linear formula follows from
$F_\ell-1=-Q^{-2\ell}(Q^{\ell-1}-1)(Q^{\ell+1}-1)$. All three colored
invariants have value one and derivative zero at $Q=1$, so the
numerators in the last two lines are divisible by $(Q-1)^2$. Moreover
\eqref{eq:ex-torus-family-symmetric-Q} gives $B_\ell(-1)=1$, so the
numerator in the last line also vanishes at $Q=-1$. Since
$(Q^2-1)(Q^2-Q)=Q(Q-1)^2(Q+1)$, this shows directly that
$c_{1,1}^{(\ell)},c_{2,0}^{(\ell)}\in A$ for every odd $\ell$.

For the trefoil, $\ell=3$ gives
\begin{equation}\label{eq:ex-trefoil-fundamental}
 F_3=Q^{-2}+Q^{-4}-Q^{-6},
\end{equation}
\begin{equation}\label{eq:ex-trefoil-adjoint}
 \begin{split}
 A_3={}&Q^{-4}+2Q^{-6}-2Q^{-7}+2Q^{-8}-2Q^{-9}
       +2Q^{-10}-4Q^{-11}\\
      &+3Q^{-12}-2Q^{-13}+2Q^{-14}-2Q^{-15}+Q^{-16},
 \end{split}
\end{equation}
and the diagonal color gives
\[
 B_3=Q^{-4}+Q^{-7}+Q^{-8}-Q^{-9}-Q^{-12}-Q^{-13}+Q^{-15}.
\]
The mixed and diagonal Newton coefficients are
\begin{equation}\label{eq:ex-trefoil-mixed-coefficient}
 \begin{split}
 c_{1,1}^{(3)}={}&Q^{-2}+2Q^{-3}+Q^{-4}-2Q^{-6}-4Q^{-7}
                   -2Q^{-8}-2Q^{-9}\\
                &+2Q^{-12}+Q^{-14}+Q^{-16}
 \end{split}
\end{equation}
and
\[
 c_{2,0}^{(3)}=
 \frac{Q^{13}+Q^{12}+Q^{11}-Q^9-2Q^8-2Q^7-2Q^6-Q^5
                         +Q^2+Q+1}{Q^{16}}.
\]
They are different Laurent polynomials, although their values at $Q=1$
agree.

The comparison holds uniformly in $\ell$. Put $x=Q-1$ and
$\kappa_\ell=\ell^2-1$. Expanding
\eqref{eq:ex-torus-family-mixed-newton}--\eqref{eq:ex-torus-family-diagonal-newton}
gives
\begin{align}
 c_{1,1}^{(\ell)}(1+x)
  &=-\frac{\kappa_\ell}{4}
    +\frac{\kappa_\ell(\ell+1)}4x
    +\frac{\kappa_\ell(83\ell^2-18\ell-57)}{48}x^2+O(x^3),
 \label{eq:ex-torus-family-mixed-jets}\\
 c_{2,0}^{(\ell)}(1+x)
  &=-\frac{\kappa_\ell}{4}
    +\frac{\kappa_\ell(2\ell+9)}8x
    +\frac{\kappa_\ell(53\ell^2-60\ell-141)}{48}x^2+O(x^3).
 \label{eq:ex-torus-family-diagonal-jets}
\end{align}
To obtain these expansions, set $Q=e^h$, use
\[
 [k]=k\left(1+\frac{k^2-1}{24}h^2
      +\frac{3k^4-10k^2+7}{5760}h^4\right)+O(h^6)
\]
in \eqref{eq:ex-torus-family-symmetric}--\eqref{eq:ex-torus-family-adjoint},
and then substitute $h=\log(1+x)$ in the Newton quotients. The
normalized Alexander polynomial
\[
 \Delta_{T(2,\ell)}(T)=
 T^{-(\ell-1)/2}\frac{T^\ell+1}{T+1}
\]
has $\Delta''(1)=\kappa_\ell/4$. The constant terms therefore agree
with Corollary~\ref{cor:cons-alexander-hessian}, whereas
\[
 (c_{2,0}^{(\ell)})'(1)-(c_{1,1}^{(\ell)})'(1)
     =\frac{7\kappa_\ell}{8}
\]
is nonzero for every nontrivial knot in the family: the diagonal and
mixed directions are already distinguished at first order in $Q-1$.

For $T(3,4)$, the third Adams operation in three variables is
\begin{align*}
 \psi_3(s_1)&=s_3-s_{21}+s_{111},\\
 \psi_3(s_2)&=s_6-s_{51}+s_{411}+s_{33}-s_{321}+s_{222},\\
 \psi_3(s_{21})&=s_{63}-s_{621}-s_{54}+s_{522}+s_{441}
                              -s_{432}+2s_{333};
\end{align*}
these decompositions, including the multiplicity two of $s_{333}$,
follow from the three-variable Schur determinants.
Formula~\eqref{eq:ex-rosso-jones} now gives
\begin{align}
 F_{3,4}:=J^{T(3,4)}_{1,0}
   &=\frac{q^{-16}[4][5]/[2]-q^{-24}[2][4]+q^{-32}}{[3]}
 \notag\\
   &=Q^{-6}+Q^{-8}+Q^{-9}-Q^{-11}-Q^{-13}-Q^{-14}+Q^{-15},
 \label{eq:ex-torus-34-fundamental}\\
 B_{3,4}:=J^{T(3,4)}_{2,0}
   &=\frac{q^{-32}[7][8]/[2]-q^{-48}[5][7]
       +2q^{-64}[4][5]/[2]-q^{-72}[2][4]+q^{-80}
          }{[3][4]/[2]},
 \label{eq:ex-torus-34-symmetric}\\
 A_{3,4}:=J^{T(3,4)}_{1,1}
   &=\frac{q^{-32}[4]^2[8]/[2]-2q^{-40}[5][7]
       +2q^{-56}[4][5]/[2]-q^{-64}[2][4]+2q^{-72}
          }{[2][4]}.
 \label{eq:ex-torus-34-adjoint}
\end{align}
The first Newton coefficient is $(F_{3,4}-1)/(Q-1)$, and the
quadratic coefficients are given by
\eqref{eq:ex-torus-family-mixed-newton}--\eqref{eq:ex-torus-family-diagonal-newton}
with $F_\ell,B_\ell,A_\ell$ replaced by $F_{3,4},B_{3,4},A_{3,4}$.
Their first terms are
\begin{align*}
 c_{1,1}^{T(3,4)}(1+x)&=-10+70x+2590x^2+O(x^3),\\
 c_{2,0}^{T(3,4)}(1+x)&=-10+105x+1335x^2+O(x^3).
\end{align*}
Here
\[
 \Delta_{T(3,4)}(T)=T^3-T^2+1-T^{-2}+T^{-3},\qquad
 \Delta_{T(3,4)}''(1)=10,
\]
so the classical coefficient again agrees with the Alexander Hessian.
The divisions in \eqref{eq:ex-torus-34-symmetric} and
\eqref{eq:ex-torus-34-adjoint}, as well as the Newton divisions, are
exact and give Laurent polynomials with integer coefficients. The
presentation $T(4,3)$ of the same knot, evaluated with $\psi_4$, gives
the same three Laurent polynomials, which confirms the conventions in
\eqref{eq:ex-rosso-jones}.

The first quantum coefficients for these four knots are
\[
\begin{array}{c|r|rr|rr}
 K&\Delta_K''(1)&c_{1,1}(1)&c_{1,1}'(1)&c_{2,0}(1)&c_{2,0}'(1)\\ \hline
 T(2,3)&2&-2&8&-2&15\\
 T(2,5)&6&-6&36&-6&57\\
 T(2,7)&12&-12&96&-12&138\\
 T(3,4)&10&-10&70&-10&105
\end{array}
\]

\subsection{Torus knots in type \texorpdfstring{$A_3$}{A3}}
\label{subsec:ex-sl4}

For $\mathfrak{sl}_4$, the color $r_1\omega_1+r_2\omega_2+r_3\omega_3$
corresponds to the partition $(r_1+r_2+r_3,r_2+r_3,r_3,0)$. The six
colors
\[
 \omega_1,\quad\omega_2,\quad2\omega_1,\quad2\omega_2,
 \quad\omega_1+\omega_2,\quad\omega_1+\omega_3
\]
therefore correspond to the partitions
$(1),(1,1),(2),(2,2),(2,1),(2,1,1)$, and the Dynkin diagram involution
supplies the colors obtained by interchanging $\omega_1$ and
$\omega_3$.

For $T(2,\ell)$ with $\ell\ge1$ odd, we use
\eqref{eq:ex-rosso-jones} with $n=4$ and
\[
 \dim_q\lambda=
 \prod_{1\le i<j\le4}
       \frac{[\lambda_i-\lambda_j+j-i]}{[j-i]},
 \qquad \lambda_i=0\quad(i>\operatorname{length}(\lambda)).
\]
The required Adams identities in four variables are
\begin{align*}
 \psi_2(s_1)&=s_2-s_{11},\\
 \psi_2(s_{11})&=s_{22}-s_{211}+s_{1111},\\
 \psi_2(s_2)&=s_4-s_{31}+s_{22},\\
 \psi_2(s_{22})&=s_{44}-s_{431}+s_{422}+s_{3311}-s_{3221}+s_{2222},\\
 \psi_2(s_{21})&=s_{42}-s_{411}-s_{33}+s_{3111}+s_{222}-s_{2211},\\
 \psi_2(s_{211})&=s_{422}-s_{4211}-s_{332}+s_{3311}+s_{2222}.
\end{align*}
Unlike the three-variable case, terms with four rows now occur and
must be kept. Together with the quantum-dimension product, these
decompositions determine the six invariants for every odd $\ell$.

To extract Newton coefficients, put
$F_i=J_{\omega_i}$, $B_i=J_{2\omega_i}$ and
$A_{ij}=J_{\omega_i+\omega_j}$.  Then
\[
 c_{2\mathbf e_i}=
    \frac{B_i-(Q+1)F_i+Q}{Q(Q-1)^2(Q+1)},\qquad
 c_{\mathbf e_i+\mathbf e_j}=
    \frac{A_{ij}-F_i-F_j+1}{(Q-1)^2}.
\]
Expanding the finite sums with $Q=e^h$ gives, for the six partitions
in the order listed above,
\[
 J_\lambda^{T(2,\ell)}(e^h)
  =1-a_\lambda(\ell^2-1)h^2
       +\frac{4\ell}{3}a_\lambda(\ell^2-1)h^3+O(h^4),
 \qquad
 (a_\lambda)=
 \left(\frac{15}{8},\frac52,\frac92,6,\frac{39}{8},4\right).
\]
Substituting $h=\log(1+x)$ proves the first two columns of the table
below for every odd $\ell$. The last two columns give the triples
$(c_\alpha(1),c_\alpha'(1),[x^2]c_\alpha(1+x))$ for the indicated
knots, obtained by expanding the same exact sums one order further.
In the table $\ell\ge3$ and $\kappa=\ell^2-1$; for $\ell=1$ all four
Newton coefficients vanish.
\[
\begin{array}{c|cc|rr}
 \alpha&c_\alpha(1)/\kappa&c_\alpha'(1)/\kappa
                    &T(2,3)&T(2,5)\\\hline
 2\mathbf e_1
   &-3/8&(4\ell+15)/8&(-3,27,116)&(-9,105,1455)\\
 2\mathbf e_2
   &-1/2&(4\ell+15)/6&(-4,36,247)&(-12,140,2773)\\
 \mathbf e_1+\mathbf e_2
   &-1/2&(4\ell+3)/6&(-4,20,451)&(-12,92,4089)\\
 \mathbf e_1+\mathbf e_3
   &-1/4&(4\ell+3)/12&(-2,10,306)&(-6,46,2748)
\end{array}
\]
Since $\Delta_{T(2,\ell)}''(1)=\kappa/4$, the constant terms reproduce
the $A_3$ inverse Cartan pattern of
Corollary~\ref{cor:cons-alexander-hessian}. The middle diagonal
direction $2\mathbf e_2$ and the adjacent mixed direction
$\mathbf e_1+\mathbf e_2$ have the same classical coefficient
$-\kappa/2$, but their first quantum coefficients differ by $2\kappa$.

\subsection{Cyclotomic expansions and coefficient blocks}
\label{subsec:ex-cyclotomic}

We return to the figure-eight series
\eqref{eq:ex-figure-eight-series}, with $n=2$ and
$\Delta(T)=3-T-T^{-1}$. Its local components are obtained by grouping
the summands into complete cyclotomic blocks. All-order pole estimates
for the symmetric-color expansion are given in
\cite[Theorem~9.8]{FangNewton}.

\begin{samepage}
\begin{proposition}\label{prop:ex-figure-eight-all-roots}
Let $\zeta$ be a primitive $m$th root of unity, including $m=1$,
and set
\[
 f_j(t)=\zeta t+\zeta^{-1}t^{-1}-\zeta^j-\zeta^{-j},\qquad
 S_m(t)=\sum_{b=0}^{m-1}\prod_{j=1}^b f_j(t).
\]
Then the leading local coefficient is
\begin{equation}\label{eq:ex-figure-eight-leading}
 g_{4_1,2,\zeta}(t,0)=\frac{S_m(t)}{\Delta(t^m)}.
\end{equation}
More generally, with $x=Q-\zeta$, one has the uniform bound
\begin{equation}\label{eq:ex-figure-eight-pole-bound}
 \Delta(t^m)^{h+\lfloor h/2\rfloor+1}
      [x^h]g_{4_1,2,\zeta}(t,x)
       \in\mathbb Z[\zeta][t^{\pm1}]
       \qquad(h\ge0).
\end{equation}
\end{proposition}
\end{samepage}

\begin{proof}
For $w=\zeta t$, each $f_j$ factors as
$(w-\zeta^j)(1-w^{-1}\zeta^{-j})$, so a complete period has product
\begin{equation}\label{eq:ex-figure-eight-period}
 R(t):=\prod_{j=1}^m f_j(t)
      =(w^m-1)(1-w^{-m})=t^m+t^{-m}-2.
\end{equation}
At $Q=\zeta$ the factors of \eqref{eq:ex-figure-eight-series}
are periodic in $j$ with period $m$, so the summand with index $mk+b$
is $R(t)^k\prod_{j=1}^b f_j(t)$. In the simultaneous cyclotomic
neighborhood with $z=t^m-1$ one has $R=z^2/(1+z)$, so the geometric
series in $R$ converges there. Since $1-R=\Delta(t^m)$, this proves
\eqref{eq:ex-figure-eight-leading}.

For the higher orders we first put $Q=\zeta e^u$ and hold $y=Qt$
fixed. The complete block with index $k$ is then
\[
 G_k(y,u)=\prod_{r=1}^m
   \bigl(y+y^{-1}-\zeta^r e^{(mk+r)u}
                         -\zeta^{-r}e^{-(mk+r)u}\bigr).
\]
Its value at $u=0$ is $R_y=y^m+y^{-m}-2$, and its coefficient of
$u^a$ is a polynomial in $k$ of degree at most $a$. For odd $a$ the
degree is at most $a-1$, because the part of degree $a$ in $k$ is
obtained by replacing $(mk+r)u$ by $mku$ in every factor, and
\[
 \prod_{r=1}^m
   \bigl(y+y^{-1}-\zeta^r e^v-\zeta^{-r}e^{-v}\bigr)
       =y^m+y^{-m}-2\cosh(mv)
\]
has no odd coefficients in $v$. A final partial block of length $b<m$
has, in order $u^a$, degree at most $a$ in $k$.

Consider a term of total order $h$ in which $s$ complete blocks have
positive order, $o$ of them odd, and the final partial block has order
$a_0$. Before summation over the block indices its degree in the
indices is at most $h-o$, and summing over $s$ ordered indices raises
the degree by at most $s$. Since a block of even positive order has
order at least two,
\[
 2s-o\le h-a_0,\qquad
 h-o+s\le h+\lfloor h/2\rfloor.
\]
The resulting term has the form $R_y^{k-s}A(k)$, where
$\deg A\le h-o+s$ and $A(k)=0$ for $k=0,\ldots,s-1$, because no $s$
blocks can be chosen there. Expanding $A$ in falling factorials,
starting in degree $s$, and using
\[
 \sum_{k\ge s}(k)_{\underline r}R_y^{k-s}
       =\frac{r!R_y^{r-s}}{(1-R_y)^{r+1}}
       \qquad(r\ge s),
\]
we obtain no negative powers of $R_y$ and a denominator dividing
$\Delta(y^m)^{h+\lfloor h/2\rfloor+1}$.

Finally we substitute $y=\zeta e^ut$. Each color derivative raises
the pole order by at most one, so the contribution of the fixed-$y$
coefficient of order $j\le h$ to total order $h$ has pole order at
most
\[
 j+\lfloor j/2\rfloor+1+(h-j)
     \le h+\lfloor h/2\rfloor+1.
\]
The triangular change of coordinates $u=\log(1+x/\zeta)$ preserves
this bound. This constructs a rational local representative. Its joint
Taylor germ is the integral Newton germ and $\Delta(1)=1$, so
Lemma~\ref{lem:rel-integral-saturation} gives the integral
Laurent numerator, and uniqueness identifies the representative with
$g_{4_1,2,\zeta}$.
\end{proof}

Equation~\eqref{eq:ex-figure-eight-leading} is the
figure-eight case of the relation between Habiro's series and the ADO
invariants \cite{ADO} established in \cite[Theorems~2.1--2.2]{BH}. The
corresponding unified invariant is constructed in \cite{Willetts},
and the finite numerator is the case $n=2$ of
\cite[Theorem~9.7]{FangNewton}. The first three cases are
\begin{align}
 g_{4_1,2,1}(t,0)&=\Delta(t)^{-1},
 \label{eq:ex-figure-eight-m1}\\
 g_{4_1,2,-1}(t,0)&=\frac{\Delta(t)}{\Delta(t^2)},
 \label{eq:ex-figure-eight-m2}\\
 g_{4_1,2,\zeta_3}(t,0)
   &=\frac{\zeta_3^2t^2+\zeta_3t^{-2}
             +3\zeta_3t+3\zeta_3^2t^{-1}+5}{\Delta(t^3)},
 \label{eq:ex-figure-eight-m3}
\end{align}
where $\zeta_3^2+\zeta_3+1=0$.

The coefficient criterion of Section~\ref{subsec:rel-blocks} becomes
particularly explicit at $Q=-1$. Put $z=t^2-1$ and
$E(z)=1+z-z^2$, so that $\Delta(1+z)=E(z)/(1+z)$. Writing
\eqref{eq:ex-figure-eight-m2} in the basis $1,t-1$ gives
\begin{equation}\label{eq:ex-figure-eight-blocks}
 C_{-1,0,0}(z)=\frac{1+2z}{E(z)},\qquad
 C_{-1,1,0}(z)=-\frac{2+z}{E(z)}.
\end{equation}
Here $C_{-1,b,0}=\sum_{k\ge0}c_{2k+b}(-1)z^k$ are the raw blocks of
\eqref{eq:rel-raw-block-series}, which agree with the transported
blocks in quantum order zero. The initial values are
\[
\begin{array}{c|rrrr}
 k&0&1&2&3\\\hline
 c_{2k}(-1)&1&1&0&1\\
 c_{2k+1}(-1)&-2&1&-3&4
\end{array}
\]
and both rows satisfy the recurrence
\begin{equation}\label{eq:ex-figure-eight-block-recurrence}
 u_k+u_{k-1}-u_{k-2}=0\qquad(k\ge2).
\end{equation}
At the identity root the full sequence satisfies
\[
 \sum_{k\ge0}c_k(1)z^k=\frac{1+z}{E(z)}
       =1+z^2-z^3+2z^4-3z^5+5z^6+\cdots.
\]
The same Alexander recurrence polynomial thus governs the coefficients
at different roots, as predicted by
Corollary~\ref{cor:rel-block-recurrence}.

Differentiating the factors gives the first quantum coefficient. Put
$s=t+t^{-1}$ and $d=t-t^{-1}$. At $Q=-1$ the odd and even factors in
\eqref{eq:ex-figure-eight-series} are $2-s$ and $-2-s$,
respectively, and the $Q$-derivative of every factor is $d$.
Differentiating the geometrically summed series therefore gives
\begin{equation}\label{eq:ex-figure-eight-first-jet}
 g_{4_1,2,-1}(t,x)
   =\frac{3-s}{5-s^2}
      +x\,\frac{d(s-1)(s-5)}{(5-s^2)^2}+O(x^2),
 \qquad x=Q+1.
\end{equation}
Equivalently, the coefficient of $x$ is
$-t\partial_t\bigl(\Delta(t)/\Delta(t^2)\bigr)$, and it has double
poles along $\Delta(t^2)=0$.

The corresponding raw blocks are obtained from the first-order
transport operator. With $T=1+z$ and primes denoting $z$-derivatives,
this operator is
\[
 \mathcal L_{-1,1}
   =\begin{pmatrix}
       2z\partial_z^2&z\partial_z\\
       -\partial_z&2\partial_z+2z\partial_z^2
     \end{pmatrix},
\]
which is \eqref{eq:rel-first-block-operator} for $m=2$ after the
reduction $(t-1)^2=z-2(t-1)$. If $A=C_{-1,0,0}$ and $B=C_{-1,1,0}$,
the coefficient of $x$ in \eqref{eq:ex-figure-eight-first-jet}
has transported coordinates $-2TA'-B$ and $-B-2TB'$. Subtracting
$\mathcal L_{-1,1}(A,B)^{\mathsf T}$ yields
\begin{align}
 C_{-1,0,1}(z)
    &=\frac{z(4z^4-7z^3-17z^2-21z-4)}{E(z)^3},
 \label{eq:ex-figure-eight-raw-first-even}\\
 C_{-1,1,1}(z)
    &=\frac{-z^5-8z^4+19z^3+21z^2+30z-1}{E(z)^3}.
 \label{eq:ex-figure-eight-raw-first-odd}
\end{align}
The raw blocks of order one have denominator $E^3$, while their
combination under transport has denominator $E^2$; this is the double
pole along $\Delta(t^2)=0$ in \eqref{eq:ex-figure-eight-first-jet}.

\subsection{Alexander divisors and pole orders}
\label{subsec:ex-divisors}

The trace denominators of Section~\ref{sec:trace} produce
the Alexander divisors directly from a braid presentation. We use the
diagonally conjugate unreduced blocks
\[
 \widehat B_1=
 \begin{pmatrix}1-T^{-1}&T^{-1}&0\\1&0&0\\0&0&1\end{pmatrix},
 \qquad
 \widehat B_2=
 \begin{pmatrix}1&0&0\\0&1-T^{-1}&T^{-1}\\0&1&0\end{pmatrix};
\]
a diagonal conjugation does not change the principal minors, so it
preserves the determinants of the cut blocks. For the figure-eight
braid $\beta=(\sigma_1\sigma_2^{-1})^2$, direct multiplication gives
\[
 M_{4_1}(T)=
 ((\widehat B_1\widehat B_2^{-1})^2)^{\rm tr}
 =\begin{pmatrix}0&T^{-1}\\T(1-T)&(T-1)^2\end{pmatrix},
\]
whose regulated determinant is
\begin{equation}\label{eq:ex-figure-eight-burau}
 \det(I-zM_{4_1})
 =1-(T-1)^2z+(T-1)z^2,
 \qquad \det(I-M_{4_1})=T(3-T-T^{-1}).
\end{equation}
Here $z=s^{m\operatorname{ht}(\alpha)}$ for the root block of
$\alpha$. In $\mathfrak{sl}_3$ the three root blocks have
$T=t_1^m,t_2^m,(t_1t_2)^m$, so this single computation produces all
three factors of $D_{K,3}(\mathbf t^m)$, including the factor
$\Delta_K(t_1^mt_2^m)$ of the composite root.

For $T(3,4)$, represented by $(\sigma_1\sigma_2)^4$, the corresponding
cut block is
\[
 M_{3,4}(T)=
 \begin{pmatrix}
 (T-1)/T^2&(T-1)/T^3\\
 (T^2-T+1)/T^3&(T-1)/T^3
 \end{pmatrix},
\]
and
\[
 \det(I-M_{3,4})
 =T^{-3}\bigl(T^3-T^2+1-T^{-2}+T^{-3}\bigr).
\]
The monomial factors $T$ and $T^{-3}$ agree with
\eqref{eq:trace-exact-burau-alexander}, since the two braids have
$(b,w)=(3,0)$ and $(3,8)$. At $T=1$ the two cut matrices become
$\left(\begin{smallmatrix}0&1\\0&0\end{smallmatrix}\right)$ and
$\left(\begin{smallmatrix}0&0\\1&0\end{smallmatrix}\right)$, both
nilpotent. Their regulated determinants are therefore identically one
at that color, as used to remove the regulator in the proof of
Proposition~\ref{prop:trace-all-rank-trace-closure}.

For the figure-eight knot with $n=2$, the pole orders at $Q=1$ can be
determined in every quantum degree. Set $Q=e^{\hbar}$ and first hold
the dimension variable $y=Qt$ fixed in
\eqref{eq:ex-figure-eight-series}. With $a=y+y^{-1}-2$, each
factor of the $k$th summand is
\[
 y+y^{-1}-2\cosh(j\hbar)=a-j^2\hbar^2+O(\hbar^4).
\]
The series is therefore even in $\hbar$ at fixed $y$, and geometric
summation gives
\begin{equation}\label{eq:ex-figure-eight-dimension-jet}
 \mathcal J_{4_1,2}(e^{-\hbar}y;e^{\hbar})
   =\frac1{\Delta(y)}
     -\hbar^2\frac{y+y^{-1}-1}{\Delta(y)^4}
     +O(\hbar^4).
\end{equation}
Here the coefficient of $-\hbar^2$ has been summed as
\[
 \sum_{k\ge1}a^{k-1}\sum_{j=1}^k j^2
  =\frac1{1-a}\sum_{j\ge1}j^2a^{j-1}
  =\frac{1+a}{(1-a)^4},
\]
and $1-a=\Delta(y)$. All sums converge in the joint color--quantum
topology, so this is an identity of formal germs and of their rational
representatives. With $Q=q^2$, the coefficient of $(q-1)^2$ in
\eqref{eq:ex-figure-eight-dimension-jet} is the case $n=2$ of
\cite[equation~(9.8)]{FangNewton}, written there in the variable
$y+y^{-1}-2$.

To return to the Dynkin coordinate, put $\eta=Q-1$, $s=t+t^{-1}$ and
$d=t-t^{-1}$, so that $y=(1+\eta)t$ and $\hbar=\log(1+\eta)$. This
gives
\begin{equation}\label{eq:ex-figure-eight-identity-second-jet}
 \begin{split}
 g_{4_1,2,1}(t,\eta)
  ={}&\frac1{\Delta(t)}
       +\eta\frac{d}{\Delta(t)^2}\\
     &+\eta^2\left(
        \frac{t^{-1}}{\Delta(t)^2}
        +\frac{d^2}{\Delta(t)^3}
        -\frac{s-1}{\Delta(t)^4}\right)+O(\eta^3).
 \end{split}
\end{equation}
The first two terms in the parentheses come from expanding
$\Delta((1+\eta)t)^{-1}$, and the last one is the quantum correction
in \eqref{eq:ex-figure-eight-dimension-jet}. The coefficient of
$\eta^2$ has a pole of order four: after multiplication by
$\Delta(t)^4$ and reduction modulo $\Delta(t)$ the numerator is $-2$,
because $s=3$ on that divisor. The exponent four in
\eqref{eq:ex-figure-eight-pole-bound} is therefore attained in
quantum order two. At $t=1$ both coefficients of positive order in
\eqref{eq:ex-figure-eight-identity-second-jet} vanish, as they
must, since every colored value equals one at $Q=1$.

The same calculation determines the highest pole in every quantum
order.

\begin{samepage}
\begin{proposition}\label{prop:ex-figure-eight-exact-poles}
Put
\[
 G_N(t)=[\eta^N]g_{4_1,2,1}(t,\eta),\qquad
 \nu_N=N+\lfloor N/2\rfloor+1,\qquad
 a_g=(-1)^g\frac{(3g)!}{3^g g!}.
\]
Then $\Delta(t)^{\nu_N}G_N(t)\in\mathbb Z[t^{\pm1}]$, and
$G_N$ has a pole of order exactly $\nu_N$ at each zero of $\Delta$.
More precisely, the leading Laurent numerators satisfy
\begin{align}
 \Delta(t)^{3g+1}G_{2g}(t)&\equiv a_g\pmod{\Delta(t)},
 \label{eq:ex-figure-eight-even-top-pole}\\
 \Delta(t)^{3g+2}G_{2g+1}(t)
     &\equiv (3g+1)a_g(t-t^{-1})\pmod{\Delta(t)}.
 \label{eq:ex-figure-eight-odd-top-pole}
\end{align}
\end{proposition}
\end{samepage}

\begin{proof}
Write $H_{2g}(y)$ for the coefficient of $\hbar^{2g}$ at fixed $y$ in
the expansion \eqref{eq:ex-figure-eight-dimension-jet},
continued to all orders. In a summand with $k$ factors, choose $s\le g$
factors of positive order. The coefficient of $\hbar^{2r}$ in factor
$j$ is $-2j^{2r}/(2r)!$, so summing over the $s$ distinct indices
gives a polynomial in $k$ of degree at most $2g+s$, multiplied by
$a^{k-s}$. As in the proof of
Proposition~\ref{prop:ex-figure-eight-all-roots}, this
polynomial vanishes for $k<s$, so its falling-factorial expansion
starts in degree $s$, and geometric summation gives a denominator
dividing $(1-a)^{2g+s+1}$, hence $\Delta(y)^{3g+1}$.

Only $s=g$ can produce a pole of this order. In that case every
selected factor has order two, and the polynomial in $k$ is
$(-1)^g e_g(1^2,\ldots,k^2)$, where $e_g$ is the elementary symmetric
polynomial. Its leading term is
\[
 e_g(1^2,\ldots,k^2)
     =\frac{k^{3g}}{3^g g!}+O(k^{3g-1}).
\]
Indeed, by Newton's identities $e_g$ equals
$(\sum_{j\le k}j^2)^g/g!$ plus products of power sums with fewer than
$g$ factors, and such a product with $s<g$ factors has degree
$2g+s<3g$. The top falling factorial therefore sums to
\[
 \frac{(-1)^g}{3^g g!}
     \sum_{k\ge g}(k)_{\underline{3g}}a^{k-g}
   =a_g\frac{a^{2g}}{(1-a)^{3g+1}},
\]
and since $a\equiv1$ modulo $\Delta(y)=1-a$, the top numerator is
$a_g$.

Now we return to fixed $t$, keeping the logarithmic coordinate
$\hbar$ for the moment. With $\theta=t\partial_t$, the coefficient of
$\hbar^N$ is
\[
 \sum_{g=0}^{\lfloor N/2\rfloor}
   \frac{\theta^{N-2g}H_{2g}(t)}{(N-2g)!}.
\]
The term indexed by $g$ has pole order at most $N+g+1$, so the maximum
$\nu_N$ is attained only at $g=\lfloor N/2\rfloor$. Its highest term
is the one in \eqref{eq:ex-figure-eight-even-top-pole} or
\eqref{eq:ex-figure-eight-odd-top-pole}, because
$\theta\Delta=-(t-t^{-1})$. The sequence $\nu_N$ is strictly
increasing, so the change of coordinates from $\hbar$ to
$\eta=e^{\hbar}-1$, whose linear coefficient is one, does not change
the highest pole in any order. Lemma~\ref{lem:rel-integral-saturation}
gives the integral numerator. Finally, the zeros of $\Delta$ are
simple, $a_g\ne0$, and $(t-t^{-1})^2=5$ on $\Delta=0$, so the poles
occur at both zeros with exactly the stated orders.
\end{proof}

\subsection{Classical color coefficients}
\label{subsec:ex-classical}

To extend the quadratic formulas of
Corollary~\ref{cor:cons-alexander-hessian} to higher degree, put
\[
 y=T+T^{-1}-2,\qquad
 \Delta_K(T)=1+a y+b y^2+c y^3+O(y^4).
\]
Every symmetric Laurent polynomial with integer coefficients is an
integral polynomial in $T+T^{-1}$, so $a,b,c$ are integers; in
particular $a=\Delta_K''(1)/2$. Introduce logarithmic color
coordinates $t_i=e^{x_i}$ and the interval-root power sums
\[
 S_{2j}(\mathbf x)=\sum_{1\le p<q\le n}
             (x_p+\cdots+x_{q-1})^{2j}.
\]
With
\[
 \beta=b+\frac a{12}-\frac{a^2}{2},\qquad
 \gamma=c+\frac b6+\frac a{360}-ab-\frac{a^2}{12}
                                      +\frac{a^3}{3},
\]
the complete classical color expansion through degree six is
\begin{equation}\label{eq:ex-log-six}
 \begin{split}
 g_{K,n,1}(e^{\mathbf x},0)
 ={}&1-aS_2+\frac{a^2}{2}S_2^2-\beta S_4\\
 &-\frac{a^3}{6}S_2^3+a\beta S_2S_4-\gamma S_6
                                      +O((\mathbf x)^8).
 \end{split}
\end{equation}
To prove this, substitute $e^x+e^{-x}-2=x^2+x^4/12+x^6/360+O(x^8)$
in $\log\Delta_K$, which gives
\[
 \log\Delta_K(e^x)=ax^2+\beta x^4+\gamma x^6+O(x^8);
\]
summing over the positive roots and exponentiating the negative gives
\eqref{eq:ex-log-six}. By Alexander symmetry no odd
degrees occur in logarithmic coordinates; in the ordinary coordinates
$t_i-1$ they do.

For type $A_2$, set $t_1=1+u$ and $t_2=1+v$. Direct expansion gives
\begin{equation}\label{eq:ex-a2-four}
 \begin{split}
 g_{K,3,1}((1+u,1+v),0)={}&
 1-2a(u^2+uv+v^2)\\
 &+a(2u^3+u^2v+uv^2+2v^3)\\
 &+(3a^2-2a-2b)(u^4+v^4)\\
 &+(6a^2-a-4b)(u^3v+uv^3)\\
 &+(9a^2-a-6b)u^2v^2+O((u,v)^5);
 \end{split}
\end{equation}
for example, one may use
\[
 \Delta_K(1+w)^{-1}
 =1-aw^2+aw^3+(a^2-a-b)w^4+O(w^5)
\]
in the three factors with $w=u,v,u+v+uv$. At $Q=1$ the Newton basis is
the monomial basis, so \eqref{eq:ex-a2-four} computes every
$c_{r,s}(1)$ with $r+s\le4$. The following table records the
independent values; the coefficients with the two indices interchanged
are equal.
\begin{center}\small
\begin{tabular}{lrr|rrrrrr}
\toprule
$K$&$a$&$b$&$c_{2,0}=c_{1,1}$&$c_{3,0}$&$c_{2,1}$
 &$c_{4,0}$&$c_{3,1}$&$c_{2,2}$\\
\midrule
$3_1$&1&0&$-2$&2&1&1&5&8\\
$4_1$&$-1$&0&2&$-2$&$-1$&5&7&10\\
$T(2,5)$&3&1&$-6$&6&3&19&47&72\\
$T(2,7)$&6&5&$-12$&12&6&86&190&288\\
$T(3,4)$&5&5&$-10$&10&5&55&125&190\\
$3_1\#3_1$&2&1&$-4$&4&2&6&18&28\\
$3_1\#3_1\#3_1$&3&3&$-6$&6&3&15&39&60\\
$3_1\#4_1$&0&$-1$&0&0&0&2&4&6\\
\bottomrule
\end{tabular}
\end{center}
All coefficients in the table are evaluated at $Q=1$, and the
Alexander polynomials used are
$1+y$, $1-y$, $1+3y+y^2$, $1+6y+5y^2+y^3$,
$1+5y+5y^2+y^3$, $(1+y)^2$, $(1+y)^3$ and $1-y^2$,
respectively, with $y=T+T^{-1}-2$. The knots $T(2,5)$ and
$3_1\#3_1\#3_1$ have the same quadratic and cubic classical color
terms, and the quartic terms distinguish them through their different
values of $b$; formula~\eqref{eq:ex-log-six} gives the
corresponding expansions in every rank. For $3_1\#4_1$ one has $a=0$
and $b=-1$, so the quadratic and cubic color terms vanish, and the
first nonzero classical color term has degree four.

\subsection{Connected sums}
\label{subsec:ex-connected-sums}

Connected sums give examples whose entire classical color germs agree.
Let $T$ denote the trefoil in the convention of
\eqref{eq:ex-trefoil-fundamental}, and put $G=T\#T$ and
$S=T\#\overline T$. Both have normalized Alexander polynomial
$\Delta_T^2$, so for every $n\ge2$ their constant quantum coefficients
agree:
\[
 g_{G,n,1}(\mathbf t,0)=g_{S,n,1}(\mathbf t,0)
                        =D_{T,n}(\mathbf t)^{-2}.
\]
This equality includes all classical Newton coefficients, in every
degree, and it exhibits repeated Alexander poles.

With $f(Q)=Q^{-2}+Q^{-4}-Q^{-6}$,
Proposition~\ref{prop:cons-knot-operations} gives their fundamental
invariants in type $A_2$:
\begin{align}
 J_{G,\omega_1}^{(3)}(Q)
 &=f(Q)^2
   =Q^{-4}+2Q^{-6}-Q^{-8}-2Q^{-10}+Q^{-12},
 \label{eq:ex-granny}\\
 J_{S,\omega_1}^{(3)}(Q)
 &=f(Q)f(Q^{-1})=3-Q^4-Q^{-4}.
 \label{eq:ex-square}
\end{align}
With $h=Q-1$ these become
\[
 J_{G,\omega_1}^{(3)}=1-16h^2+64h^3-108h^4+O(h^5),\qquad
 J_{S,\omega_1}^{(3)}=1-16h^2+16h^3-36h^4+O(h^5),
\]
so their first Newton coefficients satisfy
\[
 c^G_{1,0}(Q)-c^S_{1,0}(Q)=48(Q-1)^2+O((Q-1)^3).
\]
The two Newton elements are separated by the term of order two in
$Q-1$ of this coefficient.

\subsection{The Kinoshita--Terasaka knot}
\label{subsec:ex-alexander-trivial}

The Kinoshita--Terasaka knot gives a nonconstant lift whose classical
color germ is constant. Let $K$ be the representative of $11n42$
specified by the planar diagram in \cite{KnotAtlasKT}. Its reduced
fundamental invariant is
\begin{equation}\label{eq:ex-kt-jones}
 \begin{split}
 J_{K,\omega_1}^{(2)}(Q)
 ={}&Q^{-6}-2Q^{-5}+2Q^{-4}-2Q^{-3}+Q^{-2}
          +2Q-2Q^2+2Q^3-Q^4\\
 ={}&1-Q^{-6}(Q-1)^2(Q^2+1)(Q^6-1),
 \qquad \Delta_K(T)=1.
 \end{split}
\end{equation}
The eleven crossing tuples of this diagram are
\[
 \begin{gathered}
 (4,2,5,1),\ (8,4,9,3),\ (12,5,13,6),\ (2,8,3,7),\\
 (9,18,10,19),\ (11,21,12,20),\ (6,13,7,14),\\
 (15,10,16,11),\ (17,22,18,1),\ (19,15,20,14),\ (21,16,22,17).
 \end{gathered}
\]
At a tuple $(a,b,c,d)$ the two bracket smoothings pair
$(a,b),(c,d)$ and $(a,d),(b,c)$. Their weights are $A$ and $A^{-1}$,
respectively, and each additional state circle contributes a factor
$-A^2-A^{-2}$. The sum over all $2^{11}$ states is
\[
 \langle K\rangle
 =A^{-19}-2A^{-15}+2A^{-11}-2A^{-7}-A^5
                     +2A^9-2A^{13}+2A^{17}-A^{21}.
\]
Orient the diagram by traversing the arc labels cyclically from
$1$ to $22$; its writhe is then $w=-1$. Multiplication by
$(-A^3)^{-w}$ followed by $A=Q^{-1/4}$ gives
\eqref{eq:ex-kt-jones}, with the bracket normalization of
\cite[Section~2, Theorems~2.6 and~2.8]{KauffmanState}; in this
convention the three-crossing trefoil has value
$Q^{-1}+Q^{-3}-Q^{-4}$.

For the Alexander polynomial, identify the two over-arcs $b,d$ at each
crossing, order the resulting Wirtinger generators by their least arc
label, and keep the displayed order of the crossings. At a crossing of
sign $\sigma$ we use the relation $y^\sigma xy^{-\sigma}z^{-1}$, where
$x,y,z$ are the incoming under-generator, the over-generator and the
outgoing under-generator. Its abelianized Fox row has entries
$T^\sigma$, $1-T^\sigma$ and $-1$ in these three columns. After each
row of a negative crossing is multiplied by $T$, the minor obtained by
deleting the last row and column is $-T^5$, so the normalized
Alexander polynomial is one.

Theorem~\ref{main:type-a} now gives, simultaneously for every rank,
\[
 g_{K,n,1}(\mathbf t,0)=1,\qquad
 c_{\boldsymbol\alpha}(1)=0\quad
                    (\boldsymbol\alpha\ne\mathbf0),
\]
and every local quantum coefficient is a Laurent polynomial in the
colors. Nevertheless the lift is not constant: its first Newton
coefficient for $n=2$ is
\begin{equation}\label{eq:ex-kt-newton}
 \begin{split}
 c_1(Q)&=\frac{J_{K,\omega_1}^{(2)}(Q)-1}{Q-1}
       =-Q^{-6}(Q-1)(Q^2+1)(Q^6-1)\\
       &=-12(Q-1)^2+30(Q-1)^3+O((Q-1)^4).
 \end{split}
\end{equation}

\subsection{Recurrence symbols along Alexander divisors}
\label{subsec:ex-recurrence}

The Alexander divisors also constrain the symbols of difference
equations. For $K=3_1$ or $4_1$, put $\tau=1$ or $3$, respectively,
and write
\[
 \delta_1=t_1^2-\tau t_1+1,\quad
 \delta_2=t_2^2-\tau t_2+1,\quad
 \delta_{12}=(t_1t_2)^2-\tau t_1t_2+1.
\]
Up to a Laurent unit, $D_{K,3}=\delta_1\delta_2\delta_{12}$.
For a recurrence satisfying the hypotheses of
Theorem~\ref{thm:cons-logarithmic-annihilator}, let $A_1,A_2$ be the
first symbol derivatives defined there. The theorem gives the explicit
necessary conditions
\begin{equation}\label{eq:ex-three-divisibilities}
 \delta_1\mid A_1,\qquad
 \delta_2\mid A_2,\qquad
 \delta_{12}\mid A_1+A_2
 \quad\text{in }\Q[t_1^{\pm1},t_2^{\pm1}].
\end{equation}
Indeed, each of the three quadratics is irreducible in the Laurent
ring and coprime to its nonzero logarithmic derivative. The first two
depend on one color only, whereas
$\theta_1\delta_{12}=\theta_2\delta_{12}$; the factorwise conclusion
$\delta\mid V(\delta)$ of the theorem therefore gives
\eqref{eq:ex-three-divisibilities}. Conversely, these
divisibilities imply that the vector field $A_1\theta_1+A_2\theta_2$
preserves $(D_{K,3})$.

For example, the symbol
\[
 p_0=\delta_1(L_1-1)+\delta_2(L_2-1)
\]
passes the two simple-root tests but fails the composite-root test:
setting $t_1=1$ turns $\delta_{12}$ into $\delta_2$ and
$\delta_1+\delta_2$ into $(2-\tau)+\delta_2$, which is not divisible by
$\delta_2$ because $2-\tau\ne0$. Hence no choice of the higher-order
terms in $Q-1$ of the coefficients turns this symbol into an
annihilator satisfying the hypotheses of the theorem.

For comparison, put
\[
 A=\delta_1\delta_2,\qquad
 B=\delta_2(t_1^2-1)-\delta_1(t_2^2-1).
\]
The vector field $V=A(\theta_1-\theta_2)$ annihilates $t_1t_2$, so it
kills the composite-root factor, and it satisfies
$V(D_{K,3})=BD_{K,3}$, because
$\theta_i\Delta_K(t_i)/\Delta_K(t_i)=(t_i^2-1)/\delta_i$. By the
computation in the proof of
Theorem~\ref{thm:cons-logarithmic-annihilator},
\begin{equation}\label{eq:ex-first-order-recurrence}
 \bigl(A(E_1-E_2)+(Q-1)B\bigr)g_{K,3,1}
       =O((Q-1)^2),
\end{equation}
that is, this operator annihilates the local component modulo
$(Q-1)^2$.

\begingroup
\renewcommand{\addcontentsline}[3]{}
\section*{Declaration of AI Use}
GPT-6.2 Astra was used to help check the mathematical arguments and calculations.
\par
\endgroup

\begingroup\raggedright

\endgroup
\end{document}